\documentclass[12pt]{colt2024} % Include author names

\usepackage{hyperref}

\SetCommentSty{mycommfont}

\usepackage{amsmath}
\usepackage{amssymb}
\usepackage{mathtools}
\usepackage{tikz}
\usepackage{comment}
\usepackage{enumitem}
\usepackage{booktabs}

\usetikzlibrary{arrows.meta,backgrounds,patterns,positioning,shapes.geometric,quotes}

\usepackage[capitalize,noabbrev]{cleveref}

\newtheorem{assumption}[theorem]{Assumption}
\newtheorem{property}[theorem]{Property}
\Crefname{assumption}{Assumption}{Assumptions}
\Crefname{assumption}{Property}{Properties}

\let\originalleft\left
\let\originalright\right
\renewcommand{\left}{\mathopen{}\mathclose\bgroup\originalleft}
\renewcommand{\right}{\aftergroup\egroup\originalright}
\newcommand{\norm}[1]{\left\lVert#1\right\rVert}
\newcommand{\abs}[1]{\left\lvert#1\right\rvert}

\newcommand{\decayfactor}{\rho}
\newcommand{\bP}[2][]{\Pr\ifthenelse{\isempty{#1}}{}{_{#1}}\left[#2\right]}
\newcommand{\bE}[2][]{\mathop\mathbb{E}\ifthenelse{\isempty{#1}}{}{_{#1}}\left[#2\right]}
\newcommand{\bI}[2][]{\mathop\mathbb{I}\ifthenelse{\isempty{#1}}{}{_{#1}}\left[#2\right]}
\newcommand{\Var}[2][]{\mathbf{Var}\ifthenelse{\isempty{#1}}{}{_{#1}}\left[#2\right]}

\newcommand{\dist}{\mathsf{dist}}

\newcommand{\ALG}{\pi}

\allowdisplaybreaks

\title[Online Policy Optimization in Unknown Nonlinear Systems]{Online Policy Optimization in Unknown Nonlinear Systems}
\usepackage{times}
\coltauthor{%
 \Name{Yiheng Lin} \Email{yihengl@caltech.edu}\\
 \addr California Institute of Technology
 \AND
 \Name{James A. Preiss} \Email{japreiss@caltech.edu}\\
 \addr California Institute of Technology
 \AND
 \Name{Fengze Xie} \Email{fxxie@caltech.edu}\\
 \addr California Institute of Technology
 \AND
 \Name{Emile Anand} \Email{emilea@andrew.cmu.edu}\\
 \addr Carnegie Mellon University
 \AND
 \Name{Soon-Jo Chung} \Email{sjchung@caltech.edu}\\
 \addr California Institute of Technology
 \AND
 \Name{Yisong Yue} \Email{yyue@caltech.edu}\\
 \addr California Institute of Technology
 \AND
 \Name{Adam Wierman} \Email{adamw@caltech.edu}\\
 \addr California Institute of Technology
}

\begin{document}

\maketitle

\begin{abstract}
We study online policy optimization in nonlinear time-varying dynamical systems where the true dynamical models are unknown to the controller. This problem is challenging because, unlike in linear systems, the controller cannot obtain globally accurate estimations of the ground-truth dynamics using local exploration. We propose a meta-framework that combines a general online policy optimization algorithm (\texttt{ALG}) with a general online estimator of the dynamical system's model parameters (\texttt{EST}). We show that if the hypothetical joint dynamics induced by \texttt{ALG} with \emph{known} parameters satisfies several desired properties, the joint dynamics under \emph{inexact} parameters from \texttt{EST} will be robust to errors. Importantly, the final policy regret only depends on \texttt{EST}'s predictions on the visited trajectory, which relaxes a bottleneck on identifying the true parameters globally. To demonstrate our framework, we develop a computationally efficient variant of Gradient-based Adaptive Policy Selection, called Memoryless GAPS (M-GAPS), and use it to instantiate \texttt{ALG}. Combining \mbox{M-GAPS} with online gradient descent to instantiate \texttt{EST} yields (to our knowledge) the first local regret bound for online policy optimization in nonlinear time-varying systems with unknown dynamics.
\end{abstract}

\section{Introduction}
\label{sec:intro}
We consider a class of discrete-time policy optimization problems with unknown time-varying nonlinear dynamics (called nonautonomous systems in nonlinear control theory~\cite{slotine1991applied}).
Our setting specifies a particular functional form of the dynamics and the parameterized policy class that is broad enough to capture applications from drone control to robotic manipulation \citep{dawson2023safe,o2022neural,KADIRKAMANATHAN1995245,shi_adaptive_2020}.
Our goal is to optimize the control policy online to minimize the total cost even if the online agent cannot obtain a globally accurate model of the true dynamics.

Online policy optimization and the broader field of learning-based control have received significant attention over the last several years due to their ability to leverage data and adapt to time-varying dynamical systems \citep{fazel2018global,lin2023online,arous2021online,mokhtari2016online,zhao2022non,zhou2023efficient,hazan2007adaptive,10.1145/1553374.1553425,gradu2020adaptive,baby2021optimal}. Online policy optimization faces two major challenges in practice. The first comes from the unknown dynamical model, which increases the difficulties of deciding the right directions for policy improvement. %\yisong{this point is very subtle for the COLT community, as one can also do direct policy learning without a model}\yiheng{We need to double-check this.} 
The second comes from the generality of dynamics/policy classes, which requires the algorithm to apply broadly to nonlinear time-varying dynamics and general policy classes. The early works in this field considered the linear-quadratic regulator (LQR) \citep{fazel2018global} and linear time-invariant dynamics with adversarial disturbances \citep{agarwal2019online,cohen2019learning,chen2021black,simchowitz2018learning,li2019online}, where the dynamics are known, linear, and time-invariant. Since then, much progress has been made to address challenges arising from considering either the unknown dynamical models or the general nonlinear dynamics \citep{muthirayan2022adaptive,minasyan2021online,yu2022online,dogan2021regret}, but addressing the two challenges simultaneously is still open.

One line of work focuses on online policy optimization under increasingly general classes of dynamical systems and policies, but under the assumption that the true dynamical models are known \citep{chen2023regret,zhou2023efficient,agarwal2019logarithmic}. For example, \cite{lin2023online} proposes an algorithm with provable regret guarantees that can be applied to nonlinear time-varying dynamical systems with general policy classes. However, assuming exact knowledge of the dynamical systems can be particularly restrictive in many applications when the system is nonlinear and time-varying. Even if the online agent has oracle access to a good dynamical model estimator, it is unclear whether the model estimation errors will accumulate in the policy update.

Another line of work about online policy optimization focuses on learning the unknown dynamical models but is generally restricted to linear systems with specific policy classes \citep{qu2021exploiting,minasyan2021online}. A common approach in the literature is random local exploration, where the controller either sets the control input to be a random perturbation or a randomly added perturbation to obtain sufficiently accurate estimations of the true dynamical model with high probability \citep{dean2018regret, pmlr-v151-lale22a}. However, in a more general nonlinear dynamical model, the online agent can no longer rely on small random perturbations to identify the true dynamical model well globally. This is because how the system responds to a small perturbation (approximately) depends on its linearization at the current state, and a nonlinear dynamical model can have different linearizations at different states.

The adaptive control literature also studies a similar problem \citep{ANNASWAMY202118, slotine1991applied, ioannou2012robust, 1657677}.  Typically, a set of linear parameters over a set of basis functions~\citep{shi_adaptive_2020, KADIRKAMANATHAN1995245,o2022neural} is dynamically adjusted online to compensate for the effect of time-varying disturbances, in order to stabilize the system and improve its trajectory tracking performance. \cite{yu2022online} proposes an adaptive stabilizing algorithm for unknown linear systems without any system identification. \cite{shi2021metaadaptive} and \cite{o2022neural} recently proposed a meta-online adaptive control algorithm to address the time-varying prediction errors; however, their methods do not optimize the gains in the control policy.
Another related work from control theory is online robust control, where the goal is to ensure stability \citep{yu2022online,li2023online} or staying in the safety sets \citep{ho2021online}, subject to model uncertainty.

The motivating insight that we take from adaptive control is that the controller does not need to learn the true dynamical model to stabilize a system. The controller only needs to focus on ``fitting'' the actual trajectory it visited rather than ``actively exploring'' with the purpose of identifying the true model parameter. This idea of ``lazy learning'' is shared by some works on online robust control \citep{pmlr-v144-boffi21a, ho2021online,yu2023online}, which maintains a set of possible dynamical models that are consistent with past observations. In this work, we are interested in developing a general approach for online policy optimization that can address the challenges of dealing with both unknown dynamical models and general nonlinear dynamics/policy classes.

\textbf{Contributions:} We make three main contributions. First, we develop a meta-framework that combines an online policy optimization algorithm (\texttt{ALG}) with an online parameter estimator (\texttt{EST}), where \texttt{ALG} focuses on optimizing the policy parameters while \texttt{EST} focuses on estimating the unknown component in the dynamics. We specify a set of properties that, if satisfied, implies that our meta-framework can mimic the behavior of applying \texttt{ALG} with known dynamical models up to an error that depends on the predictions of \texttt{EST}. 
This setup enables us to reason about how to use existing results in online policy optimization and online regression (for model learning) as subroutines.

Second, we provide a theoretical analysis of our meta-framework, establishing conditions under which we can derive regret guarantees.
We study the behavior of our meta-framework in two steps:
The first step (\Cref{sec:meta-online-policy-selection}) focuses on the behaviors of \texttt{ALG} and treats the predictions of \texttt{EST} as external inputs. We specify a set of properties that make the joint dynamics of applying \texttt{ALG} to the original system robust to the errors injected by using \texttt{EST} instead of the true dynamical models. The second step (\Cref{sec:online-param-est}) formulates the task of \texttt{EST} as an online regression problem, where the states visited by \texttt{ALG} are treated as adversarial inputs that can adapt to the history, i.e., the adversary is non-oblivious).
Compared to a standard online regression problem of minimizing the errors, \texttt{EST} faces the additional challenge of minimizing the errors of the model's partial derivatives with respect to the state.
We address this challenge by showing a reduction from the regret of predicting these partial derivatives to the regret of predicting the unknown component when the original dynamics contain a certain level of randomness.

Third, we provide a concrete instantiation of our meta-framework on matched-disturbance dynamics. 
For \texttt{ALG}, we develop Memoryless Gradient-based Adaptive Policy Selection (M-GAPS), which extends the GAPS algorithm for online policy optimization \citep{lin2023online} to utilize only $O(1)$ computational/memory complexity per time step and may be of independent interest.
For \texttt{EST}, we utilize standard online regression.
Combining these components, we obtain a bound on the local regret, an online analog of the stationary point conditions in nonconvex optimization.
To our knowledge, this is the first local regret bound for online policy selection in nonlinear time-varying systems with unknown dynamics.

\section{Problem Setting}
\label{sec:problem}
We consider online policy optimization in a discrete-time dynamical system
that varies over time
with dynamics
$x_{t+1} = g_t \big(x_t, u_t, f_t(x_t, a_t^*) \big) + w_t$,
where $x_t \in \mathbb{R}^n$ denotes the system state, $u_t \in \mathbb{R}^m$ denotes the control input, and $g_t$ is the dynamical function. Here, $f_t(x_t, a_t^*) \in \mathbb{R}^k$ is a nonlinear residual term of which the online agent can make (noisy) observations. It has a known function form $f_t$ and an unknown parameter $a_t^* \in \mathcal{A} \subseteq \mathbb{R}^p$.
The disturbance term $w_t \in \mathcal{W} \subseteq \mathbb{R}^n$ does not depend on the states or the control inputs.

To control this system, the online agent adopts a time-varying control policy $\pi_t$ that is parameterized by a policy parameter $\theta_t \in \Theta \subseteq \mathbb{R}^d$, where $\Theta$ is a closed convex subset of $\mathbb{R}^d$. Specifically, the online agent picks the control input from the policy class
$u_t = \pi_t \big(x_t, \theta_t, f_t(x_t, \hat{a}_t))$.
Here, function $f_t(\cdot, \hat{a}_t)$ reflects the online agent's current estimation of the ground true nonlinear residual function $f_t(\cdot, a_t^*)$ at time step $t$. Intuitively, we assume the policy class $\pi_t$ cares about predicting the true nonlinear residual $f_t(x_t, a_t^*)$ rather than the unknown model parameter $a_t^*$.
The objective of the online agent is to minimize the total cost $\sum_{t=0}^{T-1} c_t$ incurred over a finite horizon, where the stage cost at time step $t$ is given by
$c_t = h_t(x_t, u_t, \theta_t)$.

We provide a simple nonlinear control example that can be captured by our online policy optimization framework to help the readers understand the concepts we discussed.
\begin{example}\label{example:nonlinear-control}
Consider the problem of controlling a scalar discrete-time nonlinear system:
\begin{align}\label{equ:nonlinear-dynamics-example}
    x_{t+1} = x_t + \Delta \left(u_t + f_t(x_t, a_t^*) + w_t\right), \text{ where }f_t(x_t, a_t^*) = \phi(x_t) \cdot a_t^*.
\end{align}
In this equation, $\Delta $ is the discretization step size. The nonlinear residual takes the form $\phi(x_t) \cdot a_t^*$, where $\phi: \mathbb{R} \to \mathbb{R}^k$ is a (nonlinear) feature map and $a_t^*$ is the unknown model parameter. To control this system, the online agent with an estimated model parameter $\hat{a}_t$ can adopt the policy class:
\begin{align}\label{equ:policy-class-example}
    u_t = - k_t x_t - f_t(x_t, \hat{a}_t), \text{ where } f_t(x_t, \hat{a}_t) = \phi(x_t) \cdot \hat{a}_t, \text{ and } k_t = \theta_t.
\end{align}
Here, the goal of the second term $- f_t(x_t, \hat{a}_t)$ is to cancel out the true nonlinear residual $f_t(x_t, a_t^*)$. In an ideal case where the online agent has access to the true model parameter $a_t^*$, policy \eqref{equ:policy-class-example} achieves the effect of removing the nonlinear residual and directly doing feedback control, resulting in the closed-loop dynamics
\(x_{t+1} = x_t + \Delta \left(- k_t x_t + w_t\right).\)
In this case, the problem reduces to finding the optimal policy parameters (gains) $\{\theta_t\}$ in a known time-varying dynamical system.
\end{example}

\subsection{Performance Metrics}\label{sec:metrics}

In the literature of online optimization, \textit{regret} is a common performance metric that directly compares the total cost $\sum_{t=0}^{T-1} c_t$ incurred by the online policy optimization algorithm against the optimal total cost one can achieve in hindsight.
Before introducing variants of the regret we study, we first introduce the concept of the \textit{surrogate cost}.
The present formulation extends the similarly named concept of \cite{lin2023online} to the setting where the true dynamical models are unknown.
\begin{definition}[Surrogate Cost]\label{def:surrogate-cost}
The surrogate cost function is $F_t(\theta) \coloneqq h_t(\tilde{x}^*_t(\theta), \tilde{u}^*_t(\theta), \theta)$, where $\tilde{x}^*_t(\theta)$ and $\tilde{u}^*_t(\theta)$ are the system state and control input at time $t$ if the agent applies the control input $\tilde{u}^*_\tau(\theta) \coloneqq \pi_\tau\left(\tilde{x}^*_\tau(\theta), \theta, f_\tau(x_\tau, a_\tau^*)\right)$ at time step $\tau = 0, 1, \ldots, t$.
\end{definition}
Intuitively, the surrogate cost $F_t(\theta)$ evaluates how good a policy parameter $\theta$ is at an intermediate time step $t$ by eliminating the interference of inexact estimations $\hat{a}_{0:t}$ of model parameters and any other policy parameters $\theta_{0:t-1}$ in the history that may be different with $\theta$. This concept is useful for defining different regret metrics. For example, the \textit{static regret} is a widely-used metric in the literature of online policy optimization (\cite{Cesa-Bianchi_Lugosi_2006,10.1145/1553374.1553425}) that compares the total cost of an online agent with the best static policy parameter in hindsight can be written as
$R^S(T) \coloneqq \sum_{t=0}^{T-1} c_t - \min_{\theta^* \in \Theta} \sum_{t=0}^{T-1} F_t(\theta^*).$

However, as noted in previous works \citep{lin2023online,hazan2017efficient}, regret metrics that directly compare the cost difference (like $R^S(T)$) are not always suitable for nonconvex cost functions, because gradient-based online optimization algorithms may easily get stuck in local
minima even when the cost functions are time-invariant. 
Therefore, the metric of \textit{local regret} is used. For a policy parameter sequence $\theta_{0:T-1}$, {\color{black} the local regret is defined as $R_\eta^L(T) \coloneqq \sum_{t=0}^{T-1} \norm{\nabla_{\eta, \Theta} F_t(\theta_t)}^2$. Here, the projected gradient $\nabla_{\eta, \Theta} F_t(\theta_t)$ (parameterized by $\eta$) is a surrogate of the original gradient $\nabla F_t(\theta_t)$ that also considers the constraint set $\Theta \subseteq \mathbb{R}^d$. Specifically, an update step with the projected gradient is equivalent to projecting the output of the original gradient descent step back onto $\Theta$, i.e., for any $\theta \in \Theta$, $\theta - \eta \nabla_{\eta, \Theta} F_t(\theta_t) = \Pi_\Theta(\theta - \eta \nabla F_t(\theta_t))$ \footnote{$\Pi_\Theta$ is the Euclidean projection to $\Theta$.}. This notion of local regret is first introduced by \citet{hazan2017efficient}, and we provide a formal definition in Definition \ref{def:proj-grad} in Appendix \ref{appendix:notation-and-defs}.} Intuitively, the local regret measures how well the policy parameter sequence $\theta_{0:T-1}$ tracks the changing stationary points of the surrogate cost functions $F_{0:T-1}$.
It is also helpful to understand local regret by drawing connections to the time-invariant setting, where a sublinear local regret implies convergence to a stationary point.

Although local regret is useful for measuring the performance of an online policy optimization algorithm under nonconvex surrogate costs, a limitation of applying it alone to our setting with unknown dynamical models is that the surrogate cost $F_t$ is defined in terms of \texttt{ALG}'s behavior with known true dynamics.
To address this limitation, in addition to bounding the local regret of the policy parameters $\theta_{0:T-1}$, we also bound the distance between the actual trajectory of the online agent and the trajectory it would achieve with the same policy parameters $\theta_{0:T-1}$ and exact knowledge of true model parameters $a_{0:T-1}^*$.

\section{Main Results}
\label{sec:main}
Our approach is outlined in \Cref{alg:meta-framework}, where {\color{black} two modules} \texttt{ALG} and \texttt{EST} work together to update the policy and estimated model parameter at each time step {\color{black} (see Figure \ref{fig:meta-framework} for an illustration). \texttt{ALG} and \texttt{EST} are responsible for optimizing the policy parameters $\theta_{0:T-1}$ and learning the unknown model parameters $a_{0:T-1}^*$ respectively:
\begin{itemize}[nolistsep, leftmargin=*]
\item \textbf{\texttt{ALG}:} At time step $t$, \texttt{ALG} receives the current state $x_t$, policy parameter $\theta_t$, and the known part of the time-varying system $\pi_t, g_t, h_t, f_t$. It also receives the current estimation $\hat{a}_t$ of the unknown model parameter $a_t^*$. Then, \texttt{ALG} outputs the new policy parameter $\theta_{t+1}$. Note that we allow \texttt{ALG} to leverage/memorize historical inputs by maintaining an internal state $y_t$.
\item \textbf{\texttt{EST}:} At time step $t$, \texttt{EST} receives the current state $x_t$ and a (noisy) observation $\tilde{f}_t$ of the unknown component $f_t(x_t, a_t^*)$. Then, \texttt{EST} outputs the new estimation $\hat{a}_{t+1}$. Like \texttt{ALG}, we allow \texttt{EST} to keep internal state/memory (e.g., to memorize historical input data). We require \texttt{EST} to minimize the \textit{trajectory-dependent model mismatches}:
\begin{subequations}\label{equ:online-param-est-obj}
\begin{align}
    &\text{Zeroth-order model mismatch: } \varepsilon_t(x_t, \hat{a}_t, a_t^*) \coloneqq \norm{f_t(x_t, \hat{a}_t) - f_t(x_t, a_t^*)},\label{equ:online-param-est-obj:e0}\\
    &\text{First-order model mismatch: }\varepsilon_t'(x_t, \hat{a}_t, a_t^*) \coloneqq \norm{\nabla_x f_t(x_t, \hat{a}_t) - \nabla_x f_t(x_t, a_t^*)}_F.\label{equ:online-param-est-obj:e1}
\end{align}
\end{subequations}
We adopt the shorthand $\varepsilon_t = \varepsilon_t(x_t, \hat{a}_t, a_t^*)$ and $\varepsilon_t' = \varepsilon_t'(x_t, \hat{a}_t, a_t^*)$ when the context is clear.
\end{itemize}
}

\begin{algorithm2e}
    \caption{Meta-Framework}\label{alg:meta-framework}
    %\LinesNumbered
    \DontPrintSemicolon
    {\bfseries Require:} \texttt{ALG} and \texttt{EST}\\
    {\bfseries Require:} Knowing functions $\{\pi_t, g_t, h_t, f_t\}$ at each time step $t$\\
    % , sequence of policies $\pi_{0:T-1}$ and residual functions $f_{0:T-1}$.
    {\bfseries Initialize:} State $x_0$; Policy parameter $\theta_0$; Model parameter estimation $\hat{a}_0$.\\
    \For{$t = 0, 1, \ldots, T-1$}{

    Decide control input $u_t = \pi_t(x_t, \theta_t, f_t(x_t, \hat{a}_t))$.

    Incur stage cost $h_t(x_t, u_t, \theta_t)$.

    $\theta_{t+1} \gets \texttt{ALG}.\text{update}(x_t, \theta_t, \pi_t, g_t, h_t, f_t, \hat{a}_t)$. \tcc*{\texttt{ALG} can have memory.}

    System evolves to $x_{t+1} = g_t(x_t, u_t, f_t(x_t, a_t^*)) + w_t.$

    Receive a (noisy) observation $\tilde{f}_t$ of $f_t(x_t, a_t^*)$.

    $\hat{a}_{t+1} \gets \texttt{EST}.\text{update}(x_t, \tilde{f}_t)$. \tcc*{\texttt{EST} can have memory.}
    }
\end{algorithm2e}

The key idea in analyzing our meta-framework (\Cref{alg:meta-framework}) is to characterize how the inexact model estimations generated by \texttt{EST} affect the behavior \texttt{ALG}.
We start by considering the ``ideal'' dynamics of applying \texttt{ALG} with exact model parameters $a_{0:T-1}^*$, which we denote as $\texttt{ALG}^*$, and compare them with the actual dynamics of \texttt{ALG} that performs the update with estimated model parameters $\hat{a}_{0:T-1}$. {\color{black} We state the key insight of our analysis in the informal lemma below, which connects the performance of the meta-framework with $\texttt{ALG}^*$ and the model mismatches.}

{\color{black}
\begin{lemma}[Informal]\label{lemma:informal-insight}
Suppose $\texttt{ALG}^*$ satisfies the desired properties in \Cref{sec:meta-online-policy-selection}. Then, the meta-framework (\Cref{alg:meta-framework}) generates the same policy parameters as $\texttt{ALG}^*$ with perturbation $\zeta_t$ on the update of $\theta_{t+1}$ (see Figure \ref{fig:exact-model}). Further, $\sum_{t=0}^{T-1} \norm{\zeta_t} = O\left(\sum_{t=0}^{T-1}\varepsilon_t + \sum_{t=0}^{T-1}\varepsilon_t'\right)$.
\end{lemma}
The formal statement of Lemma \ref{lemma:informal-insight} can be found in Theorem \ref{thm:meta-dynamics-regret-and-stability}.
}

{\color{black} The rest of this section is organized as following: In \Cref{sec:meta-online-policy-selection}, we specify the properties of $\texttt{ALG}^*$ that enables the meta-framework to be robust against inexact model parameters in the policy parameter update.
Then, in \Cref{sec:online-param-est}, we formulate \texttt{EST}'s task of learning $f_t(x_t, a_t^*)$ as an online optimization problem, where we view the state $x_t$ as picked by an adaptive adversary. We also discuss how this problem reduces to existing results on online optimization.}

\usetikzlibrary{calc}

\begin{figure}
    \begin{minipage}[c]{0.46\linewidth}
    \centering
    \begin{tikzpicture}
\draw[blue, thick] (0,3.6) rectangle (1.2,4.4);
\node at (0.6, 4.0) {Policy};

\coordinate (A1) at (0.6, 3.2);
\filldraw[black] (A1) circle (2pt);
\node at ($(A1) - (0, 0.4)$) {$u_t$};

\coordinate (A2) at (7, 3.2);
\filldraw[black] (A2) circle (2pt);
\node at ($(A2) - (0, 0.4)$) {$\hat{f}_t$};

\draw[blue, thick] (1.6,1.6) rectangle (2.8,2.4);
\node at (2.2, 2.0) {Env};
\coordinate (D3) at (2.2, 2.8);
%\filldraw[blue] (D3) circle (2pt);

\draw[blue, thick] (3.2,1.6) rectangle (4.4,2.4);
\node at (3.8, 2.0) {ALG};

\draw[blue, thick] (4.8,1.6) rectangle (6,2.4);
\node at (5.4, 2.0) {EST};

\coordinate (A3) at (2.2, 6.0);
\filldraw[black] (A3) circle (2pt);
\node at ($(A3) + (0, 0.4)$) {$x_t$};
\coordinate (A4) at (3.8, 6.0);
\filldraw[black] (A4) circle (2pt);
\node at ($(A4) + (0, 0.4)$) {$\theta_t$};
\coordinate (A5) at (5.4, 6.0);
\filldraw[black] (A5) circle (2pt);
\node at ($(A5) + (0, 0.4)$) {$\hat{a}_t$};

\coordinate (B3) at (2.2, 0.4);
\filldraw[black] (B3) circle (2pt);
\node at ($(B3) - (0, 0.4)$) {$x_{t+1}$};
\coordinate (B4) at (3.8, 0.4);
\filldraw[black] (B4) circle (2pt);
\node at ($(B4) - (0, 0.4)$) {$\theta_{t+1}$};
\coordinate (B5) at (5.4, 0.4);
\filldraw[black] (B5) circle (2pt);
\node at ($(B5) - (0, 0.4)$) {$\hat{a}_{t+1}$};

\coordinate (C) at (0.6, 4.4);

\draw[-stealth,thick] (A3) -- (2.2, 2.4);
\draw[-stealth,thick] (A3) -- (3.8, 2.4);
\draw[-stealth,thick] (A3) -- (C);
\draw[-stealth,thick] (A4) -- (3.8, 2.4);
\draw[-stealth,thick] (A4) -- (C);
\draw[-stealth,thick] (A5) -- (3.8, 2.4);
\draw[-stealth,thick] (A5) -- (C);
\draw[-stealth,thick] (A5) -- (5.4, 2.4);

\draw[-stealth,thick] (0.6, 3.6) -- (A1);
\draw[-stealth,thick] (A1) -- (2.2, 2.4);
\draw[-stealth,thick] (A2) -- (5.4, 2.4);

\draw[-stealth,thick] (2.2, 1.6) -- (B3);
\draw[-stealth,thick] (3.8, 1.6) -- (B4);
\draw[-stealth,thick] (5.4, 1.6) -- (B5);
\end{tikzpicture}
    \caption{The meta-framework.}
    \label{fig:meta-framework}
    \end{minipage}
    \hfill
    \begin{minipage}[c]{0.46\linewidth}
    \centering
    %\usetikzlibrary{calc}

\begin{tikzpicture}
\draw[blue, thick] (0,3.6) rectangle (1.2,4.4);
\node at (0.6, 4.0) {Policy};

\coordinate (A1) at (0.6, 3.2);
\filldraw[black] (A1) circle (2pt);
\node at ($(A1) - (0, 0.4)$) {$u_t$};

\draw[blue, thick] (1.6,1.6) rectangle (2.8,2.4);
\node at (2.2, 2.0) {Env};
\coordinate (D3) at (2.2, 2.8);
%\filldraw[blue] (D3) circle (2pt);

\draw[blue, thick] (3.2,1.6) rectangle (4.4,2.4);
\node at (3.8, 2.0) {ALG};

\coordinate (A3) at (2.2, 6.0);
\filldraw[black] (A3) circle (2pt);
\node at ($(A3) + (0, 0.4)$) {$x_t$};
\coordinate (A4) at (3.8, 6.0);
\filldraw[black] (A4) circle (2pt);
\node at ($(A4) + (0, 0.4)$) {$\theta_t$};
\coordinate (A5) at (5.4, 6.0);
\filldraw[black] (A5) circle (2pt);
\node at ($(A5) + (0, 0.4)$) {$a_t^*$};

\coordinate (B3) at (2.2, 0.4);
\filldraw[black] (B3) circle (2pt);
\node at ($(B3) - (0, 0.4)$) {$x_{t+1}$};
\coordinate (B4) at (3.8, 0.4);
\filldraw[black] (B4) circle (2pt);
\node at ($(B4) - (0, 0.4)$) {$\theta_{t+1}$};
\coordinate (B5) at (5.4, 0.4);

\coordinate (C) at (0.6, 4.4);

\draw[-stealth,thick] (A3) -- (2.2, 2.4);
\draw[-stealth,thick] (A3) -- (3.8, 2.4);
\draw[-stealth,thick] (A3) -- (C);
\draw[-stealth,thick] (A4) -- (3.8, 2.4);
\draw[-stealth,thick] (A4) -- (C);
\draw[-stealth,thick] (A5) -- (3.8, 2.4);
\draw[-stealth,thick] (A5) -- (C);

\draw[-stealth,thick] (0.6, 3.6) -- (A1);
\draw[-stealth,thick] (A1) -- (2.2, 2.4);

\filldraw[color=black, fill=white, thick] (3.8, 1) circle (0.25);

\draw[-stealth,thick] (2.2, 1.6) -- (B3);
\draw[-stealth,thick] (3.8, 1.6) -- (B4);
\draw[thick] (3.55, 1) -- (4.05, 1);

\draw[-stealth,thick] (5, 1) -- (4.05, 1);
\filldraw[black] (5, 1) circle (2pt);
\node at (5.4, 1) {$\zeta_t$};

\end{tikzpicture}
    \caption{$\texttt{ALG}^*$ with perturbations $\zeta_{0:T-1}$ on policy parameter updates.}
    \label{fig:exact-model}
    \end{minipage}
\end{figure}

\subsection{Online Policy Optimization}\label{sec:meta-online-policy-selection}
In this section, we take a perspective that views the updates performed by \texttt{ALG} as part of a joint dynamics formed together with the original dynamical system.
Compared to the common approach of analyzing \texttt{ALG} separately from the dynamical system to which it applies, our dynamical view enables us to compare the differences of applying \texttt{ALG} under different external inputs (i.e. different $\hat a_t$ estimates) more efficiently. % by studying the properties of the joint dynamics.

We consider the class of online policy optimization algorithms whose joint dynamics with the original system can be written in the following form: {\color{black} When the model parameter $a_t$ is given as the input to \texttt{ALG} at time step $t$, the joint dynamics can be written as}
\begin{align}\label{equ:meta-dynamics-general-pred}
    \begin{pmatrix}
        x_{t+1}\\
        y_{t+1}\\
        \theta_{t+1}
    \end{pmatrix} = q_t(x_t, y_t, \theta_t, a_t) = \begin{pmatrix}
        q_t^x(x_t, y_t, \theta_t, a_t)\\
        q_t^y(x_t, y_t, \theta_t, a_t)\\
        q_t^\theta(x_t, y_t, \theta_t, a_t)
    \end{pmatrix}, \text{ for } x_t \in \mathbb{R}^n, y_t \in \mathbb{R}^p, \theta_t \in \Theta \subset \mathbb{R}^d.
\end{align}
Here, $y_t \in \mathbb{R}^p$ is an auxiliary state that \texttt{ALG} can use to store something besides the system state $x_t$ and the policy parameter $\theta_t$ to help it perform the update. For example, $y_t$ can be a finite memory buffer that stores information from the past. It can also be the integral of past states in an integral controller. Thus, we introduce $y_t$ to allow broader classes of online policy optimization algorithms, and we will provide a concrete example of $y_t$ later in \Cref{sec:application:MGAPS}.

{\color{black} The goal of formulating joint dynamics \eqref{equ:meta-dynamics-general-pred} is to compare the behaviors of the meta-framework and $\texttt{ALG}^*$ with perturbations on policy parameter updates. Specifically, recall that $\hat{a}_{0:T-1}$ denote the estimated model parameters of \texttt{EST}. The actual trajectory of the meta-framework is
\begin{align}\label{equ:meta-dynamics-inexact-pred}
    \text{Meta-framework: } (x_{t+1}, y_{t+1}, \theta_{t+1})^\top = q_t(x_t, y_t, \theta_t, \hat{a}_t).
\end{align}
We compare it with the joint dynamics of $\texttt{ALG}^*$ (see Figure \ref{fig:exact-model}). Recall that $\texttt{ALG}^*$ denotes the scenario when \texttt{ALG} has access to exact model parameters $a_{0:T-1}^*$:
\begin{align}\label{equ:meta-dynamics-exact-pred}
    \texttt{ALG}^* \text{ with perturbations: } (x_{t+1}, y_{t+1}, \theta_{t+1})^\top = q_t(x_t, y_t, \theta_t, a_t^*) + (0, 0, \zeta_t)^\top.
\end{align}
Here, $\zeta_t$ is an additive perturbation on the update equation of policy parameter $\theta_{t+1}$. To understand \eqref{equ:meta-dynamics-exact-pred} intuitively, it is helpful to draw connections with the process of using a gradient-based optimizer to update the parameter $\theta_t$ in ML, where $\zeta_t \equiv 0$ corresponds to the case when exact gradients are available. In contrast, nonzero perturbations correspond to the more practical case when the optimizer can only use biased estimations of the gradient, which still performs well in general.
}

Note that the estimated model parameters $\hat{a}_{0:T-1}$ generated by \texttt{EST} may also depend on the state $x_t$ and other parts of the dynamical system. Thus, a natural question is whether we should also incorporate the update rule of \texttt{EST} into the joint dynamical system in \eqref{equ:meta-dynamics-inexact-pred}, where we include $\hat{a}_t$ as another element of the joint state. However, we still choose to model $\hat{a}_t$ as an external input in \eqref{equ:meta-dynamics-inexact-pred} and handle the update of $\hat{a}_t$ separately in \Cref{sec:online-param-est}. This is because our approach requires comparing the actual joint dynamics with \eqref{equ:meta-dynamics-exact-pred}. Since $a_t^*$ is an external input decided by the environment in \eqref{equ:meta-dynamics-exact-pred}, keeping the joint state space identical in \eqref{equ:meta-dynamics-inexact-pred} makes the comparison easier. Further, a strength of our proof framework based on the joint dynamics is that we can show the actual trajectory \eqref{equ:meta-dynamics-inexact-pred} will stay close to \eqref{equ:meta-dynamics-exact-pred}. However, we know that the estimated model parameter sequence $\{\hat{a}_t\}$ will not converge to the true sequence $\{a_t^*\}$ in general.

We state three important properties of the joint dynamics induced by \texttt{ALG}. The first property is about the Lipschitzness with respect to the model mismatches $\varepsilon_t$ and $\varepsilon_t'$.

\begin{property}\label{assump:Lipschitzness-a}[Lipschitzness]
For any $x_t, y_t, \theta_t, \hat{a}_t$ that satisfies
$\norm{x_t} \leq R_x, \norm{y_t} \leq R_y, \theta_t \in \Theta, \hat{a}_t \in \mathcal{A},$
the following Lipschitzness conditions hold:
\begin{align*}
    \norm{q_t^x(x_t, y_t, \theta_t, a_t^*) - q_t^x(x_t, y_t, \theta_t, \hat{a}_t)} &\leq \alpha_x \varepsilon_t(x_t, \hat{a}_t, a_t^*) + \beta_x \varepsilon_t'(x_t, \hat{a}_t, a_t^*),\\
    \norm{q_t^y(x_t, y_t, \theta_t, a_t^*) - q_t^y(x_t, y_t, \theta_t, \hat{a}_t)} &\leq \alpha_y \varepsilon_t(x_t, \hat{a}_t, a_t^*) + \beta_y \varepsilon_t'(x_t, \hat{a}_t, a_t^*),\\
    \norm{q_t^\theta(x_t, y_t, \theta_t, a_t^*) - q_t^\theta(x_t, y_t, \theta_t, \hat{a}_t)} &\leq \alpha_\theta \varepsilon_t(x_t, \hat{a}_t, a_t^*) + \beta_\theta \varepsilon_t'(x_t, \hat{a}_t, a_t^*).
\end{align*}
Further, $q_t^\theta(x, y, \theta, a_t^*)$ is $(L_{\theta, x}, L_{\theta, y})$-Lipschitz in $(x, y)$.
\end{property}

Intuitively, Property \ref{assump:Lipschitzness-a} says that the error brought by the inexact model parameters only ``distort'' the ideal joint dynamics \eqref{equ:meta-dynamics-exact-pred} in the form of zeroth-order and first-order prediction errors. Therefore, to bound the error injected into the joint dynamics at every step, \texttt{EST} only needs to minimize $\varepsilon_t$ and $\varepsilon_t'$ on the actual state trajectory $x_{0:T-1}$ that the online agent visits.
Note that this property can be viewed as a standard assumption about Lipschitzness if \texttt{ALG} is a gradient-based algorithm. This is because all terms that involve the unknown model parameter will take the form $f_t(x_t, \hat{a}_t)$ and $\nabla_x f_t(x_t, \hat{a}_t)$ in the joint dynamics.

The second property is about contraction stability of $x_t$ and $y_t$ under exact model parameters $a_{0:T-1}^*$. As we show in \Cref{thm:meta-dynamics-regret-and-stability}, this property guarantees that the dynamical updates of states $x_t$ and $y_t$ in the joint dynamics are robust to the model mismatches $\{\varepsilon_t, \varepsilon_t'\}_{0:T-1}$.

\begin{property}\label{assump:contraction-and-stability}[Contraction Stability]
For any sequence $\theta_{0:T-1}$ that satisfies the slowly time-varying constraint that $\norm{\theta_t - \theta_{t-1}} \leq \epsilon_\theta$ for all time step $t$, the partial dynamical system
\begin{align}\label{equ:partial-dynamics}
    x_{t+1} = q_t^x(x_t, y_t, \theta_t, a_t^*),
    \quad \quad
    y_{t+1} = q_t^y(x_t, y_t, \theta_t, a_t^*)
\end{align}
satisfies that
$\norm{x_t} \leq R_x^* < R_x$
and
$\norm{y_t} \leq R_y^* < R_y$ always hold
if the system starts from $(x_\tau, y_\tau) = (0, 0)$.
Further, for some function $\gamma: \mathbb{Z}_{\geq 0} \to \mathbb{R}_{\geq 0}$ that satisfies $\sum_{t=0}^\infty \gamma(t) \leq C$, from any initial states $(x_\tau, y_\tau), (x_\tau', y_\tau')$ that satisfy $\norm{x_\tau}, \norm{x_\tau'} \leq R_x$ and $\norm{y_\tau}, \norm{y_\tau'} \leq R_y$, the trajectory satisfies
$\norm{(x_{\tau+t}, y_{\tau+t}) - (x_{\tau+t}', y_{\tau+t}')} \leq \gamma(t) \cdot \norm{(x_\tau, y_\tau) - (x_\tau', y_\tau')}.$
\end{property}
Note that Property \ref{assump:contraction-and-stability} is different with the contraction assumption of \citet{lin2023online} because it also considers the internal state $y_t$ of \texttt{ALG} besides the system state $x_t$. The requirement that $\sum_{t=0}^\infty \gamma(t) \leq C$ is also weaker than the exponential decay rate in \citet{lin2023online}.

Intuitively, Property \ref{assump:contraction-and-stability} guarantees that when the exact model parameters $\{a_t^*\}$ are replaced by inexact $\{\hat{a}_t\}$, the resulting trajectory $\{(x_t, y_t)\}$ still stays close to the trajectory that $\theta_{0:T-1}$ would achieve under exact predictions once the mismatch errors $\varepsilon_t, \varepsilon_t'$ are small or bounded.
Property \ref{assump:contraction-and-stability} can be viewed as an extension of the time-varying stability and contractive perturbation property in \cite{lin2023online} to include state $y_t$ maintained by \texttt{ALG}. This is required in our framework because $y_t$ can be affected by the prediction errors and it is involved in the dynamics of updating $\theta_t$.

The third property we need is the robustness of the update rule of the policy parameter $\theta_t$. Specifically, it requires the regret guarantee achieved by \texttt{ALG} to be robust against a certain level of adversarial disturbances $\{\zeta_t\}$ on the update dynamics of $\theta_t$.
\begin{property}\label{assump:robustness}[Robustness]
Consider the joint dynamics in \eqref{equ:meta-dynamics-exact-pred}.
When $\norm{\zeta_t} \leq \bar{\zeta}$ holds for all $t$, the resulting $\{\theta_t\}$ satisfies the slowly-time-varying constraint $\norm{\theta_t - \theta_{t-1}} \leq \epsilon_\theta$ for all time $t$. Further, $\texttt{ALG}^*$ with perturbations \eqref{equ:meta-dynamics-exact-pred} can achieve a regret guarantee $R(T, \sum_{t=0}^{T-1} \norm{\zeta_t})$ that depends on the total magnitude of the perturbation sequence $\zeta_{0:T-1}$.
\end{property}
To understand Property \ref{assump:robustness}, we can think about online gradient descent (OGD) in online optimization problems without state or dynamics. It is known that this approach is robust to (biased) disturbances on the gradient estimation, and the total amount of added disturbances will affect the final regret bound (see, for example, \Cref{thm:nonconvex-biased-OGD}).

Now, we present our main results about the stability of applying \texttt{ALG} with inexact model parameters and the regret bound in \Cref{thm:meta-dynamics-regret-and-stability}. Besides, \Cref{thm:meta-dynamics-regret-and-stability} also bounds the distances between the actual trajectory and the trajectory achieved by applying the same policy parameter sequence with the exact model parameter sequence.

\begin{theorem}\label{thm:meta-dynamics-regret-and-stability}
Suppose Properties \ref{assump:Lipschitzness-a}, \ref{assump:contraction-and-stability}, and \ref{assump:robustness} hold. Let $\xi = \{x_t, y_t, \theta_t\}_{0:T-1}$ be the trajectory of the meta-framework (\Cref{alg:meta-framework}). If the prediction errors $\{\varepsilon_t, \varepsilon_t'\}_{0:T-1}$ are uniformly bounded such that the following inequalities hold for all time step $t$: $\alpha_\theta \varepsilon_t + \beta_\theta \varepsilon_t' \leq \bar{\zeta} / 2$, and
\begin{align*}
    (\alpha_x + \alpha_y) \varepsilon_t + (\beta_x + \beta_y)\varepsilon_t' \leq{}& \min\left\{\frac{\sqrt{2}\bar{\zeta}}{4(L_{\theta, x} + L_{\theta, y})C}, \frac{\min\{R_x - R_x^*, R_y - R_y^*\}}{C}\right\},
\end{align*}
%\soonjo{did you clearly define what tilde means? $\tilde{x}$?}
then the trajectory $\xi$ satisfies $\norm{x_t} \leq R_x$, $\norm{y_t} \leq R_y$, and $\norm{\theta_t - \theta_{t-1}} \leq \epsilon_\theta$ for all time steps~$t$. Further, define $\tilde{\xi} \coloneqq \{\tilde{x}_t, \tilde{y}_t, \theta_t\}_{0:T-1}$, where $\{\tilde{x}_t, \tilde{y}_t\}_{0:T-1}$ are obtained by implementing the policy parameters $\theta_{0:T-1}$ with exact model parameters $a_{0:T-1}^*$, i.e., the trajectory of partial joint dynamics \eqref{equ:partial-dynamics}. The trajectory $\tilde{\xi}$ achieves the regret $R(T, \sum_{t=0}^{T-1}\norm{\zeta_t})$ with $\sum_{t=0}^{T-1} \norm{\zeta_t}$ upper bounded by
\begin{align*}
    \left(\alpha_\theta + \sqrt{2} C (L_{\theta, x} + L_{\theta, y})(\alpha_x + \alpha_y)\right) \sum_{t=0}^{T-1}\varepsilon_t + \left(\beta_\theta + \sqrt{2} C (L_{\theta, x} + L_{\theta, y})(\beta_x + \beta_y)\right) \sum_{t=0}^{T-1}\varepsilon_t'.
\end{align*}
The total distances between the states on the trajectories $\xi$ and $\tilde{\xi}$ satisfies that
\begin{align*}
    \sum_{t=1}^{T} \norm{(x_{t}, y_{t}) - (\tilde{x}_{t}, \tilde{y}_{t})} \leq C \left((\alpha_x + \alpha_y)\sum_{t=0}^{T-1}\varepsilon_t + (\beta_x + \beta_y)\sum_{t=0}^{T-1}\varepsilon_t'\right).
\end{align*}
\end{theorem}

We defer the proof of \Cref{thm:meta-dynamics-regret-and-stability} to \Cref{appendix:meta-dynamics-regret-and-stability}. Intuitively, \Cref{thm:meta-dynamics-regret-and-stability} states that when the prediction error terms $\{\varepsilon_t, \varepsilon_t'\}_{0:T-1}$ are uniformly bounded, the actual trajectory $\xi$ of applying \texttt{ALG} with inexact model parameters $\hat{a}_{0:T-1}$ will be uniformly bounded. Further, if the actual parameter sequence of $\theta_{0:T-1}$ is applied with exact model parameters $a_{0:T-1}^*$, the resulting trajectory $\tilde{\xi}$ achieves a regret guarantee that depends on the magnitudes of the prediction errors. It is worth noticing that the regret in \Cref{thm:meta-dynamics-regret-and-stability} can be any regret that depends on the trajectory $\tilde{\xi}$. And as we discussed in \Cref{sec:metrics}, we evaluate the regret on trajectory $\tilde{\xi}$ rather than $\xi$ because the metrics like the local regret are designed for evaluating the actual policy parameters $\theta_{0:T-1}$ rather than the whole trajectory $\xi$. The distances between $\xi$ and $\tilde{\xi}$ are bounded in the last inequality in \Cref{thm:meta-dynamics-regret-and-stability}.

To show \Cref{thm:meta-dynamics-regret-and-stability}, the key idea is to fit the trajectory $\tilde{\xi}$ into the dynamical equation \eqref{equ:meta-dynamics-exact-pred}, where we design $\zeta_t$ to compensate the difference between the update rules $q_t(\tilde{x}_t, \tilde{y}_t, \theta_t, a_t^*)$ and $q_t(x_t, y_t, \theta_t, \hat{a}_t)$. To leverage Property \ref{assump:robustness}, we show the perturbations $\zeta_{0:T-1}$ we constructed are uniformly bounded by $\bar{\zeta}$. We bound $\zeta_t$ and the distances between $\xi$ and $\tilde{\xi}$ by induction. The induction is important because the magnitude of $\zeta_t$ depends on the distance between $\{x_t, y_t\}_{0:T-1}$ and $\{\tilde{x}_t, \tilde{y}_t\}_{0:T-1}$ in the past time steps. On the other hand, to bound the distance between $\{x_t, y_t\}$ and $\{\tilde{x}_t, \tilde{y}_t\}$, we need to leverage the contraction property in Property \ref{assump:contractive-and-stability}, which relies on $\norm{\zeta_t} \leq \bar{\zeta}$ so that $\theta_{0:T-1}$ is slowly time-varying. Lastly, we conclude the proof with the bounds on the distance between $\xi$ and $\tilde{\xi}$ as well as the norm of $\zeta_t$ that depend on the model mismatches $\{\varepsilon_t, \varepsilon_t'\}_{0:T-1}$.

\subsection{Online Parameter Estimation}\label{sec:online-param-est}
The second part of our meta framework focuses on predicting the unknown model parameter based on possibly noisy observations of the true nonlinear residual $f_t(x_t, a_t^*)$. A critical difference with prior works on system identification or model-based learning (e.g., \cite{dean2020sample}) is that we only seek to optimize the zeroth-order and first-order model mismatches $\{\varepsilon_t, \varepsilon_t'\}$ (defined in \eqref{equ:online-param-est-obj}) on the actual trajectory that the online agent experiences. It is worth noticing that, although learning the ground-truth model parameter $a_t^*$ is impossible for a general nonlinear residual, minimizing the sum of zeroth-order model mismatches incurred on the actual trajectory can be formulated as a classic online regression problem, which we discuss below:

\noindent\textbf{Online regression problem:} At the beginning, the environment commits a sequence of error functions $e_t: \mathbb{R}^n \times \mathcal{A} \to \mathbb{R}, t = 0, \ldots, T-1$, which are defined as
$e_t(x, a) \coloneqq f_t(x, a_t^*) - f_t(x, a)$
for
$t = 0, \ldots, T-1$. \footnote{Thus, the error functions $e_{0:T-1}$ will not adapt to the inputs and online decisions.}
The relationship between the error function $e_t$ and the model mismatches $\{\varepsilon_t, \varepsilon_t'\}$ is
$\varepsilon_t = \norm{e_t(x_t, \hat{a}_t)},$ and  $\varepsilon_t' = \norm{\nabla_x e_t(x_t, \hat{a}_t)}.$
At each time step $t$, the online parameter estimator \texttt{EST} predicts
$\hat{a}_t = \texttt{EST}(x_{0:t-1}, \hat{a}_{0:t-1}) \in \mathcal{A}$, which means the estimation $\hat{a}_t$ can be a general function of the historical states and estimations. Then, the environment reveals $x_t \in B_n(0, R_x)$ that can depend on the history $x_{0:t-1}$ and $\hat{a}_{0:t-1}$. We define the stage loss of \texttt{EST} as $\ell_t = \norm{e_t(x_t, \hat{a}_t)}^2$, which is equal to the squared $\ell_2$-norm of the model mismatch $e_t(x_t, \hat{a}_t)$.

Under different sets of assumptions on the error functions and the sequence of true model parameters $\{a_t^*\}$, existing online algorithms can achieve regret guarantees. We consider a general form of expected regret bound:
$\mathbb{E}\left[\sum_{t=1}^T \ell_t\right] \leq R_0^\ell(T),$
where the expectation is taken over the randomness of implementing \texttt{EST} and generating $x_t$.
While different assumptions and designs of \texttt{EST} can achieve different bounds on $R_0^\ell(T)$, an example we provide in \Cref{sec:application} shows that a simple gradient estimator can achieve sublinear $R_0^\ell(T)$ when the nonlinear residual can be decomposed as $f_t(x, a) = \phi(x)\cdot a$ under the path length constraint $\sum_{t=1}^{T-1} \norm{a_{t+1}^* - a_t^*} \leq C$ (see \Cref{sec:application:grad-est} for detailed discussions).

While most prior works focus on minimizing the magnitude of the zeroth-order model mismatch $e_t(x_t, \hat{a}_t)$, we also need to bound the first-order model mismatch $\nabla_x e_t(x_t, \hat{a}_t)$ because it contributes to the regret bound in \Cref{thm:meta-dynamics-regret-and-stability} (recall that $\norm{\nabla_x e_t(x_t, \hat{a}_t)}_F = \varepsilon_t'$).
Our main result in this section is about an automatic reduction from the regret bound $R_0^\ell(T)$ to a bound on the expected sum of the squared gradients $\mathbb{E}\left[\sum_{t=1}^T \norm{\nabla_x e_t(x_t, \hat{a}_t)}_F^2\right].$

\begin{remark}
Besides the online policy optimization problem for control, the regret bound that concerns $\norm{\nabla_x e_t(x_t, \hat{a}_t)}$ can be of independent interest for the problem of online regression, because it characterizes how sensitive the regression loss is to any perturbations on the input sequence $x_{0:T-1}$ under the same estimations $\hat{a}_{0:T-1}$. Intuitively, if gradients of the error functions are small, the estimations $\hat{a}_{0:T-1}$ will be robust to small perturbations on the input sequence.
\end{remark}

To enable a reduction from the regret bound $R_0^\ell(T)$ to the gradient error bound, we employ Property \ref{assump:perturb-the-input} about the dynamical system that generates the state $x_t$.
Specifically, we require there to be at least a small level of randomness when choosing $x_t$. Recall that $\hat{a}_{t+1}$ is decided based on the history $x_{0:t}$ and $\hat{a}_{0:t}$. We define the filtrations $\mathcal{F}_t \coloneqq \sigma\left(x_{1:t}, \hat{a}_{1:t}\right)$ and $\mathcal{F}_t' \coloneqq \sigma\left(x_{1:t}, \hat{a}_{1:t+1}\right)$, which satisfy $\mathcal{F}_t \subseteq \mathcal{F}_t' \subseteq \mathcal{F}_{t+1}$.

\begin{property}\label{assump:perturb-the-input}
There is a certain level of random disturbances when generating each state $x_t$, i.e., for some $\bar{\epsilon} > 0$ and $\underline{\sigma} > 0$, one can find a $\sigma$-algebra $\mathcal{G}_t$ such that $\mathcal{F}_t' \subseteq \mathcal{G}_t \subseteq \mathcal{F}_{t+1}$ and
\(
    \norm{x_{t+1} - \mathbb{E}[x_{t+1}\mid \mathcal{G}_t]} \leq \overline{\epsilon}, \ \text{Cov}(x_{t+1} \mid \mathcal{G}_t)\succeq \underline{\sigma} I.
\)
\end{property}

Intuitively, the randomness enforced by Property \ref{assump:perturb-the-input} will ``force'' \texttt{EST} to also minimize the gradient of the error functions. To see this, suppose an input state $x_{t}$ is given by $\bar{x}_{t} + v_{t}$, where $\bar{x}_t$ is the mean and $v_t$ is a random disturbance. When the disturbance $v_t$ is sufficiently small, we know that $e_t(x_t, \hat{a}_t) \approx e_t(\bar{x}_t, \hat{a}_t) + \nabla_x e_t(\bar{x}_t, \hat{a}_t) \cdot v_t$ by Taylor's expansion. Since we can pick $v_t$ randomly in different directions, we know the zeroth-order loss $\mathbb{E}[e_t(x_t, \hat{a}_t)^2]$ cannot converge to zero unless the magnitude of the gradient $\nabla_x e_t(\bar{x}_t, \hat{a}_t)$ converges to zero. We follow this intuition to show the reduction from the regret bound $R_0^\ell(T)$ to the total gradient error in \Cref{thm:gradient-error-bound}.

\begin{theorem}\label{thm:gradient-error-bound}
Suppose that for all time $t$, each dimension $i \in [k]$ of the error function satisfies 
\[\norm{\nabla_x e_t(x, a)_i} \leq \beta_e, \text{ and } \norm{\nabla_x^2 e_t(x, a)_i} \leq \gamma_e,
\text{ for any } x \in B(0, R_x) \text{ and } a \in \mathcal{A}.\] Suppose Property \ref{assump:perturb-the-input} holds with $\bar{\epsilon} \leq \min\{\frac{1}{4}, \frac{1}{2\gamma_e}, \frac{1}{4\beta_e \gamma_e}\}$ and $\underline{\sigma} > 0$. If \texttt{EST} achieves the zeroth-order regret $\mathbb{E}\left[\sum_{t=1}^T \ell_t\right] \leq R_0^\ell(T) \leq \overline{\epsilon}^3 T$, the expected total squared gradient loss satisfies that
\begin{align*}
    \mathbb{E}\left[\sum_{t=1}^T \norm{\nabla_x e_t(x_t, \hat{a}_t)}_F^2\right] \leq \frac{2k}{\underline{\sigma}}(1 + \gamma_e + \beta_e \gamma_e) \bar{\epsilon}^3 T + 2 k \gamma_e^2 \overline{\epsilon}^2 T.
\end{align*}
\end{theorem}
Recall that $k$ is the dimension of the unknown component $f_t(x_t, a_t^*)$. We defer the proof of \Cref{thm:gradient-error-bound} to Appendix \ref{appendix:gradient-error-bounds}. We provide the following corollary to help the readers understand this result in a special case when $R_0^\ell(T)$ is $O(\sqrt{T})$. For example, a gradient estimator can achieve this regret bound if $\ell_t$ is convex in $a$ (see \Cref{sec:application:grad-est}).

\begin{corollary}
Under the same assumptions as \Cref{thm:gradient-error-bound}, if \texttt{EST} achieves $R_0^\ell(T) = O(\sqrt{T})$ and Property \ref{assump:perturb-the-input} holds with $\bar{\epsilon} = \theta(T^{1/6})$ and $\underline{\sigma} = \Omega(\bar{\epsilon}^2)$, then the expected total squared gradient loss is bounded by $\mathbb{E}\left[\sum_{t=1}^T \norm{\nabla_x e_t(x_t, \hat{a}_t)}_F^2\right] = O(k T^{5/6})$.
\end{corollary}

In summary, with the help of \Cref{thm:gradient-error-bound}, we reduce the problem of bounding the total squared first-order prediction errors $\sum_{t=0}^{T-1} (\varepsilon_t')^2$ to the standard online optimization problem. By substituting the bounds on $\varepsilon_{0:T-1}$ and $\varepsilon_{0:T-1}'$ into \Cref{thm:meta-dynamics-regret-and-stability} in \Cref{sec:meta-online-policy-selection}, one can derive the local regret bound for the actual joint dynamics and bound the distance between trajectories $\xi$ and $\tilde{\xi}$.

\section{Application: Matched Disturbance}
\label{sec:application}
In this section, we consider an instantiation of our setting to demonstrate the effectiveness of our meta-framework. Specifically, we study the matched-disturbance dynamics \citep{8788577,9867457,6225272}, where the controller can choose a control input to ``cancel out'' the nonlinear residual term $f_t(x_t, a_t^*)$ when the exact model parameter $a_t^*$ is available. The dynamics have the form
\begin{align}\label{equ:dynamics-fully-actuated}
    x_{t+1} = g_t(x_t, u_t, f_t(x_t, a_t^*)) + w_t = \phi_t(x_t, u_t + f_t(x_t, a_t^*)) + w_t.
\end{align}
A ubiquitous application of the matched disturbance dynamics is the general joint-space dynamics of robotic manipulators
\citep{siciliano-robotics} when the system has actuators for every joint.
The matched-disturbance structure also appears in tilted-rotor rotorcraft \citep{rajappa-tilthex,9144371}, which can move in six degrees of freedom.
In both cases, due to the second-order structure of the rigid-body dynamics, all external disturbances are equivalent to additional joint torque (resp. rotor tilt/thrust) inputs. Our \Cref{example:nonlinear-control} also fits into the framework of \eqref{equ:dynamics-fully-actuated}.
To control a matched-disturbance system, a natural policy class is to first cancel out the nonlinear residual with $- f_t(x_t, \hat{a}_t)$ and then apply an actuation term $\psi_t(x_t, \theta_t)$ to achieve the optimal costs. This policy class can be expressed as
\begin{align}\label{equ:policy-cancel-out}
    u_t = \pi_t(x_t, \theta_t, f_t(x_t, \hat{a}_t)) = - f_t(x_t, \hat{a}_t) + \psi_t(x_t, \theta_t).
\end{align}
{\color{black} To derive local regret for the meta-framework, we need assumptions (Assumptions \ref{assump:Lipschitz-and-smoothness}-\ref{assump:estimator-environment}) on the system that includes the dynamics, policy classes, and costs, which we discuss in detail in Appendix \ref{appendix:application-matched-disturbance-main}. Note that the matched-disturbance dynamics/policy class we consider can recover the setting \citep{lin2023online} as a special case when $f_t$ and $w_t$ are always zero (so there is no need to estimate $a_t^*$). We recover the same regret bound as \cite{lin2023online} in that special case (see Lemma \ref{lemma:recover-exact-regret}).}

\subsection{Online Policy Optimization: M-GAPS}\label{sec:application:MGAPS}

This section introduces a general online policy optimization algorithm, Memoryless Gradient-based Adaptive Policy Selection (M-GAPS, \Cref{alg:M-GAPS}), which can serves as \texttt{ALG} in our meta-framework. %Since M-GAPS is a general algorithm with applications not limited to disturbance-matched dynamics, we provide the pseudo-code of M-GAPS in the general form in \Cref{alg:M-GAPS}.
M-GAPS use $\hat{a}_t$ to estimate how the current state $x_t$ and policy parameter $\theta_t$ would affect the next state $x_{t+1}$ and the current cost. The estimations are characterized by
\begin{subequations}\label{equ:g_hat_h_hat}
\begin{align}
    \hat{g}_{t+1\mid t}(x_t, \theta_t) &\coloneqq g_t(x_t, \pi_t(x_t, \theta_t, f_t(x_t, \hat{a}_t)), f_t(x_t, \hat{a}_t)), \text{ and }\\
    \hat{h}_{t\mid t}(x_t, \theta_t) &\coloneqq h_t\left(x_t, \ALG_t\left(x_t, \theta_{t}, f_t(x_t, \hat{a}_t)\right), \theta_t\right)
\end{align}
\end{subequations}
Although M-GAPS can be applied to any online policy optimization problems that fit into the setting we discussed in \Cref{sec:problem}, we focus on its application to disturbance-matched dynamics and policy class for theoretical analysis. We verify that the joint dynamics of M-GAPS satisfy the required properties of our meta-framework to derive a concrete regret bound in Appendix \ref{appendix:application-details}.

\begin{algorithm2e}
    \caption{Memoryless Gradient-based Adaptive Policy Selection (M-GAPS, for \texttt{ALG})}\label{alg:M-GAPS}
    %\LinesNumbered
    \DontPrintSemicolon
    {\bfseries Parameters:} Learning rate $\eta$.
    
    {\bfseries Initialize:} Internal state $y_0 = \mathbf{0}$.
    
    \For{$t = 0, 1, \ldots, T-1$}{
    Take inputs $x_t, \theta_t, g_t, \pi_t, h_t, f_t$ and $\hat{a}_t$. \tcc*{Given when meta-framework calls $\texttt{ALG}.\text{update}$.}

    Use $\hat{a}_t$ to obtain $\hat{g}_{t+1\mid t}$ and $\hat{h}_{t\mid t}$. \tcc*{$\hat{g}_{t+1\mid t}$ and $\hat{h}_{t\mid t}$ are defined in \eqref{equ:g_hat_h_hat}.}
    
    Update $y_{t+1} \gets \left.\frac{\partial \hat{g}_{t+1\mid t}}{\partial x_t}\right|_{x_t, \theta_t} \cdot y_t + \left.\frac{\partial \hat{g}_{t+1\mid t}}{\partial \theta_t}\right|_{x_t, \theta_t}.$ \tcc*{Update partial derivatives accumulator.}
    
    Let $G_t \gets \left.\frac{\partial \hat{h}_{t\mid t}}{\partial x_t}\right|_{x_t, \theta_t} \cdot y_t + \left.\frac{\partial \hat{h}_{t\mid t}}{\partial \theta_t}\right|_{x_t, \theta_t}.$
    
    Output $\theta_{t+1} \gets \prod_\Theta\left(\theta_t - \eta G_t\right)$. \tcc*{$\Pi_\Theta$ is the Euclidean projection to $\Theta$.}
    }
\end{algorithm2e}

The design of M-GAPS takes inspiration from the Gradient-based Adaptive Policy Selection (GAPS) algorithm \citep{lin2023online}, but it significantly improves computational efficiency and generality. Specifically, in the setting where the online agent has exact knowledge of the time-varying dynamical models, M-GAPS can achieve the same regret guarantees as GAPS, while the memory/computational complexities are reduced from $O(\log T)$ to $O(1)$. To understand why this improvement is possible, note that the core problem of GAPS is to design an efficient estimation of $\nabla F_t(\theta_t)$. GAPS does two steps of approximations: first replacing the imaginary trajectory with the actual trajectory and then doing a bounded-memory truncation. In comparison, M-GAPS only keeps the first approximation step of GAPS but greatly simplifies the computation by introducing the auxiliary internal state $y_t$ that accumulates past partial derivatives. Intuitively, the estimation of M-GAPS is even closer to $\nabla F_t(\theta_t)$ than GAPS, so it can achieve the same regret guarantees as GAPS. A more detailed comparison between M-GAPS and GAPS can be found in Appendix \ref{appendix:useful-lemmas}.

{\color{black}
A key step of our proof shows that, when exact model parameters $a_{0:T-1}^*$ are available, M-GAPS is robust against perturbations on policy parameter updates as required by Property \ref{assump:robustness} in \Cref{sec:meta-online-policy-selection}.

\begin{lemma}\label{lemma:recover-exact-regret}
Under Assumptions \ref{assump:Lipschitz-and-smoothness} and \ref{assump:contractive-and-stability}, Property \ref{assump:robustness} holds when $\eta \leq \bar{\eta}$ for some positive constant $\bar{\eta}$ and
$R_\eta^L(T, \sum_{t=0}^{T-1} \norm{\zeta_t}) = O\left(\frac{1}{\eta} (1 + V_{\text{sys}} + V_w) + \eta T + \eta^3 T + \frac{1}{\eta}\sum_{t=1}^{T-1} \norm{\zeta_t}\right),$ where the variation intensities are defined as $V_{w} = \sum_{t=1}^{T-1} \norm{w_t - w_{t-1}}$ and
\begin{align*}
    V_{\text{sys}} ={}& \sum_{t = 1}^{T-1}
    \Big(\sup_{x \in \mathcal{X}, u \in \mathcal{U}}
        \norm{\phi_t(x, u) - \phi_{t-1}(x, u)}
    + \sup_{x \in \mathcal{X}, \theta \in \Theta}
        \norm{\psi_t(x, \theta) - \psi_{t-1}(x, \theta)}\\
    &+ \sup_{x \in \mathcal{X}, u \in \mathcal{U}, \theta \in \Theta}
        \abs{h_t(x, u, \theta) - h_{t-1}(x, u, \theta)}\Big).
\end{align*}
\end{lemma}
The formal statement and proof of Lemma \ref{lemma:recover-exact-regret} can be found in Appendix \ref{appendix:fully-actuated-assumptions}. Note that in the special case of \cite{lin2023online}, we have $\Theta = \mathbb{R}^d$, $V_w = 0$, and $\sum_{t=1}^{T-1} \norm{\zeta_t} = 0$. The local regret bound $R_\eta^L(T, 0)$ of M-GAPS given by Lemma \ref{lemma:recover-exact-regret} matches the local regret bound of GAPS in \cite{lin2023online}, because the projected gradients are identical with the gradients when $\Theta = \mathbb{R}^d$.
}

\subsection{Online Parameter Estimation: Gradient Estimator}\label{sec:application:grad-est}

In the application of matched-disturbance dynamics, we assume the online parameter estimator \texttt{EST} can make a noisy observation $\tilde{f}_t$ of the true nonlinear residual $f_t(x_t, a_t^*)$ after it decides $\hat{a}_t$ at each time step $t$. Recall that the prediction error function is defined as $e_t(x, a) \coloneqq f_t(x, a) - f_t(x, a_t^*)$ and the true prediction loss at time step $t$ as $\ell_t \coloneqq \norm{e_t(x_t, \hat{a}_t)}^2$. We instantiate \texttt{EST} with the gradient estimator (\Cref{alg:grad-est}), where $\tilde{f}_t$ is a (noisy) observation of $f_t(x_t, a_t^*)$ provided by the environment. It performs online gradient descent on an estimated prediction loss function constructed from $\tilde{f}_t$.

\begin{algorithm2e}
    \caption{Gradient Estimator (for \texttt{EST})}\label{alg:grad-est}
    %\LinesNumbered
    \DontPrintSemicolon
    {\bfseries Parameters:} Learning rate $\iota$; {\bfseries Initialize:} Model parameter estimation $\hat{a}_0$.
    
    \For{$t = 0, 1, \ldots, T-1$}{
    Take inputs $x_t$ and $\tilde{f}_t$. \tcc*{Inputs given when meta-framework calls $\texttt{EST}.\text{update}$.}

    Incur loss $\tilde{\ell}_t(x_t, \hat{a}_t, \tilde{f}_t) \coloneqq \|f_t(x_t, \hat{a}_t) - \tilde{f}_t\|^2$.
    
    Update and output $\hat{a}_{t+1} \gets \prod_{\mathcal{A}}\left(\hat{a}_t - \iota \cdot \left. \partial \ell_t / \partial a_t \right|_{x_t, \hat{a}_t, \tilde{f}_t}\right)$.}
\end{algorithm2e}

Using Theorems \ref{thm:meta-dynamics-regret-and-stability} and \ref{thm:gradient-error-bound}, we show a local regret guarantee of our meta-framework in \Cref{thm:application-matched-disturbance-main} and test it numerically in the setting of Example \ref{example:nonlinear-control}. Due to space limit, we defer the proof of \Cref{thm:application-matched-disturbance-main} to Appendix \ref{appendix:application-matched-disturbance-main} and the simulation results to Appendix \ref{appendix:simulation}.

\begin{theorem}\label{thm:application-matched-disturbance-main}
Under Assumptions \ref{assump:Lipschitz-and-smoothness}-\ref{assump:estimator-environment}, if we use M-GAPS for \texttt{ALG} and Gradient Estimator for \texttt{EST}, the trajectory $\xi = \{x_t, y_t, \theta_t\}$ achieves an expected local regret \footnote{We change the gradient $\nabla F_t(\theta_t)$ to the projected gradient $\nabla_{\eta, \Theta} F_t(\theta_t)$ (see Definition \ref{def:proj-grad} in Appendix \ref{appendix:notation-and-defs}) in the local regret. This metric is introduced by \cite{hazan2017efficient} for online nonconvex optimization with constraints.} of
\[O\left(\eta^{-1}(1 + V_{sys} + \bar{\epsilon} \cdot T) + \eta T + (\sqrt{m \bar{\epsilon}} + m \bar{\epsilon}) \cdot T\right),\]
where $V_{sys}$ is the total variation of the system and $\bar{\epsilon}$ is the magnitude of the random disturbance $w_t$ (see Appendix \ref{appendix:application-details} for detailed definitions). 
Under the same definition of $\tilde{\xi}$ as \Cref{thm:meta-dynamics-regret-and-stability}, the expected total distance between $\xi$ and $\tilde{\xi}$ is bounded by
\[\mathbb{E}\left[\sum_{t=0}^{T-1} \left(\norm{x_t - \tilde{x}_t} + \norm{y_t - \tilde{y}_t}\right)\right] = O\left(T^{3/4} + \sqrt{m \bar{\epsilon}} \cdot T\right).\]
\end{theorem}

\section{Conclusion}
\label{sec:conclusion}
In this work, we study online policy optimization with unknown dynamics.
We propose a meta-framework combining an online policy optimization algorithm \texttt{ALG} with an online parameter estimator \texttt{EST}.
We specify a set of properties that, if satisfied, imply that \texttt{ALG} can act as if \texttt{EST} is providing the true dynamics models, while \texttt{EST} only needs to minimize the prediction errors on the actual trajectory visited by \texttt{ALG}.
To demonstrate our framework, we propose an efficient candidate for \texttt{ALG} called M-GAPS (Algorithm \ref{alg:M-GAPS}). When our meta-framework is instantiated with M-GAPS (as \texttt{ALG}) and the gradient estimator (as \texttt{EST}), it achieves the first local regret bound (\Cref{thm:application-matched-disturbance-main}) for online policy optimization in a class of nonlinear time-varying systems with unknown dynamics.

Our meta-framework also motivates interesting future work:
The structural properties of \texttt{ALG} and \texttt{EST} that we identify can serve as guidelines for the design of improved online policy optimization and dynamics estimations methods.
For example, tools from online nonconvex optimization and adaptive control could be leveraged to handle more general classes of dynamics.

%\section{Proof Outline}
%\label{sec:outline}
%\input{proof_outline}

\newpage


\bibliography{main}

%%%%%%%%%%%%%%%%%%%%%%%%%%%%%%%%%%%%%%%%%%%%%%%%%%%%%%%%%%%%%%%%%%%%%%%%%%%%%%%
%%%%%%%%%%%%%%%%%%%%%%%%%%%%%%%%%%%%%%%%%%%%%%%%%%%%%%%%%%%%%%%%%%%%%%%%%%%%%%%
% APPENDIX
%%%%%%%%%%%%%%%%%%%%%%%%%%%%%%%%%%%%%%%%%%%%%%%%%%%%%%%%%%%%%%%%%%%%%%%%%%%%%%%
%%%%%%%%%%%%%%%%%%%%%%%%%%%%%%%%%%%%%%%%%%%%%%%%%%%%%%%%%%%%%%%%%%%%%%%%%%%%%%%
\newpage
\appendix

\noindent{\textbf{Outline of the appendices.}}
\begin{itemize}[leftmargin=*]
    \item In Appendix \ref{appendix:simulation}, we provide the simulation results for a nonlinear control application (Example \ref{example:nonlinear-control}).
    \item In Appendix \ref{appendix:notation-and-defs}, we provide a notation table and list important definitions used in the proofs.
    \item In Appendix \ref{appendix:application-details}, we present the details about the application of matched-disturbance dynamics.
    \item In Appendix \ref{appendix:meta-dynamics-regret-and-stability}, we prove \Cref{thm:meta-dynamics-regret-and-stability} on the joint dynamics of \texttt{ALG} in our meta-framework.
    \item In Appendix \ref{appendix:gradient-error-bounds}, we prove \Cref{thm:gradient-error-bound} on the regret of \texttt{EST} in our meta-framework.
    \item In Appendix \ref{appendix:fully-actuated-assumptions}, we show M-GAPS satisfies the properties for \texttt{ALG} in the application of matched-disturbance dynamics.
    \item In Appendix \ref{appendix:grad-est}, we show the gradient estimator satisfies the properties for \texttt{EST} in the application of matched-disturbance dynamics.
    \item In Appendix \ref{appendix:application-matched-disturbance-main}, we prove \Cref{thm:application-matched-disturbance-main} about instantiating our meta-framework with M-GAPS (for \texttt{ALG}) and the gradient estimator (for \texttt{EST}) in the application of matched-disturbance dynamics.
    \item In Appendix \ref{appendix:nonconvex-biased-OGD}, we show online gradient descent with inexact updates can achieve local regret bounds in online nonconvex optimization, which is used in the proof of M-GAPS.
    \item In Appendix \ref{appendix:useful-lemmas}, we make a detailed comparison between our M-GAPS algorithm with the GAPS algorithm proposed by \cite{lin2023online}. We summarize some results from \cite{lin2023online} that are useful for us to analyze M-GAPS.
\end{itemize}

\section{Simulation Results}\label{appendix:simulation}
\begin{figure}[h]
\centering
\begin{minipage}{.5\textwidth}  
\centering
  \includegraphics[width=0.8\linewidth]{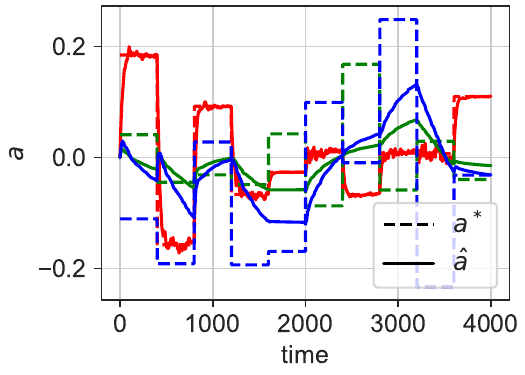}
  \caption{$a^*$ and $\hat{a}$}
  \label{fig:example_a}
    \includegraphics[width=0.8\linewidth]{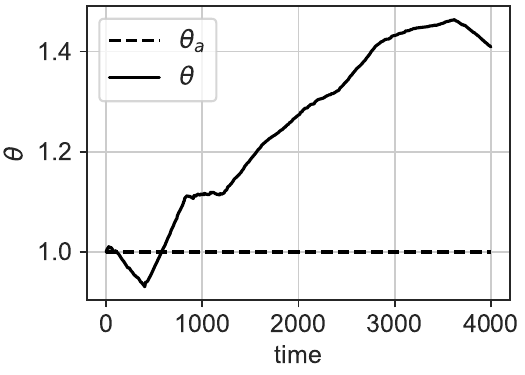}
  \caption{$\theta_a$ and $\hat{\theta}$}
  \label{fig:example_theta}
\end{minipage}%
\begin{minipage}{.5\textwidth} 
\centering
  \includegraphics[width=0.8\linewidth]{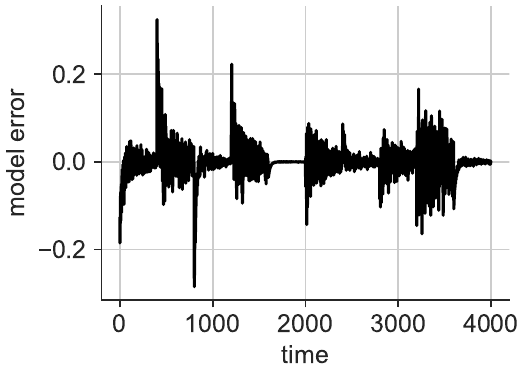}
  \caption{Model error $\phi(x)\hat{a} - \phi(x)a^*$}
  \label{fig:example_error}
  \includegraphics[width=0.8\linewidth]{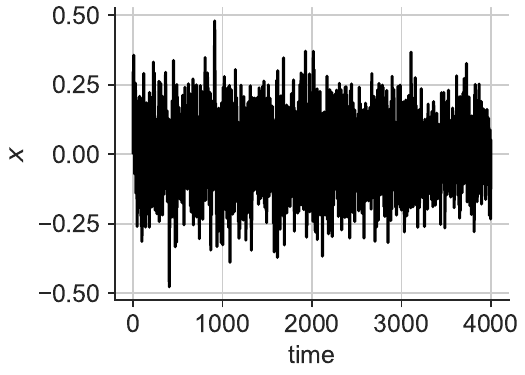}
  \caption{Tracking error $x$}
  \label{fig:example_x}
\end{minipage}
\begin{minipage}{1.\textwidth} 
\centering
    \includegraphics[width=0.6\linewidth]{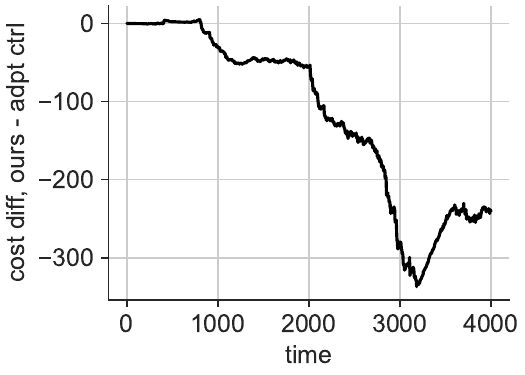}
    \caption{Regret difference between M-GAPS and adaptive control}
    \label{fig:example_costs}
\end{minipage}
\end{figure}

In this section, we show our simulation results in the setting of Example \ref{example:nonlinear-control}, where the dynamics model, control policy and the update law of the online parameter estimator \texttt{EST} are given by
\begin{subequations}\label{appendix:simulation:dynamics-e1}
\begin{align}
    x_{t + 1} &= x_{t} + \Delta(u_t + \phi(x_t)a^*_t) + w_t\\
    u_t &= -\theta_t x_t - \phi(x_t)\hat{a}_t\\
    \hat{a}_{t + 1} &= \hat{a}_t - (\phi(x_t)\hat{a}_t - \phi(x_t)a^*_t)\Delta P\phi(x_t), \label{appendix:simulation:dynamics-e1:s2}
\end{align}
\end{subequations}
where $\phi(x) = \begin{bmatrix}1 & \sin(x) & \sin(2x)\end{bmatrix}$ is a basis function, $\Delta = 0.01$ is the time interval between steps and $P$ is a constant gain.
The disturbances $w_{0:T-1}$ are generated by Ornstein-Uhlenbeck random walk. 
The Ornstein-Uhlenbeck random walk updates the random disturbance $w_t$ following the dynamics $w_{t+1} = \gamma w_t + \delta_t$, where parameter $\gamma \in [0, 1)$ is a constant decay factor and $\delta_t$ is Gaussian-distributed with zero mean and variance $\sigma^2$. In the simulation, we set $\gamma = 0.95$ and $\sigma = 0.1$. We consider a time-varying dynamical system where the true model parameter $a_t^*$ changes at the end of each period of $40000$ time steps. When a period ends, we resample $a_t^*$ from a uniform distribution over the closed interval $[-0.25, 0.25]^3$.
\Cref{fig:example_a} contains a plot of time-varying $a_t^*$.
The stage cost function is given by $h_t(x_t, u_t) = x_t^2 + \beta u_t^2$, where we set $\beta = 0.1$.

We test our meta-framework by
1) instantiating \texttt{ALG} using M-GAPS with the learning rate $\eta = 1 \times 10^{-4}$,
and 2) instantiating \texttt{EST} with the gradient estimator (\Cref{alg:grad-est}), whose update is given by \eqref{appendix:simulation:dynamics-e1:s2}. For \texttt{EST}, we make a minor modification to the prediction loss function by changing it from the $\ell_2$-squared prediction error $\norm{\phi(x_t)\hat{a}_t - \phi(x_t)a^*_t}^2$ to $\norm{\phi(x_t)\hat{a}_t - \phi(x_t)a^*_t}_P^2$, where $\norm{\cdot}_P$ denotes the norm induced by matrix $P$.
Note that the update law of the gradient estimator in this setting
g is identical with a classic adaptive controller \citep[][\S 8.7.3]{slotine1991applied}.
We compare our algorithm with a common benchmark that uses the adaptive controller to learn $\hat{a}_t$ while the policy parameter (the feedback gain) $\theta_a$ is fixed:
\begin{subequations}\label{appendix:simulation:dynamics-e2}
\begin{align}
    x_{t + 1} &= x_{t} + \Delta(u_t + \phi(x_t)a^*_t) + w_t\\
    u_t &= -\theta_a x_t - \phi(x_t)\hat{a}_t\\
    \hat{a}_{t + 1} &= \hat{a}_t - (\phi(x_t)\hat{a}_t - \phi(x_t)a^*_t)\Delta P\phi(x_t).
\end{align}
\end{subequations}
The result is shown in Figure~\ref{fig:example_a}-\ref{fig:example_costs}. Similar to traditional adaptive controllers, adaptive parameters $\hat{a}_t$ do not necessarily converge to the real value $a_t^*$. Still, our model and tracking errors converge to a small error ball. Our algorithm optimizes the cost function with control input and improves the accumulative costs significantly compared with the baseline adaptive controller \eqref{appendix:simulation:dynamics-e2}. 

\section{Notations and Definitions}\label{appendix:notation-and-defs}
We provide a notation table (Table \ref{table:notations}) that summarizes the important notations in this paper.

\begin{table}
  \centering
  \caption{Important notations in this paper.}\label{table:notations}
  \begin{tabular}{c|l}
    \specialrule{1.5pt}{0pt}{0pt}
    \textbf{Notation} & \textbf{Meaning} \\
    \specialrule{0.4pt}{0pt}{0pt}
    $t_1:t_2$ & The integer sequence $\{t_1, \ldots, t_2\}$;  \\
    $a_{t_1:t_2}$ & A sequence of variables $\{a_t\}_{t = t_1, \ldots, t_2}$; \\
    $\|\cdot\|$ & $\ell_2$ (Euclidean) norm; \\
    $\|\cdot\|_F$ & Frobenius norm; \\
    $\|\cdot\|_P$ & Norm induced by matrix $P$; \\
    $\mathbb{Z}_{\geq 0}$ & The set of non-negative integers; \\
    $\mathbb{R}_{\geq 0}$ & The set of non-negative reals; \\
    $\sigma(z_{1:t}, z'_{1:t})$ & Product sigma-algebra generated by sequences $z_{1:t}$ and $z'_{1:t}$; \\
    $x_t$ & $x_t \in\mathbb{R}^n$ is the system state; \\
    $u_t$ & $u_t \in \mathbb{R}^m$ is the control input; \\
    $w_t$ & $w_t \in \mathcal{W} \subseteq \mathbb{R}^n$ is a disturbance term;\\
    $f_t(x_t,a_t^*)$ & $f_t$ is a nonlinear residual term where the online agent makes (noisy) observations;\\
    $a_t^*$ & $a_t^* \in \mathcal{A} \subseteq \mathbb{R}^p$ is the unknown parameter in $f_t$ \\
    $f_t(\cdot,\hat{a}_t)$ & An estimation of the true nonlinear residual function $f_t(\cdot,a_t^*)$;\\
    $q_t(x_t,y_t,\theta_t,a_t^*)$ & The joint dynamics of the system at time $t$;\\
    $\Pi_\Theta(y)$ & Euclidean projection of $y$ to set $\Theta$;\\
    \specialrule{1.5pt}{0pt}{0pt}
  \end{tabular}
\end{table}

A key concept that we explore in this paper is how to compare the actual trajectory of our meta-framework with an ``ideal'' trajectory that the agent could achieve with exact knowledge of the true model parameters $a_{0:T-1}^*$, we introduce the important notations of multi-step dynamics/cost that characterize how the system would evolve under a sequence of policy parameters $\theta_{0:T-1}$ when $a_{0:T-1}^*$ is known. The concepts of multi-step dynamics/cost are first introduced in \cite{lin2023online}, which studies online policy selection with known dynamical systems. In this work, we replace all estimated $\hat{a}_t$ in the policy classes with true $a_t^*$ to reproduce the same definition as \cite{lin2023online}.

\begin{definition}[Multi-Step Dynamics and Cost]\label{def:multi-step-dynamics}
The multi-step dynamics $g_{t\mid \tau}^*$ between two time steps $\tau \leq t$ specifies the state $x_t$ as a function of the previous state $x_\tau$ and previous policy parameters $\theta_{\tau: t-1}$ under exact predictions $\{a_t^*\}$. It is defined recursively, with the base case $g_{\tau\mid \tau}^*(x_\tau) \coloneqq x_\tau$ and the recursive case 
\[
    g_{t+1\mid \tau}^*(x_\tau, \theta_{\tau:t}) = g_t\left(z_t, \ALG_t\left(z_t, \theta_{t}, f_t(z_t, a_t^*)\right), f_t(z_t, a_t^*)\right) + w_t,\ \forall\,t\geq \tau, 
\]
in which $z_t \coloneqq g_{t\mid \tau}^*(x_\tau, \theta_{\tau:t-1})$.\footnote{$z_t$ is an auxiliary variable to denote the state at $t$ under initial state $x_\tau$ and parameters $\theta_{\tau:t}$.}
% \james{I want to use larger parens to make the nesting more clear, but it wastes precious space...}
The multi-step cost $h_{t\mid \tau}^*$ specifies the cost $c_t$ as function of $x_\tau$ and $\theta_{\tau:t}$. It is defined as
\[h_{t\mid \tau}^*(x_\tau, \theta_{\tau:t}) \coloneqq h_t\left(z_t, \ALG_t\left(z_t, \theta_{t}, f_t(z_t, a_t^*)\right), \theta_t\right).\]
\end{definition}
It is worth emphasizing that, in our work, the concepts of multi-step dynamics/cost are only used for the theoretical analysis, because their definitions involve true model parameters that are unknown to the online agent. When doing online policy optimization, the online agent may use the estimations $\hat{g}_{t+1\mid t}$ and $\hat{h}_{t\mid t}$ (see \eqref{equ:g_hat_h_hat}) as the estimations of $g_{t+1\mid t}^*$ and $h_{t\mid t}^*$ respectively. Note that this is different than the case when true dynamics are known \citep{lin2023online}, where the online agent can directly construct multi-step dynamics/cost ($g_{t+1\mid t}^*$ and $h_{t\mid t}^*$) or compute the exact Jacobian matrices. 

Another important definition that we require is the \textit{projected gradient}, which is introduced in \cite{hazan2017efficient} to accommodate the challenge of converging to stationary points on a constrained set.

\begin{definition}[Projected gradient]\label{def:proj-grad}
Let $F: \Theta \to \mathbb{R}$ be a differentiable function on a closed convex set $\Theta \subseteq \mathbb{R}^d$. For $\eta > 0$, the $(\Theta, \eta)$-projected gradient of $F$ is defined as
\[\nabla_{\Theta, \eta}F(\theta) \coloneqq \frac{1}{\eta}\left(\theta - \Pi_\Theta(\theta - \eta \nabla F(\theta))\right).\]
\end{definition}

When $\Theta$ is equal to the whole Euclidean space $\mathbb{R}^d$ (unconstrained), the project gradient in Definition \ref{def:proj-grad} will be identical with the normal gradient $\nabla F(\theta)$. This concept of projected gradient is used to define the local regret in \Cref{thm:application-matched-disturbance-main} and Appendix \ref{appendix:nonconvex-biased-OGD} that studies online gradient descent for online nonconvex optimization with constraints.

\section{Matched-Disturbance Dynamics}
\label{appendix:application-details}
In this section, we discuss the detailed assumptions on the application of matched-disturbance dynamics and how they enable us to apply the theory for our meta-framework.

We first focus on the \texttt{ALG} part that is instantiated with M-GAPS (\Cref{alg:M-GAPS}). Note that in the setting of this application, the definition of multi-step dynamics (Definition \ref{def:multi-step-dynamics}) can be simplified to
\[g_{t+1\mid \tau}^*(x_\tau, \theta_{\tau:t}) = \phi_t(z_t, \psi_t(z_t, \theta_t)) + w_t, \text{ where } z_t \coloneqq g_{t\mid \tau}^*(x_\tau, \theta_{\tau:t-1}).\]
We see that $g_{t+1\mid \tau}^*$ has exactly the same form as the multi-step dynamics studied in \cite{lin2023online}. Since $g_{t+1\mid \tau}^*$ corresponds to the case when the nonlinear residual has been canceled out with the exact model parameter $a_t^*$, it is reasonable to make the same assumptions as \cite{lin2023online} about the stability and contraction properties of $g_{t+1\mid \tau}^*$. These assumptions rely on the key definitions of time-varying contractive perturbation and time-varying stability.

\begin{definition}\label{def:constrained-policy-parameter-seq}
We denote the set of policy parameter sequences with $\epsilon_\theta$-constrained step size by
\[S_{\epsilon_\theta}(t_1:t_2) \coloneqq \{\theta_{t_1:t_2} \in \Theta^{t_2 - t_1 + 1}\mid \norm{\theta_{\tau+1} - \theta_{\tau}} \leq \epsilon_\theta, \forall \tau \in [t_1 \mathbin{:} t_2 - 1]\}.\]
\end{definition}

\begin{definition}[$\epsilon_\theta$-Time-varying Contractive Perturbation]\label{def:epsilon-exp-decay-perturbation-property}
For $\epsilon_\theta \geq 0$, the $\epsilon_\theta$-time-varying contractive perturbation property holds for $R_C > 0, C > 0$, and $\decayfactor \in (0, 1)$ if, for any $\theta_{\tau:t-1} \in S_{\epsilon_\theta}(\tau:t-1)$, the following inequality holds for arbitrary $x_\tau, x_\tau' \in B_n(0, R_C)$ and time steps $\tau \leq t$:
\begin{align*}
    \norm{g_{t\mid \tau}^*(x_\tau, \theta_{\tau:t-1}) - g_{t\mid \tau}^*(x_\tau', \theta_{\tau:t-1})} \leq C \decayfactor^{t - \tau} \norm{x_\tau - x'_\tau}
\end{align*}
\end{definition}

\begin{definition}[$\epsilon_\theta$-Time-varying Stability]\label{def:epsilon-time-varying-stability}
For $\epsilon_\theta \geq 0$, the $\epsilon_\theta$-time-varying stability property holds for $R_S > 0$ if, for any $\theta_{\tau:t-1} \in S_{\epsilon_\theta}(\tau:t-1)$, $\norm{g_{t\mid \tau}^*(0, \theta_{\tau:t-1})} \leq R_S$ holds for any $t\geq \tau$.
\end{definition}

With the definitions of time-varying contractive perturbation and time-varying stability, we state our key assumptions below: Assumption \ref{assump:Lipschitz-and-smoothness} is about the Lipschitzness/smoothness properties of the dynamics, policy, nonlinear residual, and the cost functions.

\begin{assumption}\label{assump:Lipschitz-and-smoothness}
The dynamics $\phi_{0:T-1}$, policies $\psi_{0:T-1}$, residuals $f_{0:T-1}$, and costs $h_{0:T-1}$ are differentiable at every time step and satisfy that, for any convex compact sets $\mathcal{X} \subseteq \mathbb{R}^n, \mathcal{U} \subseteq \mathcal{R}^m$, one can find Lipschitzness/smoothness constants (can depend on $\mathcal{X}$ and $\mathcal{U}$) such that:
\begin{enumerate}[wide, labelindent=0pt, nolistsep]
    \item $\phi_t(x, u)$ is $(L_{\phi, x}, L_{\phi, u})$-Lipschitz and $(\ell_{\phi, x}, \ell_{\phi, u})$-smooth in $(x, u)$ on $\mathcal{X} \times \mathcal{U}$.
    \item $\psi_t(x, \theta)$ is $(L_{\psi, x}, L_{\psi, \theta})$-Lipschitz and $(\ell_{\psi, x}, \ell_{\psi, \theta})$-smooth in $(x, \theta)$ on $\mathcal{X} \times \Theta$.
    \item $f_t(x, a)$ is $(L_{f, x}, L_{f, a})$-Lipschitz and $(\ell_{f, x}, \ell_{f, a})$-smooth in $(x, a)$ on $\mathcal{X} \times \mathcal{A}$.
    \item $h_t(x, u, \theta)$ is $(L_{h, x}, L_{h, u}, L_{h, \theta})$-Lipschitz and $(\ell_{h, x}, \ell_{h, u}, \ell_{h, \theta})$-smooth in $(x, u, \theta)$ on $\mathcal{X} \times \mathcal{U} \times \Theta$.
\end{enumerate}
\end{assumption}

Compared with Assumption 2.1 in \cite{lin2023online}, our Assumption \ref{assump:Lipschitz-and-smoothness} additionally assumes the Lipschitzness and smoothness of the nonlinear residual function $f_t$, which is part of our dynamics and policy classes. The second assumption (Assumption \ref{assump:contractive-and-stability}) is on the contractive perturbation and the stability of the multi-step dynamics $g_{t\mid \tau}^*$.

\begin{assumption}\label{assump:contractive-and-stability}
Let $\mathcal{G}$ denote the set of all possible sequences $\{\phi_t, f_t, w_t, \psi_t\}_{t\in\mathcal{T}}$ the environment may provide. For a fixed $\epsilon_\theta \in \mathbb{R}_{\geq 0}$, the $\epsilon_\theta$-time-varying contractive perturbation holds with $(R_C, \bar{C}, \rho)$ for any sequence in $\mathcal{G}$. The $\epsilon_\theta$-time-varying stability holds with $R_S < R_C$ for any sequence in $\mathcal{G}$. We assume that the initial state satisfies $\norm{x_0} < (R_C - R_S)/\bar{C}$. Further, we assume that if $\{\phi, f, w, \psi\}$ is the dynamics/residual/disturbance/policy at an intermediate time step of a sequence in $\mathcal{G}$, then the time-invariant sequence $\{\phi, f, w, \psi\}_{\times T}$ is also in $\mathcal{G}$.\footnote{For $\{\phi, f, w, \psi\}_{\times T}$ to be in $\mathcal{G}$, it must satisfy other assumptions about contractive perturbation and stability that we impose on $\mathcal{G}$ but does not need to occur in real problem instances. This assumption can be made without the loss of generality for time-invariant dynamics and policy classes.}
\end{assumption}

Compared with Assumption 2.2 in \cite{lin2023online}, our Assumption \ref{assump:contractive-and-stability} also includes the disturbance $w_t$ as a part of the system configuration. This is because for every time $t$, $g_{t+1\mid t}^*$ is formed by $\phi_t, \pi_t,$ and $w_t$ together. While $w_t$ can also be viewed as a part of the dynamics $\phi_t$, we choose to represent it separately because we will leverage the randomness of $w_t$ to bound the first-order model mismatches of \texttt{EST}. Like \cite{lin2023online}, in Assumption \ref{assump:contractive-and-stability}, we assume there exists a positive real number $\bar{R}_C$ such that $R_C > \bar{R}_C > R_S + \bar{C}\norm{x_0}$. Here, we introduce the real constant $\bar{R}_C$ because $R_C$ can be $+\infty$ when time-varying contractive perturbation (Definition \ref{def:epsilon-exp-decay-perturbation-property}) holds globally. Similarly, to leverage the Lipschitzness/smoothness property, we require ${\mathcal{X} \supseteq B(0, R_x)}$ where $R_x \geq \bar{C} \bar{R}_C + R_S$ and $\mathcal{U} = \{- f(x, a) + \pi(x, \theta) \mid x \in \mathcal{X}, \theta \in \Theta, a \in \mathcal{A}, \pi, f \in \mathcal{G}\}$.
Since the coefficients in Assumption \ref{assump:Lipschitz-and-smoothness} depend on $\mathcal{X}$ and $\mathcal{U}$, we will set $\mathcal{X} = B_n(0, R_x)$ and $R_x = \bar{C} \bar{R}_C + R_S$ by default when presenting these constants. We also set $\mathcal{Y} = B_p(0, R_y)$ with $R_y = \frac{\bar{C} L_{\phi, u} L_{\psi, \theta}}{\rho(1 - \rho)}$, so that the internal state $y_t$ will stay in $\mathcal{Y}$.

It is straightforward to verify that the joint dynamics of M-GAPS satisfy the three properties required by the meta-framework. We state this result in Lemma \ref{lemma:fully-actuated-assumptions}.

\begin{lemma}\label{lemma:fully-actuated-assumptions}
Under Assumptions \ref{assump:Lipschitz-and-smoothness} and \ref{assump:contractive-and-stability}, M-GAPS (Algorithm \ref{alg:M-GAPS}) satisfy Properties \ref{assump:Lipschitzness-a}, \ref{assump:contraction-and-stability}, and \ref{assump:robustness} when applied to dynamics \eqref{equ:dynamics-fully-actuated} and policy class \eqref{equ:policy-cancel-out}.
\end{lemma}

\noindent We present the specific constants and the formal proof of Lemma \ref{lemma:fully-actuated-assumptions} in Appendix \ref{appendix:fully-actuated-assumptions}.

For the part of \texttt{EST} that is instantiated with the gradient estimator (\Cref{alg:grad-est}), we first introduce an assumption about the magnitude of the nonlinear residual to guarantee that (several) bad estimations of the unknown model parameters will not destabilize the system or violate the constraints of the contractive perturbation property.

\begin{assumption}\label{assump:linear-function-approximation}
The set of all possible model parameter $\mathcal{A}$ is a convex compact subset of $\mathbb{R}^k$. For any fixed $x \in B_n(0, R_x)$, $f_t(x, \cdot): \mathcal{A} \to \mathbb{R}^m$ is an affine function whose gradient is uniformly bounded, i.e., for some positive constant $D_f'$,
$\norm{\nabla_x f_t(x, a)} \leq D_f'$ hold for all $a \in \mathcal{A}$.
It also satisfies that for any $a, a' \in \mathcal{A}$,
\begin{align*}
    &\norm{f_t(x, a) - f_t(x, a')} \leq C_f, \ \norm{\nabla_x f_t(x, a) - \nabla_x f_t(x, a')}_F \leq \beta \leq C_f', \text{ and }\\
    &\norm{\nabla_x^2 f_t(x, a)_i - \nabla_x^2 f_t(x, a')_i}_F \leq \gamma, \text{ for any dimension }i \in [1:m].
\end{align*}
hold with some positive constants $\beta, \gamma$, and the upper bounds $C_f$ and $C_f'$ are given by
\begin{align*}
    C_f ={}& \min\left\{\frac{\sqrt{2}\bar{\zeta}}{4(L_{\theta, x} + L_{\theta, y})C \alpha_x}, \frac{\min\{R_x - R_x^*, R_y - R_y^*\}}{C \alpha_x}, \frac{\bar{\zeta}}{2 \alpha_\theta}\right\},\\
    C_f' ={}& \min\left\{\frac{\sqrt{2}\bar{\zeta}}{4(L_{\theta, x} + L_{\theta, y})C \beta_x}, \frac{\min\{R_x - R_x^*, R_y - R_y^*\}}{C \beta_x}, \frac{\bar{\zeta}}{2 \beta_\theta}\right\}.
\end{align*}
The expressions of $\alpha_x, \beta_x, \alpha_\theta, \beta_\theta, L_{\theta, x}, L_{\theta, y}, R_x^*, R_y^*,$ and $\bar{\zeta}$ are given in Appendix \ref{appendix:fully-actuated-assumptions}.
\end{assumption}

Note that we need Assumption \ref{assump:linear-function-approximation} to bound the prediction errors uniformly because even if an online parameter estimator performs well in the long term (e.g., achieving a regret bound on the total prediction errors), it may incur a large error at a single time step that can potentially destabilize the system especially when $a_t^*$ changes abruptly. Our simulation (Appendix \ref{appendix:simulation}) provides a good illustration of this intuition: The model prediction error may increase dramatically right after the system switches to a different true model parameter $a_t^*$; Then, the error converges back to near zero as the gradient estimator learns the model (see Figure \ref{fig:example_error}). Addressing this challenge with other assumptions like slowly time-varying $a_t^*$ is an interesting future direction.

The second assumption we need is about the randomness provided by the environment:

\begin{assumption}\label{assump:estimator-environment}
The total path length of the true model parameters satisfies
\[1 + \sum_{t=1}^{T-1}\norm{a_t^* - a_{t-1}^*} \leq C_p\] 
for some positive constant $C_p$. At every time step $t$, the noisy observation $\tilde{f}_t$ satisfies that \[\norm{\tilde{f}_t - f_t(x_t, a_t^*)} \leq e_f, \text{ and }\mathbb{E}[\tilde{f}_t\mid \mathcal{F}_t] = f_t(x_t, a_t^*).\]
Further, the random disturbance $w_t$ in the dynamical system \eqref{equ:dynamics-fully-actuated} satisfies that
\[\norm{w_t} \leq \bar{\epsilon},\ \mathbb{E}[w_t \mid \mathcal{F}_t'] = 0, \text{ and } \text{Cov}(w_t\mid \mathcal{F}_t') \succeq c \bar{\epsilon}^2 I.\]
Here, $\bar{\epsilon}$ satisfies that
\[(C_f + e_f)\left(2 D_f' \sqrt{3 C_p/T}\right)^{\frac{1}{3}} \leq \bar{\epsilon} \leq \min\{\frac{1}{4}, \frac{1}{2\gamma}, \frac{1}{4 \beta \gamma}\},\]
where $\beta$ and $\gamma$ are defined in Assumption \ref{assump:linear-function-approximation}.
\end{assumption}

Intuitively, Assumption \ref{assump:estimator-environment} put requirements on both the lower and upper bounds of the level of randomness in the system. The lower bound $\bar{\epsilon} = \Omega(T^{-1/6})$ is required due to the condition $R_0^\ell(T) \leq \bar{\epsilon}^3 T$ in \Cref{thm:gradient-error-bound}. This guarantees that the zeroth-order regret is sufficiently small to be used for bounding the first-order gradients in Taylor's expansion, which are multiplied by $\bar{\epsilon}^2$ when we take the square. The upper bound of $\bar{\epsilon}$ is required to ignore the higher-order terms in Taylor's expansion. With these assumptions, we show the following guarantee on the total prediction error achieved by the gradient estimator (\Cref{alg:grad-est}).

\begin{lemma}\label{lemma:gradient-estimator-regret}
Under Assumptions \ref{assump:linear-function-approximation} and \ref{assump:estimator-environment}, the total squared zeroth-order prediction errors of the gradient estimator can be bounded by
$\mathbb{E}\left[\sum_{t=0}^{T-1} \varepsilon_t^2\right] \leq 2\sqrt{3} (C_f + e_f)^3 D_f' \sqrt{C_p T}.$
\end{lemma}

We defer the proof of Lemma \ref{lemma:gradient-estimator-regret} to Appendix \ref{appendix:grad-est}. Note that our meta-framework only requires us to bound the total squared zeroth-order prediction error incurred by an instantiation of \texttt{EST}. Under Assumptions \ref{assump:linear-function-approximation} and \ref{assump:estimator-environment}, we can apply \Cref{thm:gradient-error-bound} in the meta-framework to bound the total squared first-order prediction error of the gradient estimator by $\mathbb{E}\left[\sum_{t=0}^{T-1} (\varepsilon_t')^2\right] = O(m\bar{\epsilon} T)$. Recall that $m$ is the dimension of the unknown component, which is identical with the control input in this application.

In Lemmas \ref{lemma:fully-actuated-assumptions} and \ref{lemma:gradient-estimator-regret}, we have shown that M-GAPS and the gradient estimator satisfy all the required properties for \texttt{ALG} and \texttt{EST} respectively in our meta-framework. Therefore, we can obtain the local regret guarantees for instantiating our meta-framework with M-GAPS and the gradient estimator in the application with matched-disturbance dynamics (\Cref{thm:application-matched-disturbance-main}).

\section{Proof of Theorem \ref{thm:meta-dynamics-regret-and-stability}}
\label{appendix:meta-dynamics-regret-and-stability}
To simplify the notation, we define
\[\bar{\varepsilon} \coloneqq \min\left\{\frac{\sqrt{2}\bar{\zeta}}{4(L_{\theta, x} + L_{\theta, y})C}, \frac{\min\{R_x - R_x^*, R_y - R_y^*\}}{C}\right\}.\]
By the assumption, we know that the following inequality holds for all time step $t$:
\begin{align}\label{thm:meta-dynamics-regret-and-stability:e0}
    (\alpha_x + \alpha_y)\varepsilon_t + (\beta_x + \beta_y)\varepsilon_t' \leq \bar{\varepsilon}.
\end{align}
Now we show that $\norm{x_t} \leq R_x$, $\norm{y_t} \leq R_y$, and $\norm{\theta_t - \theta_{t-1}} \leq \epsilon_\theta$ by induction. These inequalities hold for time step $0$. Suppose they hold for all time steps $\tau \leq t$. Then, for time step $t+1$, by \Cref{assump:contraction-and-stability}, we see that
\begin{align}\label{thm:meta-dynamics-regret-and-stability:e0-1}
    \norm{\tilde{x}_{t+1}} \leq R_x^*, \text{ and } \norm{\tilde{y}_{t+1}} \leq R_y^*.
\end{align}
By Property \ref{assump:contraction-and-stability} about the contraction of states $x_t$ and $y_t$ under policy parameters $\theta_{0:t}$, we see that
\begin{subequations}\label{thm:meta-dynamics-regret-and-stability:e1}
\begin{align}
    &\norm{(x_{t+1}, y_{t+1}) - (\tilde{x}_{t+1}, \tilde{y}_{t+1})}\nonumber\\
    \leq{}& \sum_{\tau = 0}^{t} \norm{q_{t+1\mid t+1-\tau}^{(x, y)*}(x_{t+1-\tau}, y_{t+1-\tau}, \theta_{t+1-\tau:t}) - q_{t+1\mid t-\tau}^{(x, y)*}(x_{t-\tau}, y_{t-\tau}, \theta_{t-\tau:t})}\label{thm:meta-dynamics-regret-and-stability:e1:s1}\\
    \leq{}& \sum_{\tau = 0}^{t} \gamma(\tau) \norm{(x_{t+1-\tau}, y_{t+1-\tau}) - q_{t+1-\tau\mid t-\tau}^{(x, y)*}(x_{t-\tau}, y_{t-\tau}, \theta_{t-\tau})}\label{thm:meta-dynamics-regret-and-stability:e1:s2}\\
    \leq{}& \sum_{\tau = 0}^{t} \gamma(\tau) \left((\alpha_x + \alpha_y)\varepsilon_{t-\tau} + (\beta_x + \beta_y)\varepsilon_{t-\tau}'\right)\label{thm:meta-dynamics-regret-and-stability:e1:s3}\\
    \leq{}& \sum_{\tau = 0}^{t} \gamma(\tau) \bar{\varepsilon} \leq C \bar{\varepsilon} \label{thm:meta-dynamics-regret-and-stability:e1:s4}\\
    \leq{}& \min\{R_x - R_x^*, R_y - R_y^*\}, \label{thm:meta-dynamics-regret-and-stability:e1:s5}
\end{align}
\end{subequations}
where we use the triangle inequality in \eqref{thm:meta-dynamics-regret-and-stability:e1:s1};
we use the contractive perturbation property in Property \ref{assump:contraction-and-stability}
in \eqref{thm:meta-dynamics-regret-and-stability:e1:s2};
we use the induction assumption and \Cref{assump:Lipschitzness-a} in \eqref{thm:meta-dynamics-regret-and-stability:e1:s3}; we use \eqref{thm:meta-dynamics-regret-and-stability:e0} in \eqref{thm:meta-dynamics-regret-and-stability:e1:s4} and the definition of $\bar{\varepsilon}$ in \eqref{thm:meta-dynamics-regret-and-stability:e1:s5}. By \eqref{thm:meta-dynamics-regret-and-stability:e0-1} and \eqref{thm:meta-dynamics-regret-and-stability:e1}, we see that
\begin{align}\label{thm:meta-dynamics-regret-and-stability:e2}
    \norm{x_{t+1}} \leq{}& \norm{\tilde{x}_{t+1}} + \norm{\tilde{x}_{t+1} - x_{t+1}} \leq R_x, \text{ and }\nonumber\\
    \norm{y_{t+1}} \leq{}& \norm{\tilde{y}_{t+1}} + \norm{\tilde{y}_{t+1} - y_{t+1}} \leq R_y.
\end{align}
Note that we can construct the disturbance sequence $\{\zeta_t\}$ in \Cref{assump:robustness} such that the dynamics
\begin{align*}
    \begin{pmatrix}
        \tilde{x}_{t+1}\\
        \tilde{y}_{t+1}\\
        \theta_{t+1}
    \end{pmatrix} = q_t(\tilde{x}_t, \tilde{y}_t, \theta_t, a_t^*) + \begin{pmatrix}
        0\\
        0\\
        \zeta_t
    \end{pmatrix}
\end{align*}
produce the same policy parameter sequence $\{\theta_t\}$ as the dynamics
\begin{align*}
    \begin{pmatrix}
        x_{t+1}\\
        y_{t+1}\\
        \theta_{t+1}
    \end{pmatrix} = q_t(x_t, y_t, \theta_t, \hat{a}_t).
\end{align*}
Therefore, under this construction, we see that
\begin{subequations}\label{thm:meta-dynamics-regret-and-stability:e3}
\begin{align}
    \norm{\zeta_t} \leq{}&\norm{\theta_{t+1} - q_t^\theta(\tilde{x}_t, \tilde{y}_t, \theta_t, a_t^*)}\nonumber\\
    ={}& \norm{q_t^\theta(x_t, y_t, \theta_t, \hat{a}_t) - q_t^\theta(\tilde{x}_t, \tilde{y}_t, \theta_t, a_t^*)}\nonumber\\
    \leq{}& \norm{q_t^\theta(x_t, y_t, \theta_t, \hat{a}_t) - q_t^\theta(x_t, y_t, \theta_t, a_t^*)} + \norm{q_t^\theta(x_t, y_t, \theta_t, a_t^*) - q_t^\theta(\tilde{x}_t, \tilde{y}_t, \theta_t, a_t^*)} \label{thm:meta-dynamics-regret-and-stability:e3:s1}\\
    \leq{}& \alpha_\theta \varepsilon_t + \beta_\theta \varepsilon_t' + L_{\theta, x} \norm{x_t - \tilde{x}_t} + L_{\theta, y} \norm{y_t - \tilde{y}_t} \label{thm:meta-dynamics-regret-and-stability:e3:s2}\\
    \leq{}& \alpha_\theta \varepsilon_t + \beta_\theta \varepsilon_t' + \sqrt{2}(L_{\theta, x} + L_{\theta, y}) \norm{(x_t, y_t) - (\tilde{x}_t, \tilde{y}_t)} \label{thm:meta-dynamics-regret-and-stability:e3:s3}\\
    \leq{}& \alpha_\theta \varepsilon_t + \beta_\theta \varepsilon_t' + \sqrt{2} C \bar{\varepsilon} (L_{\theta, x} + L_{\theta, y}) \leq \bar{\zeta}, \label{thm:meta-dynamics-regret-and-stability:e3:s4}
\end{align}
\end{subequations}
where we use the triangle inequality in \eqref{thm:meta-dynamics-regret-and-stability:e3:s1}; we use \Cref{assump:Lipschitzness-a} in \eqref{thm:meta-dynamics-regret-and-stability:e3:s2}; we use the inequality
\[\norm{x_t - \tilde{x}_t} + \norm{y_t - \tilde{y}_t} \leq \sqrt{2} \norm{(x_t, y_t) - (\tilde{x}_t, \tilde{y}_t)}\]
in \eqref{thm:meta-dynamics-regret-and-stability:e3:s3} and \eqref{thm:meta-dynamics-regret-and-stability:e1:s4} in \eqref{thm:meta-dynamics-regret-and-stability:e3:s4}. Thus, by \Cref{assump:robustness}, we see that $\norm{\theta_{t+1} - \theta_t} \leq \epsilon_\theta$. Therefore, we have shown that
\[\norm{x_t} \leq R_x, \norm{y_t} \leq R_y, \text{ and }\norm{\theta_t - \theta_{t-1}} \leq \epsilon_\theta\]
hold for all time step $t$ by induction.

By \eqref{thm:meta-dynamics-regret-and-stability:e3:s3} and \eqref{thm:meta-dynamics-regret-and-stability:e1:s3}, we also see that
\begin{align}\label{thm:meta-dynamics-regret-and-stability:e4}
    \norm{\zeta_t} \leq{}& \alpha_\theta \varepsilon_t + \beta_\theta \varepsilon_t' + \sqrt{2}(L_{\theta, x} + L_{\theta, y}) \norm{(x_t, y_t) - (\tilde{x}_t, \tilde{y}_t)}\nonumber\\
    \leq{}& \alpha_\theta \varepsilon_t + \beta_\theta \varepsilon_t' + \sqrt{2}(L_{\theta, x} + L_{\theta, y}) \sum_{\tau = 0}^{t-1} \gamma(\tau) \left((\alpha_x + \alpha_y)\varepsilon_{t-1-\tau} + (\beta_x + \beta_y)\varepsilon_{t-1-\tau}'\right).
\end{align}
Summing \eqref{thm:meta-dynamics-regret-and-stability:e4} over $t = 0, 1, \ldots, T-1$ gives that
\begin{align*}
    &\sum_{t=0}^{T-1} \norm{\zeta_t}\\
    \leq{}& \left(\alpha_\theta + \sqrt{2} C (L_{\theta, x} + L_{\theta, y})(\alpha_x + \alpha_y)\right) \sum_{t=0}^{T-1}\varepsilon_t + \left(\beta_\theta + \sqrt{2} C (L_{\theta, x} + L_{\theta, y})(\beta_x + \beta_y)\right) \sum_{t=0}^{T-1}\varepsilon_t'.
\end{align*}

Summing \eqref{thm:meta-dynamics-regret-and-stability:e1:s3} over $t = 0, 1, \ldots, T-1$ gives that
\begin{align*}
    \sum_{t=1}^{T} \norm{(x_{t}, y_{t}) - (\tilde{x}_{t}, \tilde{y}_{t})} \leq C \left((\alpha_x + \alpha_y)\sum_{t=0}^{T-1}\varepsilon_t + (\beta_x + \beta_y)\sum_{t=0}^{T-1}\varepsilon_t'\right).
\end{align*}

\section{Proof of Theorem \ref{thm:gradient-error-bound}}
\label{appendix:gradient-error-bounds}
To simplify the notation, we let $\check{x}_{t+1} \coloneqq \mathbb{E}[x_{t+1}\mid \mathcal{G}_t]$ and let $\iota_{t+1} \coloneqq x_{t+1} - \check{x}_{t+1}$.

We first focus on one dimension $i$ of the model mismatch. By Taylor's expansion, we see that
\begin{align}\label{thm:gradient-error-bound:e1}
    e_t(x_t, \hat{a}_t)_i = e_t(\check{x}_t, \hat{a}_t)_i + \nabla_x e_t(\check{x}_t, \hat{a}_t)_i \iota_t + \frac{1}{2} \iota_t^\top \nabla_x^2 e_t(\tilde{x}_t, \hat{a}_t)_i \iota_t,
\end{align}
where $\tilde{x}_t = \omega x_t + (1 - \omega) \check{x}_t$ for some $\omega \in [0, 1]$. Note that we have
\begin{align}\label{thm:gradient-error-bound:e2}
    \mathbb{E}\left[e_t(\check{x}_t, \hat{a}_t)_i \nabla_x e_t(\check{x}_t, \hat{a}_t)_i \iota_t \mid \mathcal{G}_{t-1}\right] = e_t(\check{x}_t, \hat{a}_t)_i \nabla_x e_t(\check{x}_t, \hat{a}_t)_i \mathbb{E}\left[\iota_t \mid \mathcal{G}_{t-1}\right] = 0.
\end{align}
Thus, we see that
\begin{subequations}\label{thm:gradient-error-bound:e3}
\begin{align}
    \mathbb{E}\left[e_t(x_t, \hat{a}_t)_i^2\mid \mathcal{G}_{t-1}\right] \geq{}& e_t(\check{x}_t, \hat{a}_t)_i^2 + \nabla_x e_t(\check{x}_t, \hat{a}_t)_i^\top \text{Cov}(\iota_t\mid \mathcal{G}_{t-1}) \nabla_x e_t(\check{x}_t, \hat{a}_t)_i \nonumber\\
    &- \overline{\epsilon}^2 \gamma_e \abs{e_t(\check{x}_t, \hat{a}_t)_i} - \overline{\epsilon}^3 \beta_e \gamma_e\\
    \geq{}& e_t(\check{x}_t, \hat{a}_t)_i^2 + \underline{\sigma} \norm{\nabla_x e_t(\check{x}_t, \hat{a}_t)_i}^2 - \overline{\epsilon}^2 \gamma_e \abs{e_t(\check{x}_t, \hat{a}_t)_i} - \overline{\epsilon}^3 \beta_e \gamma_e.
\end{align}
\end{subequations}
Summing over $t = 1, \ldots, T$ and taking expectation on both sides gives that
\begin{align}\label{thm:gradient-error-bound:e4}
    R_0^\ell(T) \geq \mathbb{E}\left[\sum_{t=1}^T e_t(\check{x}_t, \hat{a}_t)_i^2\right] + \underline{\sigma} \mathbb{E}\left[\sum_{t=1}^T \norm{\nabla_x e_t(\check{x}_t, \hat{a}_t)}_i^2\right] - \overline{\epsilon}^2 \gamma_e \mathbb{E}\left[\sum_{t=1}^T \abs{e_t(\check{x}_t, \hat{a}_t)_i}\right] - \overline{\epsilon}^3 \beta_e \gamma_e T.
\end{align}
Now we show that $\mathbb{E}\left[\sum_{t=1}^T e_t(\check{x}_t, \hat{a}_t)_i^2\right] \leq \overline{\epsilon}^2 T$. For the sake of contradiction, suppose 
\[\mathbb{E}\left[\sum_{t=1}^T e_t(\check{x}_t, \hat{a}_t)_i^2\right] > \overline{\epsilon}^2 T.\]
By \eqref{thm:gradient-error-bound:e4}, we see that
\begin{subequations}\label{thm:gradient-error-bound:e5}
\begin{align}
    R_0^\ell(T) \geq{}& \mathbb{E}\left[\sum_{t=1}^T e_t(\check{x}_t, \hat{a}_t)_i^2\right] - \overline{\epsilon}^2 \gamma_e \mathbb{E}\left[\sum_{t=1}^T \abs{e_t(\check{x}_t, \hat{a}_t)_i}\right] - \overline{\epsilon}^3 \beta_e \gamma_e T\nonumber\\
    \geq{}& \mathbb{E}\left[\sum_{t=1}^T e_t(\check{x}_t, \hat{a}_t)_i^2\right] - \overline{\epsilon}^2 \gamma_e \sqrt{\mathbb{E}\left[\left(\sum_{t=1}^T \abs{e_t(\check{x}_t, \hat{a}_t)_i}\right)^2\right]} - \overline{\epsilon}^3 \beta_e \gamma_e T \label{thm:gradient-error-bound:e5:s1}\\
    \geq{}& \mathbb{E}\left[\sum_{t=1}^T e_t(\check{x}_t, \hat{a}_t)_i^2\right] - \overline{\epsilon}^2 \gamma_e \sqrt{T \cdot \mathbb{E}\left[\sum_{t=1}^T e_t(\check{x}_t, \hat{a}_t)_i^2\right]} - \overline{\epsilon}^3 \beta_e \gamma_e T \label{thm:gradient-error-bound:e5:s2}\\
    ={}& \sqrt{\mathbb{E}\left[\sum_{t=1}^T e_t(\check{x}_t, \hat{a}_t)_i^2\right]} \cdot \left(\sqrt{\mathbb{E}\left[\sum_{t=1}^T e_t(\check{x}_t, \hat{a}_t)_i^2\right]} - \overline{\epsilon}^2 \gamma_e \sqrt{T}\right) -\overline{\epsilon}^3 \beta_e \gamma_e T \nonumber\\
    >{}& \frac{1}{4} \overline{\epsilon}^2 T, \label{thm:gradient-error-bound:e5:s3}
\end{align}
\end{subequations}
where we use Jensen's inequality in \eqref{thm:gradient-error-bound:e5:s1}; we use Cauchy-Schwarz inequality in \eqref{thm:gradient-error-bound:e5:s2}; we use the assumptions that $\overline{\epsilon} \gamma_e \leq \frac{1}{2}$ and $\overline{\epsilon} \beta_e \gamma_e \leq \frac{1}{4}$ in \eqref{thm:gradient-error-bound:e5:s3}. \eqref{thm:gradient-error-bound:e5} contradicts with our assumption that $R(T) \leq \bar{\epsilon}^3 T$. Thus, we have shown that $\mathbb{E}\left[\sum_{t=1}^T e_t(\check{x}_t, \hat{a}_t)_i^2\right] \leq \bar{\epsilon}^2 T$.

Using the same argument as \eqref{thm:gradient-error-bound:e5:s1} and \eqref{thm:gradient-error-bound:e5:s2}, we see that
\begin{align}\label{thm:gradient-error-bound:e6}
    \mathbb{E}\left[\sum_{t=1}^T \abs{e_t(\check{x}_t, \hat{a}_t)_i}\right] \leq \sqrt{T \cdot \mathbb{E}\left[\sum_{t=1}^T e_t(\check{x}_t, \hat{a}_t)_i^2\right]} \leq \bar{\epsilon} T.
\end{align}
Substituting \eqref{thm:gradient-error-bound:e6} into \eqref{thm:gradient-error-bound:e4} gives that
\begin{align*}
    \underline{\sigma} \mathbb{E}\left[\sum_{t=1}^T \norm{\nabla_x e_t(\check{x}_t, \hat{a}_t)_i}^2\right] \leq{}& R(T) + \overline{\epsilon}^2 \gamma_e \mathbb{E}\left[\sum_{t=1}^T \abs{e_t(\check{x}_t, \hat{a}_t)_i}\right] + \overline{\epsilon}^3 \beta_e \gamma_e T\\
    \leq{}& (1 + \gamma_e + \beta_e \gamma_e) \bar{\epsilon}^3 T.\\
\end{align*}
Therefore, we see that
\begin{align}\label{thm:gradient-error-bound:e7}
    &\mathbb{E}\left[\sum_{t=1}^T \norm{\nabla_x e_t(x_t, \hat{a}_t)_i}^2\right]\nonumber\\
    \leq{}& 2 \mathbb{E}\left[\sum_{t=1}^T \norm{\nabla_x e_t(\check{x}_t, \hat{a}_t)_i}^2\right] + 2 \mathbb{E}\left[\sum_{t=1}^T \norm{\nabla_x e_t(\check{x}_t, \hat{a}_t)_i - \nabla_x e_t(x_t, \hat{a}_t)_i}^2\right]\nonumber\\
    \leq{}& \frac{2}{\underline{\sigma}}(1 + \gamma_e + \beta_e \gamma_e) \bar{\epsilon}^3 T + 2 \gamma_e^2 \overline{\epsilon}^2 T.
\end{align}
Summing \eqref{thm:gradient-error-bound:e7} over dimensions $i \in [1:k]$ finishes the proof of \Cref{thm:gradient-error-bound}.

\section{Proof of Lemma \ref{lemma:fully-actuated-assumptions}}\label{appendix:fully-actuated-assumptions}
When applied to the dynamical system \eqref{equ:dynamics-fully-actuated} and the policy class \eqref{equ:policy-cancel-out}, the joint dynamics induced by applying M-GAPS with exact model parameters $a_{0:T-1}^*$ are given by
\begin{subequations}\label{equ:M-GAPS-dynamics-exact}
\begin{align}
    x_{t+1} ={}& q_t^x(x_t, y_t, \theta_t, a_t^*) = \phi_t(x_t, \psi_t(x_t, \theta_t)) + w_t,\\
    y_{t+1} ={}& q_t^y(x_t, y_t, \theta_t, a_t^*) = \left.\frac{\partial g^*_{t+1\mid t}}{\partial x_t}\right|_{x_t, \theta_t} \cdot y_t + \left.\frac{\partial g^*_{t+1\mid t}}{\partial \theta_t}\right|_{x_t, \theta_t},\\
    \theta_{t+1} ={}& q_t^\theta(x_t, y_t, \theta_t, a_t^*) = \Pi_\Theta\left(\theta_{t+1} - \eta \left(\left.\frac{\partial h_{t\mid t}^*}{\partial x_t}\right|_{x_t, \theta_t} \cdot y_t + \left.\frac{\partial h_{t\mid t}^*}{\partial \theta_t}\right|_{x_t, \theta_t}\right)\right).
\end{align}
\end{subequations}
The joint dynamics induced by applying M-GAPS with inexact parameters $\hat{a}_{0:T-1}$ are given by
\begin{subequations}\label{equ:M-GAPS-dynamics-inexact}
\begin{align}
    x_{t+1} ={}& q_t^x(x_t, y_t, \theta_t, \hat{a}_t) = \phi_t(x_t, \psi_t(x_t, \theta_t) + f_t(x_t, a_t^*) - f_t(x_t, \hat{a}_t)) + w_t,\\
    y_{t+1} ={}& q_t^y(x_t, y_t, \theta_t, \hat{a}_t) = \left.\frac{\partial g^*_{t+1\mid t}}{\partial x_t}\right|_{x_t, \theta_t} \cdot y_t + \left.\frac{\partial g^*_{t+1\mid t}}{\partial \theta_t}\right|_{x_t, \theta_t},\\
    \theta_{t+1} ={}& q_t^\theta(x_t, y_t, \theta_t, \hat{a}_t) = \Pi_\Theta\left(\theta_{t+1} - \eta \left(\left.\frac{\partial \hat{h}_{t\mid t}}{\partial x_t}\right|_{x_t, \theta_t} \cdot y_t + \left.\frac{\partial \hat{h}_{t\mid t}}{\partial \theta_t}\right|_{x_t, \theta_t}\right)\right),
\end{align}
\end{subequations}
where recall that we view $\hat{a}_{0:T-1}$ as external inputs as discussed
in \Cref{sec:meta-online-policy-selection}.

Since Lemma \ref{lemma:fully-actuated-assumptions} consists three properties, we show them separately in Lemmas \ref{lemma:matched-disturbance-Lipschitz-a}-\ref{lemma:matched-disturbance-robustness}.

\begin{lemma}\label{lemma:matched-disturbance-Lipschitz-a}
Consider the dynamical system
\begin{align*}
    x_{t+1} ={}& q_t^x(x_t, y_t, \theta_t, a_t^*) = \phi_t(x_t, \psi_t(x_t, \theta_t)) + w_t,\\
    y_{t+1} ={}& q_t^y(x_t, y_t, \theta_t, a_t^*) = \left.\frac{\partial g^*_{t+1\mid t}}{\partial x_t}\right|_{x_t, \theta_t} \cdot y_t + \left.\frac{\partial g^*_{t+1\mid t}}{\partial \theta_t}\right|_{x_t, \theta_t},\\
    \theta_{t+1} ={}& q_t^\theta(x_t, y_t, \theta_t, a_t^*) = \Pi_\Theta\left(\theta_{t+1} - \eta \left(\left.\frac{\partial h_{t\mid t}^*}{\partial x_t}\right|_{x_t, \theta_t} \cdot y_t + \left.\frac{\partial h_{t\mid t}^*}{\partial \theta_t}\right|_{x_t, \theta_t}\right)\right).
\end{align*}
For any $x_t, y_t, \theta_t, \hat{a}_t$ that satisfies
\[\norm{x_t} \leq R_x, \norm{y_t} \leq R_y, \theta_t \in \Theta, \hat{a}_t \in \mathcal{A},\]
the following Lipschitzness conditions hold:
\begin{align*}
    \norm{q_t^x(x_t, y_t, \theta_t, a_t^*) - q_t^x(x_t, y_t, \theta_t, \hat{a}_t)} &\leq \alpha_x \varepsilon_t(x_t, \hat{a}_t, a_t^*) + \beta_x \varepsilon_t'(x_t, \hat{a}_t, a_t^*),\\
    \norm{q_t^y(x_t, y_t, \theta_t, a_t^*) - q_t^y(x_t, y_t, \theta_t, \hat{a}_t)} &\leq \alpha_y \varepsilon_t(x_t, \hat{a}_t, a_t^*) + \beta_y \varepsilon_t'(x_t, \hat{a}_t, a_t^*),\\
    \norm{q_t^\theta(x_t, y_t, \theta_t, a_t^*) - q_t^\theta(x_t, y_t, \theta_t, \hat{a}_t)} &\leq \alpha_\theta \varepsilon_t(x_t, \hat{a}_t, a_t^*) + \beta_\theta \varepsilon_t'(x_t, \hat{a}_t, a_t^*),
\end{align*}
where
\begin{align*}
    \alpha_x ={}& \ell_{h, u} L_{\psi, \theta},\  \beta_x = \alpha_y = \beta_y = 0,\\
    \alpha_\theta ={}& \eta \left(R_y (\ell_{h, x} + \ell_{h, u} L_{f, x} + \ell_{h, u} L_{\psi, x}) + \ell_{h, u} L_{\psi, \theta}\right), \ \beta_\theta = \eta R_y L_{h, u}.
\end{align*}
Further, $q_t^\theta(x, y, \theta, a_t^*)$ is $(L_{\theta, x}, L_{\theta, y})$-Lipschitz in $(x, y)$, where
\begin{align*}
    L_{\theta, x} ={}& \eta R_y \left((\ell_{h, x} + \ell_{h, u} (L_{f, x} + L_{\psi, x})) (1 + L_{f, x} + L_{\psi, x}) + L_{h, u} (\ell_{f, x} + \ell_{\psi, x})\right)\\
    &+ \eta \left(\ell_{h, x} L_{\psi, \theta} + L_{h, u} \ell_{\psi, x} + \ell_{h, u} L_{\psi, \theta} (L_{f, x} + L_{\psi, x})\right),\\
    L_{\theta, y} ={}& \eta (L_{h, x} + L_{h, u} (L_{f, x} + L_{\psi, x})).
\end{align*}
\end{lemma}

The proof of Lemma \ref{lemma:matched-disturbance-Lipschitz-a} can be found in Appendix \ref{appendix:proof:lemma:matched-disturbance-Lipschitz-a}.

\begin{lemma}\label{lemma:matched-disturbance-contraction-stability}
Suppose the sequence $\theta_{0:T-1}$ is given and it satisfies the constraint that $\norm{\theta_t - \theta_{t-1}} \leq \epsilon_\theta$ for all time step $t$. Consider the dynamical system
\begin{align*}
    x_{t+1} ={}& q_t^x(x_t, y_t, \theta_t, a_t^*) = \phi_t(x_t, \psi_t(x_t, \theta_t)) + w_t,\\
    y_{t+1} ={}& q_t^y(x_t, y_t, \theta_t, a_t^*) = \left.\frac{\partial g^*_{t+1\mid t}}{\partial x_t}\right|_{x_t, \theta_t} \cdot y_t + \left.\frac{\partial g^*_{t+1\mid t}}{\partial \theta_t}\right|_{x_t, \theta_t}.
\end{align*}
%\soonjo{define $\rho\in [0,1)$.}
We have that $\norm{x_t} \leq R_x^* < R_x, \norm{y_t} \leq R_y^* < R_y$ always hold if the system starts from $(x_\tau, y_\tau) = (0, 0)$. Here,
\[R_x^* = R_S,\text{ and } R_y^* = \frac{C_{L, g, \theta}}{1 - \rho},\]
where recall that $\rho$ is the decay factor defined in Assumption \ref{assump:contractive-and-stability}.
Further, from any initial states $(x_\tau, y_\tau), (x_\tau', y_\tau')$ that satisfy $\norm{x_\tau}, \norm{x_\tau'} \leq R_x$ and $\norm{y_\tau}, \norm{y_\tau'} \leq R_y$, the trajectory satisfies
\[\norm{(x_t, y_t) - (x_t', y_t')} \leq \gamma(t-\tau) \cdot \norm{(x_\tau, y_\tau) - (x_\tau', y_\tau')},\]
where
\[\gamma(\tau) = \left(\bar{C} + C_{\ell, g, (x, x)} R_y + C_{\ell, g, (\theta, x)} \bar{C} \tau\right) \rho^{\tau}.\]
Note that $\gamma$ satisfies
\[\sum_{\tau = 0}^\infty \gamma(\tau) \leq C, \text{ where }C = \frac{\bar{C} + C_{\ell, g, (x, x)} R_y}{1 - \rho} + \frac{C_{\ell, g, (\theta, x)} \bar{C}}{(1 - \rho)^2}.\]
\end{lemma}

The definitions of the coefficients $C_{L, g, \theta}, C_{\ell, g, (x, x)}, C_{\ell, g, (\theta, x)}$ can be found in Lemma \ref{lemma:smooth-multi-step-dynamics} in Appendix \ref{appendix:useful-lemmas}. And the proof of Lemma \ref{lemma:matched-disturbance-contraction-stability} can be found in Appendix \ref{appendix:proof:lemma:matched-disturbance-contraction-stability}.\\

\begin{lemma}\label{lemma:matched-disturbance-robustness}
Consider the dynamical system
\begin{align}\label{equ:M-GAPS-biased-update}
    x_{t+1} ={}& q_t^x(x_t, y_t, \theta_t, a_t^*) = \phi_t(x_t, \psi_t(x_t, \theta_t)) + w_t,\nonumber\\
    y_{t+1} ={}& q_t^y(x_t, y_t, \theta_t, a_t^*) = \left.\frac{\partial g^*_{t+1\mid t}}{\partial x_t}\right|_{x_t, \theta_t} \cdot y_t + \left.\frac{\partial g^*_{t+1\mid t}}{\partial \theta_t}\right|_{x_t, \theta_t},\nonumber\\
    \theta_{t+1} ={}& q_t^\theta(x_t, y_t, \theta_t, a_t^*) = \Pi_\Theta\left(\theta_{t+1} - \eta \left(\left.\frac{\partial h_{t\mid t}^*}{\partial x_t}\right|_{x_t, \theta_t} \cdot y_t + \left.\frac{\partial h_{t\mid t}^*}{\partial \theta_t}\right|_{x_t, \theta_t}\right)\right) + \zeta_t.
\end{align}
Suppose the learning rate $\eta$ satisfies $\eta < \min\left\{\frac{(1 - \rho)\epsilon_\theta}{C_{L, h, \theta}}, \frac{1 - \rho}{2 C_{\ell, h, (\theta, \theta)}}\right\}$. When $\norm{\zeta_t} \leq \bar{\zeta} \coloneqq \min\{1, \epsilon_\theta - \frac{C_{L, h, \theta} \eta}{1 - \rho}\}$ holds for all $t$, the resulting $\{\theta_t\}$ satisfies the slowly-time-varying constraint $\norm{\theta_t - \theta_{t-1}} \leq \epsilon_\theta$. Further, the trajectory $\{\theta_t\}$ achieves the local regret guarantee
\begin{align*}
    R_\eta^L(T, \{\norm{\zeta_t}\}_{0\leq t\leq T-1}) \leq \frac{2}{\eta}(F_0(\theta_0) + S_0) + \frac{2}{1 - \rho}(C_{L, h, \theta} S_1 + C_{\ell, h, (\theta, \theta)} \eta S_2), \text{ where }
\end{align*}
\begin{align*}
    S_0 \coloneqq{}& \frac{2\bar{C} L_h (1 + L_{\psi, x} + L_{f, x}) (1 + L_{\phi, u})}{(1 - \rho)^2 \rho} \cdot (V_{sys} + V_w)\\
    &+ \frac{2 \bar{C} L_h (1 + L_{\psi, x} + L_{f, x})}{1 - \rho} \cdot \left(2\bar{C}\bar{R}_C + 2R_S\right),\\
    S_1 \coloneqq{}& \left(\frac{1}{\eta} + \frac{\hat{C}_3 + \hat{C}_5}{(1 - \rho)^2} + \frac{\hat{C}_4}{(1 - \rho)^3}\right) \sum_{t=0}^{T-1} \norm{\zeta_t} + \left(\frac{\hat{C}_0}{1 - \decayfactor} + \frac{\hat{C}_1 + \hat{C}_2}{(1 - \decayfactor)^{2}} + \frac{\hat{C}_2}{(1 - \decayfactor)^{3}}\right) \eta T,\\
    S_2 \coloneqq{}& \left(1 + \frac{\hat{C}_0}{1 - \decayfactor} + \frac{\hat{C}_1 + \hat{C}_2 + \hat{C}_3 + \hat{C}_5}{(1 - \decayfactor)^{2}} + \frac{\hat{C}_2 + \hat{C}_4}{(1 - \decayfactor)^{3}}\right)\cdot \\
    &\left[\left(\frac{1}{\eta^2} + \frac{\hat{C}_3 + \hat{C}_5}{(1 - \rho)^2} + \frac{\hat{C}_4}{(1 - \rho)^3}\right) \sum_{t=0}^{T-1} \norm{\zeta_t}^2 + \left(\frac{\hat{C}_0}{1 - \decayfactor} + \frac{\hat{C}_1 + \hat{C}_2}{(1 - \decayfactor)^{2}} + \frac{\hat{C}_2}{(1 - \decayfactor)^{3}}\right) \eta^2 T\right].
\end{align*}
Here, the variation intensity is defined as
\begin{align*}
    V_{\text{sys}} ={}& \sum_{t = 1}^{T-1}
    \Big(\sup_{x \in \mathcal{X}, u \in \mathcal{U}}
        \norm{\phi_t(x, u) - \phi_{t-1}(x, u)}
    + \sup_{x \in \mathcal{X}, \theta \in \Theta}
        \norm{\psi_t(x, \theta) - \psi_{t-1}(x, \theta)}\\
    &+ \sup_{x \in \mathcal{X}, u \in \mathcal{U}, \theta \in \Theta}
        \abs{h_t(x, u, \theta) - h_{t-1}(x, u, \theta)}\Big), \text{ and }\\
    V_{w} ={}& \sum_{t=1}^{T-1} \norm{w_t - w_{t-1}},
\end{align*}
The bound can be simplified to
\begin{align*}
    R_\eta^L(T, \{\norm{\zeta_t}\}_{0\leq t\leq T-1}) = O\left(\frac{1}{\eta} (1 + V_{\text{sys}} + V_w) + \eta T + \eta^3 T + \frac{1}{\eta}\sum_{t=1}^{T-1} \norm{\zeta_t}\right),
\end{align*}
where the $O(\cdot)$ notation hides dependence on $\frac{1}{1-\rho}, R_x, R_y, \bar{C},$ and the Lipchitzness/smoothness coefficients defined in Assumption \ref{assump:Lipschitz-and-smoothness}.
\end{lemma}

The definition of the coefficient $C_{L, h, \theta}$ can be found in Corollary \ref{coro:smooth-multi-step-costs} in Appendix \ref{appendix:useful-lemmas}. The proof of Lemma \ref{lemma:matched-disturbance-robustness} can be found in Appendix \ref{appendix:proof:lemma:matched-disturbance-robustness}.\\

\subsection{Proof of Lemma \ref{lemma:matched-disturbance-Lipschitz-a}}\label{appendix:proof:lemma:matched-disturbance-Lipschitz-a}
By Assumption \ref{assump:Lipschitz-and-smoothness}, we see that
\begin{align}\label{lemma:fully-actuated-assumptions:e1}
    \norm{q_t^x(x_t, y_t, \theta_t, \hat{a}_t) - q_t^x(x_t, y_t, \theta_t, a_t^*)} \leq L_{\phi, u} \norm{f_t(x_t, a_t^*) - f_t(x_t, \hat{a}_t)} = L_{\phi, u} \varepsilon_t.
\end{align}
We also have that
\begin{align}\label{lemma:fully-actuated-assumptions:e2}
    q_t^y(x_t, y_t, \theta_t, \hat{a}_t) = q_t^y(x_t, y_t, \theta_t, a_t^*).
\end{align}
Note that
\begin{align*}
    h_{t\mid t}^*(x_t, \theta_t) &= h_t(x_t, u_t^1, \theta_t), \text{ where }u_t^1 = - f_t(x_t, a_t^*) + \psi_t(x_t, \theta_t),\\
    \hat{h}_{t\mid t}(x_t, \theta_t) &= h_t(x_t, u_t^2, \theta_t), \text{ where }u_t^2 = - f_t(x_t, \hat{a}_t) + \psi_t(x_t, \theta_t).
\end{align*}
Therefore, we see that
\begin{subequations}\label{lemma:fully-actuated-assumptions:e3-0}
\begin{align}
    &\norm{\left.\frac{\partial h_{t\mid t}^*}{\partial x_t}\right|_{x_t, \theta_t} - \left.\frac{\partial \hat{h}_{t\mid t}}{\partial x_t}\right|_{x_t, \theta_t}}\nonumber\\
    \leq{}& \norm{\left.\frac{\partial h_t}{\partial x_t}\right|_{x_t, u_t^1, \theta_t} - \left.\frac{\partial h_t}{\partial x_t}\right|_{x_t, u_t^2, \theta_t}} + \norm{\left.\frac{\partial h_t}{\partial u_t}\right|_{x_t, u_t^1, \theta_t}\cdot \left.\frac{\partial f_t}{\partial x_t}\right|_{x_t, a_t^*} - \left.\frac{\partial h_t}{\partial u_t}\right|_{x_t, u_t^2, \theta_t}\cdot \left.\frac{\partial f_t}{\partial x_t}\right|_{x_t, \hat{a}_t}}\nonumber\\
    &+ \norm{\left.\frac{\partial h_t}{\partial u_t}\right|_{x_t, u_t^1, \theta_t}\cdot \left.\frac{\partial \psi_t}{\partial x_t}\right|_{x_t, \theta_t} - \left.\frac{\partial h_t}{\partial u_t}\right|_{x_t, u_t^2, \theta_t}\cdot \left.\frac{\partial \psi_t}{\partial x_t}\right|_{x_t, \theta_t}}\label{lemma:fully-actuated-assumptions:e3-0:s1}\\
    \leq{}& \ell_{h, x} \varepsilon_t + (\ell_{h, u} L_{f, x}\varepsilon_t + L_{h, u} \varepsilon_t') + \ell_{h, u} L_{\psi, x} \varepsilon_t\label{lemma:fully-actuated-assumptions:e3-0:s2}\\
    ={}& (\ell_{h, x} + \ell_{h, u} L_{f, x} + \ell_{h, u} L_{\psi, x}) \varepsilon_t + L_{h, u} \varepsilon_t',\nonumber
\end{align}
\end{subequations}
where we use the chain rule and the triangle inequality in \eqref{lemma:fully-actuated-assumptions:e3-0:s1}; we use Assumption \ref{assump:Lipschitz-and-smoothness} in \eqref{lemma:fully-actuated-assumptions:e3-0:s2}. Similarly, we also see that
\begin{subequations}\label{lemma:fully-actuated-assumptions:e3-1}
\begin{align}
    \norm{\left.\frac{\partial h_{t\mid t}^*}{\partial \theta_t}\right|_{x_t, \theta_t} - \left.\frac{\partial \hat{h}_{t\mid t}}{\partial \theta_t}\right|_{x_t, \theta_t}}
    ={}& \norm{\left.\frac{\partial h_t}{\partial u_t}\right|_{x_t, u_t^1, \theta_t}\cdot \left.\frac{\partial \psi_t}{\partial \theta_t}\right|_{x_t, \theta_t} - \left.\frac{\partial h_t}{\partial u_t}\right|_{x_t, u_t^2, \theta_t}\cdot \left.\frac{\partial \psi_t}{\partial \theta_t}\right|_{x_t, \theta_t}} \label{lemma:fully-actuated-assumptions:e3-1:s1}\\
    \leq{}& \ell_{h, u} L_{\psi, \theta} \varepsilon_t, \label{lemma:fully-actuated-assumptions:e3-1:s2}
\end{align}
\end{subequations}
where we use the chain rule in \eqref{lemma:fully-actuated-assumptions:e3-1:s1} and Assumption \ref{assump:Lipschitz-and-smoothness} in \eqref{lemma:fully-actuated-assumptions:e3-1:s2}.

For $q_t^\theta$, we see that
\begin{subequations}\label{lemma:fully-actuated-assumptions:e3}
\begin{align}
    &\norm{q_t^\theta(x_t, y_t, \theta_t, \hat{a}_t) - q_t^\theta(x_t, y_t, \theta_t, a_t^*)}\nonumber\\
    \leq{}& \eta \norm{\left(\left.\frac{\partial h_{t\mid t}^*}{\partial x_t}\right|_{x_t, \theta_t} - \left.\frac{\partial \hat{h}_{t\mid t}}{\partial x_t}\right|_{x_t, \theta_t}\right) \cdot y_t + \left(\left.\frac{\partial h_{t\mid t}^*}{\partial \theta_t}\right|_{x_t, \theta_t} - \left.\frac{\partial \hat{h}_{t\mid t}}{\partial \theta_t}\right|_{x_t, \theta_t}\right)}\label{lemma:fully-actuated-assumptions:e3:s1}\\
    \leq{}& \eta \left(R_y (\ell_{h, x} + \ell_{h, u} L_{f, x} + \ell_{h, u} L_{\psi, x}) + \ell_{h, u} L_{\psi, \theta}\right) \varepsilon_t + \eta R_y L_{h, u} \varepsilon_t', \label{lemma:fully-actuated-assumptions:e3:s2}
\end{align}
\end{subequations}
where we use the property that projection onto $\Theta$ is contractive in \eqref{lemma:fully-actuated-assumptions:e3:s1}; we use \eqref{lemma:fully-actuated-assumptions:e3-0} and \eqref{lemma:fully-actuated-assumptions:e3-1} in \eqref{lemma:fully-actuated-assumptions:e3:s2}. We also see that
\begin{subequations}\label{lemma:fully-actuated-assumptions:e4}
\begin{align}
    &\norm{q_t^\theta(x_t, y_t, \theta_t, a_t^*) - q_t^\theta(x_t', y_t', \theta_t, a_t^*)}\nonumber\\
    \leq{}& \eta \norm{\left.\frac{\partial h_{t\mid t}^*}{\partial x_t}\right|_{x_t, \theta_t} \cdot y_t - \left.\frac{\partial h_{t\mid t}^*}{\partial x_t}\right|_{x_t', \theta_t} \cdot y_t' + \left.\frac{\partial h_{t\mid t}^*}{\partial \theta_t}\right|_{x_t, \theta_t} - \left.\frac{\partial h_{t\mid t}^*}{\partial \theta_t}\right|_{x_t', \theta_t}}\label{lemma:fully-actuated-assumptions:e4:s1}\\
    \leq{}& \eta \norm{\left.\frac{\partial h_{t\mid t}^*}{\partial x_t}\right|_{x_t, \theta_t} - \left.\frac{\partial h_{t\mid t}^*}{\partial x_t}\right|_{x_t', \theta_t}} \cdot \norm{y_t} + \eta \norm{\left.\frac{\partial h_{t\mid t}^*}{\partial x_t}\right|_{x_t', \theta_t}} \cdot \norm{y_t - y_t'}\nonumber\\
    &+ \eta \norm{\left.\frac{\partial h_{t\mid t}^*}{\partial \theta_t}\right|_{x_t, \theta_t} - \left.\frac{\partial h_{t\mid t}^*}{\partial \theta_t}\right|_{x_t', \theta_t}},\label{lemma:fully-actuated-assumptions:e4:s2}
\end{align}
\end{subequations}
where we use the property that projection onto $\Theta$ is contractive in \eqref{lemma:fully-actuated-assumptions:e4:s1}, and apply the triangle inequality in \eqref{lemma:fully-actuated-assumptions:e4:s2}. Note that
\begin{align*}
    h_{t\mid t}^*(x_t, \theta_t) &= h_t(x_t, u_t, \theta_t), \text{ where }u_t = - f_t(x_t, a_t^*) + \psi_t(x_t, \theta_t),\\
    h_{t\mid t}^*(x_t', \theta_t) &= h_t(x_t', u_t', \theta_t), \text{ where }u_t' = - f_t(x_t', a_t^*) + \psi_t(x_t', \theta_t).
\end{align*}
Therefore, we see that
\begin{subequations}\label{lemma:fully-actuated-assumptions:e4-0}
\begin{align}
    &\norm{\left.\frac{\partial h_{t\mid t}^*}{\partial x_t}\right|_{x_t, \theta_t} - \left.\frac{\partial h_{t\mid t}^*}{\partial x_t}\right|_{x_t', \theta_t}}\nonumber\\
    \leq{}& \norm{\left.\frac{\partial h_t}{\partial x_t}\right|_{x_t, u_t, \theta_t} - \left.\frac{\partial h_t}{\partial x_t}\right|_{x_t', u_t', \theta_t}} + \norm{\left.\frac{\partial h_t}{\partial u_t}\right|_{x_t, u_t, \theta_t}\cdot \left.\frac{\partial f_t}{\partial x_t}\right|_{x_t, a_t^*} - \left.\frac{\partial h_t}{\partial u_t}\right|_{x_t', u_t', \theta_t}\cdot \left.\frac{\partial f_t}{\partial x_t}\right|_{x_t', a_t^*}}\nonumber\\
    &+ \norm{\left.\frac{\partial h_t}{\partial u_t}\right|_{x_t, u_t, \theta_t}\cdot \left.\frac{\partial \psi_t}{\partial x_t}\right|_{x_t, \theta_t} - \left.\frac{\partial h_t}{\partial u_t}\right|_{x_t', u_t', \theta_t}\cdot \left.\frac{\partial \psi_t}{\partial x_t}\right|_{x_t', \theta_t}}\label{lemma:fully-actuated-assumptions:e4-0:s1}\\
    \leq{}& \norm{\left.\frac{\partial h_t}{\partial x_t}\right|_{x_t, u_t, \theta_t} - \left.\frac{\partial h_t}{\partial x_t}\right|_{x_t', u_t', \theta_t}} + \norm{\left.\frac{\partial h_t}{\partial u_t}\right|_{x_t, u_t, \theta_t} - \left.\frac{\partial h_t}{\partial u_t}\right|_{x_t', u_t', \theta_t}}\cdot \norm{\left.\frac{\partial f_t}{\partial x_t}\right|_{x_t, a_t^*}}\nonumber\\
    &+ \norm{\left.\frac{\partial h_t}{\partial u_t}\right|_{x_t', u_t', \theta_t}} \cdot \norm{\left.\frac{\partial f_t}{\partial x_t}\right|_{x_t, a_t^*} - \left.\frac{\partial f_t}{\partial x_t}\right|_{x_t', a_t^*}} + \norm{\left.\frac{\partial h_t}{\partial u_t}\right|_{x_t, u_t, \theta_t} - \left.\frac{\partial h_t}{\partial u_t}\right|_{x_t', u_t', \theta_t}}\cdot \norm{\left.\frac{\partial \psi_t}{\partial x_t}\right|_{x_t, \theta_t}}\nonumber\\
    &+ \norm{\left.\frac{\partial h_t}{\partial u_t}\right|_{x_t', u_t', \theta_t}} \cdot \norm{\left.\frac{\partial \psi_t}{\partial x_t}\right|_{x_t, \theta_t} - \left.\frac{\partial \psi_t}{\partial x_t}\right|_{x_t', \theta_t}}\label{lemma:fully-actuated-assumptions:e4-0:s2}\\
    \leq{}& \ell_{h, x} \norm{x_t - x_t'} + \ell_{h, u} \norm{u_t - u_t'} + L_{f, x} \left(\ell_{h, x} \norm{x_t - x_t'} + \ell_{h, u} \norm{u_t - u_t'}\right)\nonumber\\
    & + L_{h, u} \ell_{f, x} \norm{x_t - x_t'} + L_{\psi, x} \left(\ell_{h, x}\norm{x_t - x_t'} + \ell_{h, u}\norm{u_t - u_t'}\right) + L_{h, u} \ell_{\psi, x} \norm{x_t - x_t'}\label{lemma:fully-actuated-assumptions:e4-0:s3}\\
    ={}& \left(\ell_{h, x} (1 + L_{f, x} + L_{\psi, x}) + L_{h, u} (\ell_{f, x} + \ell_{\psi, x})\right)\norm{x_t - x_t'} + \ell_{h, u} (1 + L_{f, x} + L_{\psi, x})\norm{u_t - u_t'}, \nonumber\\
    \leq{}& \left((\ell_{h, x} + \ell_{h, u} (L_{f, x} + L_{\psi, x})) (1 + L_{f, x} + L_{\psi, x}) + L_{h, u} (\ell_{f, x} + \ell_{\psi, x})\right)\norm{x_t - x_t'},\label{lemma:fully-actuated-assumptions:e4-0:s4}
\end{align}
\end{subequations}
where we use the triangle inequality in \eqref{lemma:fully-actuated-assumptions:e4-0:s1} and \eqref{lemma:fully-actuated-assumptions:e4-0:s2}; we use Assumption \ref{assump:Lipschitz-and-smoothness} in \eqref{lemma:fully-actuated-assumptions:e4-0:s3} and \eqref{lemma:fully-actuated-assumptions:e4-0:s4}. Similarly, we also see that
\begin{subequations}\label{lemma:fully-actuated-assumptions:e4-1}
\begin{align}
    &\norm{\left.\frac{\partial h_{t\mid t}^*}{\partial \theta_t}\right|_{x_t, \theta_t} - \left.\frac{\partial h_{t\mid t}^*}{\partial \theta_t}\right|_{x_t', \theta_t}}\nonumber\\
    ={}& \norm{\left.\frac{\partial h_t}{\partial u_t}\right|_{x_t, u_t, \theta_t}\cdot \left.\frac{\partial \psi_t}{\partial \theta_t}\right|_{x_t, \theta_t} - \left.\frac{\partial h_t}{\partial u_t}\right|_{x_t', u_t', \theta_t}\cdot \left.\frac{\partial \psi_t}{\partial \theta_t}\right|_{x_t', \theta_t}}\label{lemma:fully-actuated-assumptions:e4-1:s1}\\
    \leq{}& \norm{\left.\frac{\partial h_t}{\partial u_t}\right|_{x_t, u_t, \theta_t} - \left.\frac{\partial h_t}{\partial u_t}\right|_{x_t', u_t', \theta_t}} \cdot \norm{\left.\frac{\partial \psi_t}{\partial \theta_t}\right|_{x_t, \theta_t}} + \norm{\left.\frac{\partial h_t}{\partial u_t}\right|_{x_t', u_t', \theta_t}} \cdot \norm{\left.\frac{\partial \psi_t}{\partial \theta_t}\right|_{x_t, \theta_t} - \left.\frac{\partial \psi_t}{\partial \theta_t}\right|_{x_t', \theta_t}}\label{lemma:fully-actuated-assumptions:e4-1:s2}\\
    \leq{}& \left(\ell_{h, x} L_{\psi, \theta} + L_{h, u} \ell_{\psi, x}\right) \norm{x_t - x_t'} + \ell_{h, u} L_{\psi, \theta}\norm{u_t - u_t'},\label{lemma:fully-actuated-assumptions:e4-1:s3}\\
    \leq{}& \left(\ell_{h, x} L_{\psi, \theta} + L_{h, u} \ell_{\psi, x} + \ell_{h, u} L_{\psi, \theta} (L_{f, x} + L_{\psi, x})\right)\norm{x_t - x_t'}, \label{lemma:fully-actuated-assumptions:e4-1:s4}
\end{align}
\end{subequations}
where we use the chain rule in \eqref{lemma:fully-actuated-assumptions:e4-1:s1}; we use the triangle inequality in \eqref{lemma:fully-actuated-assumptions:e4-1:s2}; we use Assumption \ref{assump:Lipschitz-and-smoothness} in \eqref{lemma:fully-actuated-assumptions:e4-1:s3}. Substituting \eqref{lemma:fully-actuated-assumptions:e4-0} and \eqref{lemma:fully-actuated-assumptions:e4-1} into \eqref{lemma:fully-actuated-assumptions:e4} gives that
\begin{align}\label{lemma:fully-actuated-assumptions:e5}
    &\norm{q_t^\theta(x_t, y_t, \theta_t, a_t^*) - q_t^\theta(x_t', y_t', \theta_t, a_t^*)}\nonumber\\
    \leq{}& \eta R_y \left((\ell_{h, x} + \ell_{h, u} (L_{f, x} + L_{\psi, x})) (1 + L_{f, x} + L_{\psi, x}) + L_{h, u} (\ell_{f, x} + \ell_{\psi, x})\right)\norm{x_t - x_t'}\nonumber\\
    &+ \eta (L_{h, x} + L_{h, u} (L_{f, x} + L_{\psi, x})) \norm{y_t - y_t'}\nonumber\\
    &+ \eta \left(\ell_{h, x} L_{\psi, \theta} + L_{h, u} \ell_{\psi, x} + \ell_{h, u} L_{\psi, \theta} (L_{f, x} + L_{\psi, x})\right)\norm{x_t - x_t'}\nonumber\\
    \leq{}& L_{\theta, x} \norm{x_t - x_t'} + L_{\theta, y} \norm{y_t - y_t'}.
\end{align}

\subsection{Proof of Lemma \ref{lemma:matched-disturbance-contraction-stability}}\label{appendix:proof:lemma:matched-disturbance-contraction-stability}
Consider two trajectories $\{x_{t_1:t_2}, y_{t_1:t_2}\}$ and $\{x_{t_1:t_2}', y_{t_1:t_2}'\}$ given by
\begin{align*}
    x_{\tau+1} ={}& \phi_\tau(x_\tau, \psi_t(x_\tau, \theta_\tau)) + w_\tau,\\
    y_{\tau+1} ={}& \left.\frac{\partial g^*_{\tau+1\mid \tau}}{\partial x_\tau}\right|_{x_\tau, \theta_\tau} \cdot y_\tau + \left.\frac{\partial g^*_{\tau+1\mid \tau}}{\partial \theta_\tau}\right|_{x_\tau, \theta_\tau},
\end{align*}
and
\begin{align*}
    x_{\tau+1}' ={}& \phi_\tau(x_\tau', \psi_t(x_\tau', \theta_\tau)) + w_\tau,\\
    y_{\tau+1}' ={}& \left.\frac{\partial g^*_{\tau+1\mid \tau}}{\partial x_\tau}\right|_{x_\tau', \theta_\tau} \cdot y_\tau' + \left.\frac{\partial g^*_{\tau+1\mid \tau}}{\partial \theta_\tau}\right|_{x_\tau', \theta_\tau},
\end{align*}
where $\tau = t_1, t_1 + 1, \ldots, t_2$.
Note that by Assumption \ref{assump:contractive-and-stability}, we have that $\norm{x_{t_2}} \leq R_S$ and for any $x_{t_1}, x_{t_1}'$ whose norms are upper bounded by $R_C$
\begin{align}\label{lemma:fully-actuated-contraction:e1}
    \norm{x_{t_2} - x_{t_2}'} \leq \bar{C} \rho^{t_2 - t_1} \norm{x_{t_1} - x_{t_1}'}.
\end{align}
where $\rho$ is the decay factor of the contractive perturbation property defined in Assumption \ref{assump:contractive-and-stability}. For the $y$ sequence, note that $y_{t_2}$ and $y_{t_2}'$ can be expressed equivalently as
\begin{subequations}\label{lemma:fully-actuated-contraction:e2}
\begin{align}
    y_{t_2} ={}& \left.\frac{\partial g_{t_2\mid t_1}^*}{\partial x_{t_1}}\right|_{x_{t_1}, \theta_{t_1:t_2 - 1}} \cdot y_{t_1} + \sum_{\tau = t_1}^{t_2-1} \left.\frac{\partial g_{t_2\mid \tau}^*}{\partial \theta_\tau}\right|_{x_\tau, \theta_{\tau:t_2 - 1}},\label{lemma:fully-actuated-contraction:e2:s1}\\
    y_{t_2}' ={}& \left.\frac{\partial g_{t_2\mid t_1}^*}{\partial x_{t_1}}\right|_{x_{t_1}', \theta_{t_1:t_2 - 1}} \cdot y_{t_1}' + \sum_{\tau = t_1}^{t_2-1} \left.\frac{\partial g_{t_2\mid \tau}^*}{\partial \theta_\tau}\right|_{x_\tau', \theta_{\tau:t_2 - 1}}.\label{lemma:fully-actuated-contraction:e2:s2}
\end{align}
\end{subequations}
By Lemma \ref{lemma:smooth-multi-step-dynamics}, we see that if $y_{t_1} = 0$, then
\begin{align}\label{lemma:fully-actuated-contraction:e3}
    \norm{y_{t_2}} ={}& \norm{\sum_{\tau = t_1}^{t_2-1} \left.\frac{\partial g_{t_2\mid \tau}^*}{\partial \theta_\tau}\right|_{x_\tau, \theta_{\tau:t_2 - 1}}} \leq \sum_{\tau = t_1}^{t_2-1} \norm{\left.\frac{\partial g_{t_2\mid \tau}^*}{\partial \theta_\tau}\right|_{x_\tau, \theta_{\tau:t_2 - 1}}} \leq \sum_{\tau = t_1}^{t_2-1} C_{L, g, \theta} \rho^{t_2 - \tau} = \frac{C_{L, g, \theta}}{1 - \rho}.
\end{align}
We also see that
\begin{subequations}\label{lemma:fully-actuated-contraction:e4}
\begin{align}
    &\norm{y_{t_2} - y_{t_2}'}\nonumber\\
    ={}& \norm{\left(\left.\frac{\partial g_{t_2\mid t_1}^*}{\partial x_{t_1}}\right|_{x_{t_1}, \theta_{t_1:t_2 - 1}} - \left.\frac{\partial g_{t_2\mid t_1}^*}{\partial x_{t_1}}\right|_{x_{t_1}', \theta_{t_1:t_2 - 1}}\right) \cdot y_{t_1}} + \norm{\left.\frac{\partial g_{t_2\mid t_1}^*}{\partial x_{t_1}}\right|_{x_{t_1}', \theta_{t_1:t_2 - 1}} \cdot (y_{t_1} - y_{t_1}')}\nonumber\\
    &+ \sum_{\tau = t_1}^{t_2-1} \norm{\left.\frac{\partial g_{t_2\mid \tau}^*}{\partial \theta_\tau}\right|_{x_\tau, \theta_{\tau:t_2 - 1}} - \left.\frac{\partial g_{t_2\mid \tau}^*}{\partial \theta_\tau}\right|_{x_\tau', \theta_{\tau:t_2 - 1}}}\\
    \leq{}& C_{\ell, g, (x, x)} \rho^{t_2 - t_1} \norm{x_{t_1} - x_{t_1}'} \cdot R_y + C_{L, g, x} \rho^{t_2 - t_1} \norm{y_{t_1} - y_{t_1}'}\nonumber\\
    &+ C_{\ell, g, (\theta, x)} \sum_{\tau = t_1}^{t_2-1} \rho^{t_2 - \tau} \norm{x_\tau - x_\tau'}\\
    \leq{}& C_{\ell, g, (x, x)} \rho^{t_2 - t_1} \norm{x_{t_1} - x_{t_1}'} \cdot R_y + C_{L, g, x} \rho^{t_2 - t_1} \norm{y_{t_1} - y_{t_1}'}\nonumber\\
    &+ C_{\ell, g, (\theta, x)} \sum_{\tau = t_1}^{t_2-1} \rho^{t_2 - \tau} \cdot \bar{C} \rho^{\tau - t_1} \norm{x_{t_1} - x_{t_1}'}\\
    \leq{}& \left(C_{\ell, g, (x, x)} R_y + C_{\ell, g, (\theta, x)} \bar{C} (t_2 - t_1)\right) \rho^{t_2 - t_1} \norm{x_{t_1} - x_{t_1}'} + C_{L, g, x} \rho^{t_2 - t_1} \norm{y_{t_1} - y_{t_1}'}.
\end{align}
\end{subequations}
Therefore, we see that
\begin{subequations}\label{lemma:fully-actuated-contraction:e5}
\begin{align}
    &\norm{(x_{t_2}, y_{t_2}) - (x_{t_2}', y_{t_2}')}\nonumber\\
    \leq{}& \norm{x_{t_2} - x_{t_2}'} + \norm{y_{t_2} - y_{t_2}'} \label{lemma:fully-actuated-contraction:e5:s1}\\
    \leq{}& \bar{C} \rho^{t_2 - t_1} \norm{x_{t_1} - x_{t_1}'} + \left(C_{\ell, g, (x, x)} R_y + C_{\ell, g, (\theta, x)} \bar{C} (t_2 - t_1)\right) \rho^{t_2 - t_1} \norm{x_{t_1} - x_{t_1}'}\nonumber\\
    &+ \bar{C} \rho^{t_2 - t_1} \norm{y_{t_1} - y_{t_1}'} \label{lemma:fully-actuated-contraction:e5:s2}\\
    \leq{}& \gamma(t_2 - t_1) \norm{(x_{t_1}, y_{t_1}) - (x_{t_1}', y_{t_1}')}, \label{lemma:fully-actuated-contraction:e5:s3}
\end{align}
\end{subequations}
where we use the triangle inequality in \eqref{lemma:fully-actuated-contraction:e5:s1}; we use \eqref{lemma:fully-actuated-contraction:e1} and \eqref{lemma:fully-actuated-contraction:e4} and $\bar{C} = C_{L, g, x}$ in \eqref{lemma:fully-actuated-contraction:e5:s2}; we use the inequality that
\[\norm{x_{t_1} - x_{t_1}'} + \norm{y_{t_1} - y_{t_1}'} \leq \sqrt{2} \norm{(x_{t_1}, y_{t_1}) - (x_{t_1}', y_{t_1}')}\]
and the definition of $\gamma(\cdot)$ in \eqref{lemma:fully-actuated-contraction:e5:s3}.

\subsection{Proof of Lemma \ref{lemma:matched-disturbance-robustness}}\label{appendix:proof:lemma:matched-disturbance-robustness}

We compare the dynamical system \eqref{equ:M-GAPS-biased-update} with the Ideal OGD update rule:
\begin{align}\label{lemma:matched-disturbance-robustness:e1}
    \theta_{t+1} = \Pi_\Theta(\theta_t - \eta \nabla F_t(\theta_t)).
\end{align}
Note that the update on $\theta_t$ that the dynamical system \eqref{equ:M-GAPS-biased-update} performs can be written equivalently as
\begin{align}\label{lemma:matched-disturbance-robustness:e2}
    \theta_{t+1} = \Pi_\Theta(\theta_t - \eta G_t) + \zeta_t,
\end{align}
where
\begin{align}\label{lemma:matched-disturbance-robustness:e3}
    G_t \coloneqq \sum_{\tau = 0}^{t} \left.\frac{\partial h_{t\mid 0}^*}{\partial \theta_{t - \tau}}\right|_{x_0, \theta_{0:t}}.
\end{align}
By \Cref{thm:OCO-with-parameter-update-gradient-bias:full}, we know that
\begin{align*}
    \norm{G_t - \nabla F_t(\theta_t)} \leq{}& \left(\frac{\hat{C}_0}{1 - \decayfactor} + \frac{\hat{C}_1 + \hat{C}_2}{(1 - \decayfactor)^{2}} + \frac{\hat{C}_2}{(1 - \decayfactor)^{3}}\right) \eta\\
    &+ \sum_{\tau = 0}^{t-1}\left(\frac{\hat{C}_3}{1 - \rho}
    + \frac{\hat{C}_4}{(1 - \rho)^2} + \hat{C}_5 (t - \tau)\right) \rho^{t-\tau} \norm{\zeta_\tau},
\end{align*}
where the constants $\hat{C}_{0:5}$ are given in \Cref{thm:OCO-with-parameter-update-gradient-bias:full}. Let $\theta_{t+1}$ be the actual next policy parameter (following the update rule \eqref{lemma:matched-disturbance-robustness:e2}). By Lemma \ref{lemma:projection-decrease-distance}, we see that
\begin{align*}
    \norm{\theta_{t+1} - \Pi_\Theta(\theta_t - \eta \nabla F_t(\theta_t))} \leq{}& \eta \norm{G_t - \nabla F_t(\theta_t)} + \norm{\zeta_t}\\
    \leq{}& \norm{\zeta_t} + \left(\frac{\hat{C}_0}{1 - \decayfactor} + \frac{\hat{C}_1 + \hat{C}_2}{(1 - \decayfactor)^{2}} + \frac{\hat{C}_2}{(1 - \decayfactor)^{3}}\right) \eta^2\\
    &+ \eta \sum_{\tau = 0}^{t-1}\left(\frac{\hat{C}_3}{1 - \rho}
    + \frac{\hat{C}_4}{(1 - \rho)^2} + \hat{C}_5 (t - \tau)\right) \rho^{t-\tau} \norm{\zeta_\tau}.
\end{align*}
Then, we can apply \Cref{thm:nonconvex-biased-OGD} to obtain that
\begin{align}\label{lemma:matched-disturbance-robustness:e4}
    \sum_{t = 0}^{T-1}\norm{\nabla_{\Theta, \eta} F_t(\theta_t)}^2 \leq \frac{1}{\eta (1 - \eta \ell_F)}\left(F_0(\theta_0) + \sum_{t=1}^{T-1}\dist_s(F_t, F_{t-1})\right) + \frac{L_F S_1 + \ell_F \eta S_2}{1 - \eta \ell_F},
\end{align}
%\soonjo{Did you define $\dist_s$?}
where $\dist_S$ is a metric that measures the distance between two surrogate cost functions (see \Cref{thm:nonconvex-biased-OGD} for definition), and $S_1$ and $S_2$ are given by
\begin{align*}
    S_1 \coloneqq{}& \left(\frac{1}{\eta} + \frac{\hat{C}_3 + \hat{C}_5}{(1 - \rho)^2} + \frac{\hat{C}_4}{(1 - \rho)^3}\right) \sum_{t=0}^{T-1} \norm{\zeta_t} + \left(\frac{\hat{C}_0}{1 - \decayfactor} + \frac{\hat{C}_1 + \hat{C}_2}{(1 - \decayfactor)^{2}} + \frac{\hat{C}_2}{(1 - \decayfactor)^{3}}\right) \eta T,\\
    S_2 \coloneqq{}& \left(1 + \frac{\hat{C}_0}{1 - \decayfactor} + \frac{\hat{C}_1 + \hat{C}_2 + \hat{C}_3 + \hat{C}_5}{(1 - \decayfactor)^{2}} + \frac{\hat{C}_2 + \hat{C}_4}{(1 - \decayfactor)^{3}}\right)\cdot \\
    &\left[\left(\frac{1}{\eta^2} + \frac{\hat{C}_3 + \hat{C}_5}{(1 - \rho)^2} + \frac{\hat{C}_4}{(1 - \rho)^3}\right) \sum_{t=0}^{T-1} \norm{\zeta_t}^2 + \left(\frac{\hat{C}_0}{1 - \decayfactor} + \frac{\hat{C}_1 + \hat{C}_2}{(1 - \decayfactor)^{2}} + \frac{\hat{C}_2}{(1 - \decayfactor)^{3}}\right) \eta^2 T\right].\\
\end{align*}
By applying Lemma F.4 in \cite{lin2023online}, we can bound the total variational intensity on the surrogate costs by
\begin{align*}
    \sum_{t=1}^{T-1} \dist_s(F_t, F_{t-1}) \leq{}& \frac{2\bar{C} L_h (1 + L_{\psi, x} + L_{f, x}) (1 + L_{\phi, u})}{(1 - \rho)^2 \rho} \cdot (V_{sys} + V_w)\\
    &+ \frac{2 \bar{C} L_h (1 + L_{\psi, x} + L_{f, x})}{1 - \rho} \cdot \left(2\bar{C}\bar{R}_C + 2R_S\right).
\end{align*}
Substituting the above inequality and $L_F = \frac{C_{L, f, \theta}}{1 - \rho}, \ell_F = \frac{C_{\ell, h, (\theta, \theta)}}{1 - \rho}$ into \eqref{lemma:matched-disturbance-robustness:e4} finishes the proof.

\section{Proof of Lemma \ref{lemma:gradient-estimator-regret}}\label{appendix:grad-est}
%\subsection{Proof of Lemma \ref{lemma:gradient-estimator-regret}}
\label{appendix:gradient-estimator-regret}

By Assumptions \ref{assump:linear-function-approximation} and \ref{assump:estimator-environment}, we see that for any $a \in \mathcal{A}$,
\begin{align*}
    \tilde{\ell}_t(x_t, a, \tilde{f}_t) \leq{}& \left(f_t(x_t, a) - \tilde{f}_t\right)^2 = \left((f_t(x_t, a) - f_t(x_t, a_t^*)) - (\tilde{f}_t - f_t(x_t, a_t^*))\right)^2 \leq (C_f + e_f)^2.
\end{align*}
We also see that
\begin{align*}
    \norm{\nabla_a \tilde{\ell}_t(x_t, a, \tilde{f}_t)} \leq 2 \norm{\nabla_a f_t(x_t, a)} \cdot \abs{f_t(x_t, a) - \tilde{f}_t} \leq 2 D_f' (C_f + e_f).
\end{align*}
By Theorem 10.1 in \cite{hazan2016introduction}, we know that \Cref{alg:grad-est} with the learning rate $\iota = \frac{C_f + e_f}{D_f'}\cdot \sqrt{\frac{C_p}{T}}$ always achieves the guarantee that
\begin{align}\label{lemma:gradient-estimator-regret:e1}
    \sum_{t=0}^{T-1} \tilde{\ell}_t(x_t, \hat{a}_t, \tilde{f}_t) - \sum_{t=0}^{T-1} \tilde{\ell}_t(x_t, a_t^*, \tilde{f}_t) \leq R_0^\ell(T) \coloneqq 2\sqrt{3} (C_f + e_f)^3 D_f' \sqrt{C_p T}.
\end{align}
Let $v_t \coloneqq \tilde{f}_t - f_t(x_t, a_t^*)$. We see that
\begin{align*}
    &\mathbb{E}\left[\tilde{\ell}_t(x_t, \hat{a}_t, \tilde{f}_t) - \tilde{\ell}_t(x_t, a_t^*, \tilde{f}_t)\mid \mathcal{F}_t\right]\\
    ={}& \mathbb{E}\left[\norm{\left(f_t(x_t, \hat{a}_t) - f_t(x_t, a_t^*)\right) - v_t}^2 - \norm{v_t}^2 \mid \mathcal{F}_t\right]\\
    ={}& \norm{f_t(x_t, \hat{a}_t) - f_t(x_t, a_t^*)}^2 - 2 \left(f_t(x_t, \hat{a}_t) - f_t(x_t, a_t^*)\right)^\top \mathbb{E}\left[v_t \mid \mathcal{F}_t\right]\\
    ={}& \norm{f_t(x_t, \hat{a}_t) - f_t(x_t, a_t^*)}^2 = \varepsilon_t^2.
\end{align*}
Therefore, we obtain that
\begin{align*}
    \mathbb{E}\left[\sum_{t=0}^{T-1} \tilde{\ell}_t(x_t, \hat{a}_t, \tilde{f}_t) - \sum_{t=0}^{T-1} \tilde{\ell}_t(x_t, a_t^*, \tilde{f}_t)\right] ={}& \sum_{t=0}^{T-1}\mathbb{E}\left[\tilde{\ell}_t(x_t, \hat{a}_t, \tilde{f}_t) - \tilde{\ell}_t(x_t, a_t^*, \tilde{f}_t)\right]\\
    ={}& \sum_{t=0}^{T-1}\mathbb{E}\left[\mathbb{E}\left[\tilde{\ell}_t(x_t, \hat{a}_t, \tilde{f}_t) - \tilde{\ell}_t(x_t, a_t^*, \tilde{f}_t)\mid \mathcal{F}_t\right]\right]\\
    ={}& \sum_{t=0}^{T-1}\mathbb{E}\left[\varepsilon_t^2\right] = \mathbb{E}\left[\sum_{t=0}^{T-1} \varepsilon_t^2\right].
\end{align*}
Combining this with \eqref{lemma:gradient-estimator-regret:e1} gives that
\[\mathbb{E}\left[\sum_{t=0}^{T-1} \varepsilon_t^2\right] \leq 2\sqrt{3} (C_f + e_f)^3 D_f' \sqrt{C_p T}.\]
Then, we can apply \Cref{thm:gradient-error-bound} to conclude that
\[\mathbb{E}\left[\sum_{t=0}^{T-1} (\varepsilon_t')^2\right] \leq \frac{2m}{c}(1 + \gamma + \beta \gamma)\bar{\epsilon} T + 2 m \gamma^2 \bar{\epsilon}^2 T.\]

\section{Proof of Theorem \ref{thm:application-matched-disturbance-main}}\label{appendix:application-matched-disturbance-main}
By Lemma \ref{lemma:gradient-estimator-regret} and \Cref{thm:gradient-error-bound}, we know that the expected total prediction errors achieved by the gradient estimator satisfy that
\begin{align}\label{thm:application-matched-disturbance-main:e1}
    \mathbb{E}\left[\sum_{t=0}^{T-1} \varepsilon_t^2\right] &\leq 2\sqrt{3} (C_f + e_f)^3 D_f' \sqrt{C_p T}, \text{ and } \nonumber\\\mathbb{E}\left[\sum_{t=0}^{T-1} (\varepsilon_t')^2\right] &\leq \frac{2m}{c}(1 + \gamma + \beta \gamma)\bar{\epsilon} T + 2 m \gamma^2 \bar{\epsilon}^2 T.
\end{align}
By H\"older's inequality, we see that
\begin{align}\label{thm:application-matched-disturbance-main:e2}
    \mathbb{E}\left[\sum_{t=0}^{T-1} \varepsilon_t\right] &\leq \sqrt[4]{12} (C_f + e_f)^{\frac{3}{2}} (D_f')^{\frac{1}{2}} C_p^{\frac{1}{4}} T^{\frac{3}{4}}, \text{ and }\nonumber\\
    \mathbb{E}\left[\sum_{t=0}^{T-1} \varepsilon_t'\right] &\leq \sqrt{\frac{2}{c}(1 + \gamma + \beta \gamma) + 2 \gamma^2 \bar{\epsilon}} \cdot \sqrt{m \bar{\epsilon}} \cdot T.
\end{align}
By Lemma \ref{lemma:fully-actuated-assumptions}, we know that trajectory $\tilde{\xi}$ achieves the local regret
\begin{align}\label{thm:application-matched-disturbance-main:e3}
    R_L(T, \{\norm{\zeta_t}\}_{0\leq t\leq T-1}) = O\left(\frac{1}{\eta} (1 + V_{\text{sys}} + V_w) + \eta T + \eta^3 T + \frac{1}{\eta}\sum_{t=1}^{T-1} \norm{\zeta_t}\right),
\end{align}
By \Cref{thm:meta-dynamics-regret-and-stability} and Lemma \ref{lemma:matched-disturbance-Lipschitz-a}, we know that
\begin{align}\label{thm:application-matched-disturbance-main:e4}
    \mathbb{E}\left[\sum_{t=0}^{T-1} \norm{\zeta_t}\right] \leq{}& \left(\alpha_\theta + \sqrt{2} C (L_{\theta, x} + L_{\theta, y})(\alpha_x + \alpha_y)\right) \mathbb{E}\left[\sum_{t=0}^{T-1}\varepsilon_t\right]\nonumber\\
    &+ \left(\beta_\theta + \sqrt{2} C (L_{\theta, x} + L_{\theta, y})(\beta_x + \beta_y)\right) \mathbb{E}\left[\sum_{t=0}^{T-1}\varepsilon_t'\right]\nonumber\\
    \leq{}& C_0 \eta \mathbb{E}\left[\sum_{t=0}^{T-1}\varepsilon_t\right] + R_y L_{h, u} \cdot \eta \mathbb{E}\left[\sum_{t=0}^{T-1}\varepsilon_t'\right],
\end{align}
where $C = \frac{\bar{C} + C_{\ell, g, (x, x)} R_y}{1 - \rho} + \frac{C_{\ell, g, (\theta, x)} \bar{C}}{(1 - \rho)^2}$ by Lemma \ref{lemma:matched-disturbance-contraction-stability} and
\begin{align*}
    C_0 ={}& R_y (\ell_{h, x} + \ell_{h, u} L_{f, x} + \ell_{h, u} L_{\psi, x}) + \ell_{h, u} L_{\psi, \theta}\\
    &+ \sqrt{2} C \ell_{h, u} L_{\psi, \theta} \Bigg( R_y \left((\ell_{h, x} + \ell_{h, u} (L_{f, x} + L_{\psi, x})) (1 + L_{f, x} + L_{\psi, x}) + L_{h, u} (\ell_{f, x} + \ell_{\psi, x})\right)\\
    &+ \ell_{h, x} L_{\psi, \theta} + L_{h, u} \ell_{\psi, x} + \ell_{h, u} L_{\psi, \theta} (L_{f, x} + L_{\psi, x}) + L_{h, x} + L_{h, u} (L_{f, x} + L_{\psi, x})\Bigg).
\end{align*}
Substituting \eqref{thm:application-matched-disturbance-main:e4} into \eqref{thm:application-matched-disturbance-main:e3} and applying \eqref{thm:application-matched-disturbance-main:e1} and \eqref{thm:application-matched-disturbance-main:e2} give that
\[R_L(T, \{\norm{\zeta_t}\}_{0\leq t\leq T-1}) = O\left(\frac{1}{\eta}(1 + V_{sys} + \bar{\epsilon} \cdot T) + \eta T + (\sqrt{m \bar{\epsilon}} + m \bar{\epsilon}) \cdot T\right).\]
Further, by the last statement of \Cref{thm:meta-dynamics-regret-and-stability}, we obtain that
\[\mathbb{E}\left[\sum_{t=0}^{T-1} \left(\norm{x_t - \tilde{x}_t} + \norm{y_t - \tilde{y}_t}\right)\right] = O\left(T^{3/4} + \sqrt{m\bar{\epsilon}} \cdot T\right).\]

\section{Local Regret of Online Gradient Descent}
\label{appendix:nonconvex-biased-OGD}
\begin{theorem}\label{thm:nonconvex-biased-OGD}
Consider the parameter sequence $\{\theta_t\}$ that satisfies
\[\norm{\theta_{t+1} - (\theta_t - \eta \nabla_{\Theta, \eta}F_t(\theta_t))} \leq \eta \beta_t, \text{ for all }t \geq 0.\]
Suppose at every time $t$, $F_t$ is $\ell_F$-smooth and $L_F$-Lipschitz in $\Theta$. If the learning rate $\eta \leq \frac{1}{\ell_F}$, then the local regret $\sum_{t = 0}^{T-1}\norm{\nabla_{\Theta, \eta} F_t(\theta_t)}^2$ is upper bounded by
\begin{align*}
    &\frac{1}{\eta (1 - \eta \ell_F)}\left(F_0(\theta_0) + \sum_{t=1}^{T-1}\dist_s(F_t, F_{t-1})\right) + \frac{L_F \sum_{t=0}^{T-1} \beta_t + \ell_F \eta \sum_{t=0}^{T-1} \beta_t^2}{1 - \eta \ell_F},
\end{align*}
where $\dist_s(F, F') \coloneqq \sup_{\theta \in \Theta} \abs{F(\theta) - F'(\theta)}$.
\end{theorem}

Next, we state a property of projection onto the compact convex set $\Theta \in \mathbb{R}^d$ in Lemma \ref{lemma:projection-decrease-distance}. This is a classic result in convex optimization (see, for example, Theorem 1.2.1 in \cite{schneider2014convex}).

\begin{lemma}\label{lemma:projection-decrease-distance}
Let $q$ and $q'$ be arbitrary points in $\mathbb{R}^d$. Let $p = \Pi_\Theta(q)$ and $p' = \Pi_\Theta(q')$. Then, the following inequality holds:
\[\norm{p - p'} \leq \norm{q - q'}.\]
\end{lemma}

Now we come back to the proof of \Cref{thm:nonconvex-biased-OGD}.

Define the quantity
\[\epsilon_t \coloneqq \frac{1}{\eta} \left(\theta_{t+1} - (\theta_t - \eta \nabla_{\Theta, \eta}F_t(\theta_t))\right).\]
We see that
\begin{align}\label{thm:nonconvex-biased-OGD:e2}
    \theta_{t+1} - \theta_t = - \eta \nabla_{\Theta, \eta} F_t(\theta_t) + \eta \epsilon_t.
\end{align}

By the smoothness of $F_t(\cdot)$, we see that
\begin{subequations}\label{thm:nonconvex-biased-OGD:e3}
\begin{align}
    F_t(\theta_{t+1}) \leq{}& F_t(\theta_t) + \langle \nabla F_t(\theta_t), \theta_{t+1} - \theta_t\rangle + \frac{\ell_F}{2}\norm{\theta_{t+1} - \theta_t}^2\nonumber\\
    ={}& F_t(\theta_t) - \eta \langle \nabla F_t(\theta_t), \nabla_{\Theta, \eta} F_t(\theta_t) - \epsilon_t\rangle + \frac{\ell_F \eta^2}{2}\norm{\nabla_{\Theta, \eta} F_t(\theta_t) - \epsilon_t}^2 \label{thm:nonconvex-biased-OGD:e3:s1}\\
    ={}& F_t(\theta_t) - \eta \langle \nabla F_t(\theta_t), \nabla_{\Theta, \eta} F_t(\theta_t)\rangle + \frac{\ell_F \eta^2}{2}\norm{\nabla_{\Theta, \eta} F_t(\theta_t)}^2\nonumber\\
    &+ \eta \langle \nabla F_t(\theta_t), \epsilon_t\rangle - \ell_F \eta^2 \langle \nabla_{\Theta, \eta} F_t(\theta_t), \epsilon_t\rangle + \frac{\ell_F \eta^2}{2}\norm{\epsilon_t}^2, \label{thm:nonconvex-biased-OGD:e3:s2}
\end{align}
\end{subequations}
where we use \eqref{thm:nonconvex-biased-OGD:e2} in \eqref{thm:nonconvex-biased-OGD:e3:s1}. Recall that $\Theta$ is a closed convex subset of $\mathbb{R}^d$. Since $\theta_t - \eta \nabla_{\Theta, \eta} F_t(\theta_t)$ is the projection of $\theta_t - \eta \nabla F_t(\theta_t)$ onto $\Theta$ and $\theta_t \in \Theta$, we have
\[\langle (\theta_t - \eta \nabla F_t(\theta_t)) - (\theta_t - \eta \nabla_{\Theta, \eta} F_t(\theta_t)), \theta_t - (\theta_t - \eta \nabla_{\Theta, \eta} F_t(\theta_t))\rangle \leq 0.\]
Rearranging terms gives that
\[\langle \nabla F_t(\theta_t), \nabla_{\Theta, \eta} F_t(\theta_t)\rangle \geq \norm{\nabla_{\Theta, \eta} F_t(\theta_t)}^2.\]
Substituting this inequality into \eqref{thm:nonconvex-biased-OGD:e3} gives that
\begin{subequations}\label{thm:nonconvex-biased-OGD:e4}
\begin{align}
    F_t(\theta_{t+1}) \leq{}& F_t(\theta_t) - \eta \norm{\nabla_{\Theta, \eta} F_t(\theta_t)}^2 + \frac{\ell_F \eta^2}{2}\norm{\nabla_{\Theta, \eta} F_t(\theta_t)}^2\nonumber\\
    &+ \eta \langle \nabla F_t(\theta_t), \epsilon_t\rangle - \ell_F \eta^2 \langle \nabla_{\Theta, \eta} F_t(\theta_t), \epsilon_t\rangle + \frac{\ell_F \eta^2}{2}\norm{\epsilon_t}^2\nonumber\\
    \leq{}& F_t(\theta_t) - \eta (1 - \ell_F \eta) \norm{\nabla_{\Theta, \eta} F_t(\theta_t)}^2 + \eta \norm{\nabla F_t(\theta_t)} \cdot \norm{\epsilon_t}\nonumber\\
    &- \frac{\ell_F \eta^2}{2} \norm{\nabla_{\Theta, \eta} F_t(\theta_t) + \epsilon_t}^2 + \ell_F \eta^2 \norm{\epsilon_t}^2\label{thm:nonconvex-biased-OGD:e4:s1}\\
    \leq{}& F_t(\theta_t) - \eta (1 - \ell_F \eta) \norm{\nabla_{\Theta, \eta} F_t(\theta_t)}^2 + \eta L_F \beta_t + \ell_F \eta^2 \beta_t^2,\label{thm:nonconvex-biased-OGD:e4:s2}
\end{align}
\end{subequations}
where we rearrange the terms and use the Cauchy-Schwarz inequality in \eqref{thm:nonconvex-biased-OGD:e4:s1}; In \eqref{thm:nonconvex-biased-OGD:e4:s2}, we use the assumption $\norm{\epsilon_t} \leq \beta_t$. Summing \eqref{thm:nonconvex-biased-OGD:e4} over $t = 0, 1, \ldots, T-1$ gives that
\begin{subequations}\label{thm:nonconvex-biased-OGD:e5}
\begin{align}
    &\eta(1 - \ell_F \eta) \sum_{t=0}^{T-1} \norm{\nabla_{\Theta, \eta}F_t(\theta_t)}^2\nonumber\\
    \leq{}& \sum_{t=0}^{T-1} \left(F_t(\theta_t) - F_t(\theta_{t+1})\right) + \eta L_F \sum_{t=0}^{T-1} \beta_t + \ell_F \eta^2 \sum_{t=0}^{T-1} \beta_t^2\nonumber\\
    \leq{}& F_0(\theta_0) + \sum_{t=1}^{T-1} \left(F_t(\theta_t) - F_{t-1}(\theta_t)\right) + \sum_{t=1}^{T-1}\dist_s(F_t, F_{t-1}) + \eta L_F \sum_{t=0}^{T-1} \beta_t + \ell_F \eta^2 \sum_{t=0}^{T-1} \beta_t^2\label{thm:nonconvex-biased-OGD:e5:s1}\\
    \leq{}& F_0(\theta_0) + \sum_{t=1}^{T-1}\dist_s(F_t, F_{t-1}) + \eta L_F \sum_{t=0}^{T-1} \beta_t + \ell_F \eta^2 \sum_{t=0}^{T-1} \beta_t^2,\label{thm:nonconvex-biased-OGD:e5:s2}
\end{align}
\end{subequations}
where we rearrange the terms and use $F_{T-1}(\theta_T) \geq 0$ in \eqref{thm:nonconvex-biased-OGD:e5:s1}; we use the definition of $\dist_s(\cdot, \cdot)$ in \eqref{thm:nonconvex-biased-OGD:e5:s2}.

\section{Useful Lemmas}
\label{appendix:useful-lemmas}
In this section, we summarize some useful existing results in \cite{lin2023online} that can help us in the proof of M-GAPS (Algorithm \ref{alg:M-GAPS}). We can build our proof upon some results shown in \cite{lin2023online} because of the similarity between our M-GAPS algorithm and the GAPS algorithm proposed by \cite{lin2023online} when applied to known dynamical systems: Both algorithms are designed to efficiently approximate the gradient $\nabla F_t(\theta_t)$ of the surrogate cost. Note that $\nabla F_t(\theta_t)$ can be expressed as
\[\nabla F_t(\theta_t) = \sum_{\tau=0}^t \left.\frac{\partial h_{t\mid 0}^*}{\partial \theta_{t - \tau}}\right|_{x_0, (\theta_t)_{\times (t+1)}}.\]
M-GAPS adopts the approximation $G_t$ that is given by
\begin{align}\label{equ:M-GAPS-gradient-approx}
    G_t = \sum_{\tau=0}^t \left.\frac{\partial h_{t\mid 0}^*}{\partial \theta_{t - \tau}}\right|_{x_0, \theta_{0:t}},
\end{align}
which simplifies $\nabla F_t(\theta_t)$ by replacing the imaginary trajectory achieved by using policy parameter $\theta_t$ repeatedly with the actual trajectory. The approximator of GAPS, which we denote as $G_t'$, takes an additional step to approximate $G_t$ by truncating the summation from time $0$ to $t$ to at most $B$ time steps, i.e.,
\[G_t' = \sum_{\tau=\max\{0, t - B\}}^t \left.\frac{\partial h_{t\mid 0}^*}{\partial \theta_{t - \tau}}\right|_{x_0, \theta_{0:t}},\]
where $B$ is the buffer length parameter decided by the algorithm. Intuitively, the approximation $G_t$ adopted by M-GAPS is closer to $\nabla F_t(\theta_t)$, which allows us to show the same guarantees for M-GAPS as GAPS when the true dynamics are known.

In this section, we translate some results from \cite{lin2023online} into the settings of matched-disturbance dynamics application discussed in \Cref{sec:application}. Lemma \ref{lemma:smooth-multi-step-dynamics} is Lemma D.3 in \cite{lin2023online}. We changed the condition $x_\tau, x_\tau' \in B_n(0, R_S + C \norm{x_0})$ to $x_\tau, x_\tau' \in B_n(0, \bar{R}_C)$, where $\bar{R}_C$ can be any positive number that satisfies $\bar{R}_C < R_C$ and $C \bar{R}_C + R_S \leq R_x$. This minor change will not affect the proof provided in \cite{lin2023online}.

\begin{lemma}[Lipschitzness/Smoothness of the Multi-Step Dynamics]\label{lemma:smooth-multi-step-dynamics}
Suppose Assumptions \ref{assump:Lipschitz-and-smoothness} and \ref{assump:contractive-and-stability} hold. Given two time steps $t > \tau$, for any $x_\tau, x_\tau' \in B_n(0, \bar{R}_C)$ and $\theta_\tau, \theta_\tau' \in \Theta,$ $\theta_{\tau+1:t-1} \in S_{\varepsilon}(\tau+1:t-1)$, if $x_{\tau+1}' \coloneqq g_{\tau + 1\mid \tau}(x_\tau', \theta_\tau')$ is also in $B_n(0, \bar{R}_C)$, the multi-step dynamical function $g_{t\mid \tau}$ satisfies that
\begin{align*}
    &\norm{\left.\frac{\partial g_{t\mid \tau}^*}{\partial x_\tau} \right|_{x_\tau, \theta_{\tau:t-1}}} \leq C_{L, g, x} \decayfactor^{t - \tau}, \quad \norm{\left.\frac{\partial g_{t\mid \tau}^*}{\partial \theta_\tau} \right|_{x_\tau, \theta_{\tau:t-1}}} \leq C_{L, g, \theta} \decayfactor^{t-\tau}, \forall \theta_{\tau:t-1} \in S_{\varepsilon}(\tau:t-1),\\
    &\norm{\left.\frac{\partial g_{t\mid \tau}^*}{\partial x_\tau} \right|_{x_\tau, \theta_{\tau:t-1}} - \left.\frac{\partial g_{t\mid \tau}^*}{\partial x_\tau} \right|_{x_\tau', \theta_\tau', \theta_{\tau+1:t-1}}} \leq C_{\ell, g, (x, x)} \decayfactor^{t - \tau} \norm{x_\tau - x_\tau'} + C_{\ell, g, (x, \theta)} \decayfactor^{t - \tau} \norm{\theta_\tau - \theta_\tau'},\\
    &\norm{\left.\frac{\partial g_{t\mid \tau}^*}{\partial \theta_\tau} \right|_{x_\tau, \theta_{\tau:t-1}} - \left.\frac{\partial g_{t\mid \tau}^*}{\partial \theta_\tau} \right|_{x_\tau', \theta_\tau', \theta_{\tau+1:t-1}}} \leq C_{\ell, g, (\theta, x)} \decayfactor^{t - \tau} \norm{x_\tau - x_\tau'} + C_{\ell, g, (\theta, \theta)} \decayfactor^{t - \tau} \norm{\theta_\tau - \theta_\tau'},
\end{align*}
where $C_{L, g, x} = \bar{C}, C_{L, g, \theta} = \frac{\bar{C} L_{\phi, u} L_{\psi, \theta}}{\decayfactor}$, and
\begin{align*}
    C_{\ell, g, (x, x)} ={}& {\left((1 + L_{\psi, x}) \left(\ell_{\phi, x} + \ell_{\phi, u} L_{\psi, x}\right) + L_{\phi, x} \ell_{\psi, x}\right) C^3}{\decayfactor^{-1}(1 - \decayfactor)^{-1}},\\
    C_{\ell, g, (x, \theta)} ={}& {\left((1 + L_{\psi, x}) \left(\ell_{\phi, x} + \ell_{\phi, u} L_{\psi, x}\right) + L_{\phi, x} \ell_{\psi, x}\right) C^3 L_{\phi, u} L_{\psi, \theta}}{\decayfactor^{-1}(1 - \decayfactor)^{-1}}\\
    &+ {\left((1 + L_{\psi, x}) \ell_{\phi, u} L_{\psi, \theta} + L_{\phi, u} \ell_{\psi, \theta}\right) C}{\decayfactor^{-1}(1 - \decayfactor)^{-1}},\\
    C_{\ell, g, (\theta, x)} ={}& \left((1 + L_{\psi, x}) \left(\ell_{\phi, x} + \ell_{\phi, u} L_{\psi, x}\right) + L_{\phi, x} \ell_{\psi, x}\right) (L_{\phi, x} + L_{\phi, u} L_{\psi, x})\cdot\\
    &C^3 L_{\phi, u} L_{\psi, \theta} \decayfactor^{-2}(1 - \decayfactor)^{-1}
    + {C \left(L_{\psi, \theta}(\ell_{\phi, x} + \ell_{\phi, u} L_{\psi, x}) + L_{\phi, u} \ell_{\psi, x}\right)}{\decayfactor^{-1}},\\
    C_{\ell, g, (\theta, \theta)} ={}& {\left((1 + L_{\psi, x}) \left(\ell_{\phi, x} + \ell_{\phi, u} \cdot L_{\psi, x}\right) + L_{\phi, x} \cdot \ell_{\psi, x}\right) L_{\phi, u}^2 L_{\psi, \theta}^2 C^3}{\decayfactor^{-2}(1 - \decayfactor)^{-1}}\\
    &+ {\left(L_{\phi, u} \ell_{\psi, \theta} + \ell_{\phi, u} L_{\psi, \theta}^2\right)C}{\decayfactor^{-1}}.
\end{align*}
\end{lemma}

Corollary \ref{coro:smooth-multi-step-costs} is implied by Lemma \ref{lemma:smooth-multi-step-dynamics} and corresponds to Corollary D.4 in \cite{lin2023online}.

\begin{corollary}[Lipschitzness/Smoothness of the Multi-Step Costs]\label{coro:smooth-multi-step-costs}
Under the same assumptions as Lemma \ref{lemma:smooth-multi-step-dynamics}, the multi-step cost function $h_{t\mid \tau}$ satisfies that
\begin{align*}
    &\norm{\left.\frac{\partial h_{t\mid \tau}}{\partial x_\tau} \right|_{x_\tau, \theta_{\tau:t}}} \leq C_{L, h, x} \decayfactor^{t - \tau}, \norm{\left.\frac{\partial h_{t\mid \tau}}{\partial \theta_\tau} \right|_{x_\tau, \theta_{\tau:t}}} \leq C_{L, h, \theta} \decayfactor^{t-\tau},\\
    &\norm{\left.\frac{\partial h_{t\mid \tau}}{\partial x_\tau} \right|_{x_\tau, \theta_\tau, \theta_{\tau+1:t}} - \left.\frac{\partial h_{t\mid \tau}}{\partial x_\tau} \right|_{x_\tau', \theta_\tau', \theta_{\tau+1:t}}} \leq C_{\ell, h, (x, x)} \decayfactor^{t - \tau} \norm{x_\tau - x_\tau'} + C_{\ell, h, (x, \theta)} \decayfactor^{t - \tau} \norm{\theta_\tau - \theta_\tau'},\\
    &\norm{\left.\frac{\partial h_{t\mid \tau}}{\partial \theta_\tau} \right|_{x_\tau, \theta_\tau, \theta_{\tau+1:t}} - \left.\frac{\partial h_{t\mid \tau}}{\partial \theta_\tau} \right|_{x_\tau', \theta_\tau', \theta_{\tau+1:t}}} \leq C_{\ell, h, (\theta, x)} \decayfactor^{t - \tau} \norm{x_\tau - x_\tau'} + C_{\ell, h, (\theta, \theta)} \decayfactor^{t - \tau} \norm{\theta_\tau - \theta_\tau'},
\end{align*}
where $C_{L, h, x} = L_h C (1 + L_{\psi, x}), C_{L, h, \theta} = L_h \max\{C_{L, \phi, \theta} (1 + L_{\psi, x}), L_{\psi, \theta}\}$, and
\begin{align*}
    C_{\ell, h, (x, x)} &= L_h (1 + L_{\psi, x}) C_{\ell, \phi, (x, x)} + ((\ell_{h, x} + \ell_{h, u}L_{\psi, x})(1 + L_{\psi, x}) + L_h \ell_{\psi, x}) C_{L, \phi, x}^2,\\
    C_{\ell, h, (x, \theta)} &= L_h (1 + L_{\psi, x}) C_{\ell, \phi, (x, \theta)} + ((\ell_{h, x} + \ell_{h, u}L_{\psi, x})(1 + L_{\psi, x}) + L_h \ell_{\psi, x}) C_{L, \phi, x} C_{L, \phi, \theta},\\
    C_{\ell, h, (\theta, x)} &= L_h (1 + L_{\psi, x}) C_{\ell, \phi, (\theta, x)} + ((\ell_{h, x} + \ell_{h, u}L_{\psi, x})(1 + L_{\psi, x}) + L_h \ell_{\psi, x}) C_{L, \phi, x} C_{L, \phi, \theta},\\
    C_{\ell, h, (\theta, \theta)} &= L_h (1 + L_{\psi, x}) C_{\ell, \phi, (\theta, \theta)} + ((\ell_{h, x} + \ell_{h, u}L_{\psi, x})(1 + L_{\psi, x}) + L_h \ell_{\psi, x}) C_{L, \phi, \theta}^2.
\end{align*}
\end{corollary}

\Cref{thm:OCO-with-parameter-update-trajectory-distance} bounds the distances between the trajectory of M-GAPS with the imaginary trajectory achieved by using $\theta_t$ repeatedly from time step $0$. It can be shown using a similar approach as Theorem D.5 in \cite{lin2023online}, while a difference is that we consider an additional disturbance $\zeta_t$ in the update rule of policy parameters. We include the proof of \Cref{thm:OCO-with-parameter-update-trajectory-distance} in Appendix \ref{appendix:proof:thm:OCO-with-parameter-update-trajectory-distance} for completeness.

\begin{theorem}\label{thm:OCO-with-parameter-update-trajectory-distance}
Suppose Assumptions \ref{assump:Lipschitz-and-smoothness} and \ref{assump:contractive-and-stability} hold. Let $\{x_t, u_t, \theta_t\}_{t \in \mathcal{T}}$ denote the trajectory of 
\begin{subequations}\label{equ:M-GAPS-dynamics-exact-add-bias}
\begin{align}
    x_{t+1} ={}& q_t^x(x_t, y_t, \theta_t, a_t^*) = \phi_t(x_t, \psi_t(x_t, \theta_t)) + w_t,\\
    y_{t+1} ={}& q_t^y(x_t, y_t, \theta_t, a_t^*) = \left.\frac{\partial g^*_{t+1\mid t}}{\partial x_t}\right|_{x_t, \theta_t} \cdot y_t + \left.\frac{\partial g^*_{t+1\mid t}}{\partial \theta_t}\right|_{x_t, \theta_t},\\
    \theta_{t+1} ={}& q_t^\theta(x_t, y_t, \theta_t, a_t^*) = \Pi_\Theta\left(\theta_{t+1} - \eta \left(\left.\frac{\partial h_{t\mid t}^*}{\partial x_t}\right|_{x_t, \theta_t} \cdot y_t + \left.\frac{\partial h_{t\mid t}^*}{\partial \theta_t}\right|_{x_t, \theta_t}\right)\right) + \zeta_t.
\end{align}
\end{subequations}
Suppose $\eta$ and $\bar{\zeta}$ satisfy the constraint that $\bar{\varepsilon} \coloneqq \frac{C_{L, h, \theta}\eta}{1 - \rho} + \bar{\zeta} \leq \varepsilon$.
Then, both $\norm{G_t}$ and $\norm{\nabla F_t(\theta_t)}$ are upper bounded by $\frac{C_{L, h, \theta}}{1 - \decayfactor}$, and the following inequalities holds for any two time steps $\tau, t$ ($\tau \leq t$):
\begin{align*}
    &\norm{\theta_t - \theta_{\tau}} \leq \frac{C_{L, h, \theta}}{1 - \decayfactor}\cdot (t - \tau)\eta + \sum_{\tau' = \tau}^{t-1} \norm{\zeta_{\tau'}}, \text{ and } \norm{x_\tau - \hat{x}_\tau(\theta_t)} \leq \\
    &\frac{C_{L, h, \theta} C_{L, \phi, \theta} \decayfactor}{(1 - \decayfactor)^2} \left((t - \tau) + \frac{1}{1 - \decayfactor}\right) \cdot \eta + \frac{C_{L, \phi, \theta} \decayfactor}{1 - \rho}\cdot \left(\sum_{\tau' = \tau}^{t-1} \norm{\zeta_{\tau'}} + \sum_{\tau' = 0}^{\tau-1} \rho^{\tau - \tau'}\norm{\zeta_{\tau'}}\right),
\end{align*}
where we use the notation $\hat{x}_\tau(\theta) \coloneqq g_{\tau\mid 0}^*(x_0, \theta_{\times (\tau + 1)}), \forall \theta \in \Theta$. Further, we have that
\begin{align*}
    \abs{h_t(x_t, u_t, \theta_t) - F_t(\theta_t)} \leq{}& \frac{C_{L, h, \theta} C_{L, \phi, \theta}L_h (1 + L_{\psi, x} + L_{f, x}) \decayfactor}{(1 - \decayfactor)^3} \cdot \eta\\
    &+ \frac{C_{L, \phi, \theta} L_h (1 + L_{\psi, x} + L_{f, x}) \decayfactor}{1 - \rho} \cdot \sum_{\tau = 0}^{t-1} \rho^{t-\tau} \norm{\zeta_\tau}.
\end{align*}
\end{theorem}

Recall that we define the gradient approximation $G_t$ for M-GAPS in \eqref{equ:M-GAPS-gradient-approx}. Using this notation, the update rule of $\theta_{0:T-1}$ in joint dynamics \eqref{equ:M-GAPS-dynamics-exact-add-bias} can be simplified as
\[\theta_{t+1} = \Pi_\Theta\left(\theta_{t+1} - \eta G_t\right) + \zeta_t.\]
To compare the trajectory of M-GAPS with the trajectory achieved by the online gradient descent trajectory $\theta_{t+1} = \Pi_\Theta(\theta_t - \eta \nabla F_t(\theta_t))$, we bound the difference between $G_t$ and $\nabla F_t(\theta_t)$ in \Cref{thm:OCO-with-parameter-update-gradient-bias:full}. We provide its proof in Appendix \ref{appendix:proof:thm:OCO-with-parameter-update-gradient-bias:full} for completeness.

\begin{theorem}[Gradient Bias]\label{thm:OCO-with-parameter-update-gradient-bias:full}
Suppose Assumptions \ref{assump:Lipschitz-and-smoothness} and \ref{assump:contractive-and-stability} hold. Let $\{x_t, u_t, \theta_t\}_{t \in \mathcal{T}}$ denote the trajectory of \eqref{equ:M-GAPS-dynamics-exact-add-bias}. Suppose $\eta$ and $\bar{\zeta}$ satisfy the constraint that $\bar{\varepsilon} \coloneqq \frac{C_{L, h, \theta}\eta}{1 - \rho} + \bar{\zeta} \leq \varepsilon$. Then, the following holds for all $\tau \leq t$:
\begin{align*}
    &\norm{\left.\frac{\partial h_{t\mid 0}^*}{\partial \theta_\tau} \right|_{x_0, \theta_{0:t}} - \left.\frac{\partial h_{t\mid 0}^*}{\partial \theta_\tau} \right|_{x_0, (\theta_t)_{\times (t+1)}}}\\ 
    \leq{}& \left(\hat{C}_0 + \hat{C}_1 (t - \tau) + \hat{C}_2 (t - \tau)^2\right) \decayfactor^{t - \tau} \cdot \eta + \left(\hat{C}_3 \sum_{\tau' = \tau}^{t-1} \norm{\zeta_{\tau'}} + \hat{C}_4 \sum_{\tau' = \tau}^{t-1} (\tau'-\tau) \norm{\zeta_{\tau'}}\right) \cdot \rho^{t-\tau}\\
    &+ \hat{C}_5 \sum_{\tau' = 0}^{\tau-1} \rho^{t-\tau'} \norm{\zeta_{\tau'}}. 
\end{align*}
for
\begin{align*}
    &\hat{C}_0 = \frac{\decayfactor C_{L, h, \theta} C_{L, \phi, \theta} C_{\ell, h, (\theta, x)}}{(1 - \decayfactor)^3},\  \hat{C}_1 = \frac{(1 - \decayfactor) C_{L, h, \theta} C_{\ell, h, (\theta, x)} + \decayfactor C_{L, h, \theta} C_{L, \phi, \theta} C_{\ell, h, (\theta, \theta)}}{(1 - \decayfactor)^2},\\
    &\hat{C}_2 = \frac{C_{L, h, \theta} C_{\ell, h, (x, \theta) C_{L, \phi, \theta}}}{1 - \rho},\ \hat{C}_3 = \frac{C_{L, \phi, \theta} C_{\ell, h, (\theta, x)} \decayfactor}{1 - \rho},\\
    &\hat{C}_4 = {C_{\ell, h, (x, \theta)}C_{L, \phi, \theta}}, \ \hat{C}_5 = \frac{C_{L, \phi, \theta} C_{\ell, h, (\theta, x)} \decayfactor}{1 - \rho}.
\end{align*}
Next,
\begin{align*}
    \norm{G_t - \nabla F_t(\theta_t)} \leq{}& \left(\frac{\hat{C}_0}{1 - \decayfactor} + \frac{\hat{C}_1 + \hat{C}_2}{(1 - \decayfactor)^{2}} + \frac{\hat{C}_2}{(1 - \decayfactor)^{3}}\right) \eta\\
    &+ \sum_{\tau = 0}^{t-1}\left(\frac{\hat{C}_3}{1 - \rho}
    + \frac{\hat{C}_4}{(1 - \rho)^2} + \hat{C}_5 (t - \tau)\right) \rho^{t-\tau} \norm{\zeta_\tau}.
\end{align*}
\end{theorem}

\subsection{Proof of Theorem \ref{thm:OCO-with-parameter-update-trajectory-distance}}\label{appendix:proof:thm:OCO-with-parameter-update-trajectory-distance}
We first use induction to show that for all time step $t \in \mathcal{T}$,
\begin{equation}\label{thm:OCO-with-parameter-update-trajectory-distance:e0}
    \norm{G_t} \leq \frac{C_{L, h, \theta}}{1 - \decayfactor}, x_t \in B_n(0, R_S + \bar{C} \norm{x_0}), u_t \in \mathcal{U}, \text{ and }\norm{\theta_{t+1} - \theta_t} \leq \epsilon_\theta,
\end{equation}
where $\mathcal{U} = \{\psi(x, \theta) - f(x, a) \mid x \in B_n(0, R_x), \theta \in \Theta, a\in \mathcal{A}, (\psi, f) \in \mathcal{G}\}$.

Note that $\norm{G_0} \leq C_{L, h, \theta} \leq \frac{C_{L, h, \theta}}{1 - \decayfactor}$ by \Cref{coro:smooth-multi-step-costs}. We also have $x_0 \in B_n(0, R_S + C\norm{x_0})$ and $u_0 \in \mathcal{U}$.

Suppose $\norm{G_{t-1}} \leq \frac{C_{L, h, \theta}}{1 - \decayfactor}$ for some $t \geq 1$. Then, since $\eta \leq \frac{(1 - \rho)\epsilon_\theta}{C_{L, h, \theta}}$ and the projection onto $\Theta$ is a contraction, we see that
\[\norm{\theta_t - \theta_{t-1}} \leq \norm{\eta G_{t-1}} + \norm{\zeta_t} \leq \epsilon_\theta.\]
Suppose $\norm{\theta_\tau - \theta_{\tau-1}} \leq \epsilon_\theta$ holds for all $\tau \leq t$, i.e., $\theta_{0:t} \in S_{\epsilon_\theta}(0:t)$. By Lemma D.2 in \cite{lin2023online}, we see that
\[x_t \in B_n(0, R_S + C \norm{x_0}), \text{ and } u_t \in \mathcal{U}.\]
Therefore, by taking norm on both sides of \eqref{equ:M-GAPS-gradient-approx}, we see that
\begin{subequations}\label{thm:OCO-with-parameter-update-trajectory-distance:e1}
\begin{align}
    \norm{G_t} ={}& \norm{\sum_{\tau = 0}^{t} \left.\frac{\partial h_{t\mid t-\tau}}{\partial \theta_{t - \tau}}\right|_{x_{t - \tau}, \theta_{t-\tau:t}}}\nonumber\\
    \leq{}& \sum_{\tau = 0}^{t} \norm{\left.\frac{\partial h_{t\mid t-\tau}}{\partial \theta_{t - \tau}}\right|_{x_{t - \tau}, \theta_{t-\tau:t}}} \label{thm:OCO-with-parameter-update-trajectory-distance:e1:s1}\\
    \leq{}& \sum_{\tau = 0}^{t} C_{L, h, \theta} \decayfactor^\tau \label{thm:OCO-with-parameter-update-trajectory-distance:e1:s2}\\
    \leq{}& \frac{C_{L, h, \theta}}{1 - \decayfactor},\nonumber
\end{align}
\end{subequations}
where we use the triangle inequality in \eqref{thm:OCO-with-parameter-update-trajectory-distance:e1:s1} and \Cref{coro:smooth-multi-step-costs} in \eqref{thm:OCO-with-parameter-update-trajectory-distance:e1:s2}. Note that we can apply \Cref{coro:smooth-multi-step-costs} because $x_t \in B_n(0, R_S + C \norm{x_0})$. Therefore, we have shown \eqref{thm:OCO-with-parameter-update-trajectory-distance:e0} by induction. One can use the same technique as \eqref{thm:OCO-with-parameter-update-trajectory-distance:e1} to show $\norm{\nabla F_t(\theta_t)} \leq \frac{C_{L, f, \theta}}{1 - \decayfactor}$.

Since the projection onto the set $\Theta$ is a contraction, we obtain that for any $t > \tau$,
\begin{align}\label{thm:OCO-with-parameter-update-trajectory-distance:e2}
    \norm{\theta_t - \theta_\tau} \leq \frac{C_{L, h, \theta}}{1 - \decayfactor}\cdot (t - \tau)\eta + \sum_{\tau' = \tau}^{t-1} \norm{\zeta_{\tau'}}.
\end{align}
Now we bound the distance between $x_\tau$ and $\hat{x}_\tau(\theta_t)$ for $\tau \leq t$. We see that
\begin{subequations}\label{thm:OCO-with-parameter-update-trajectory-distance:e3}
\begin{align}
    \norm{x_\tau - \hat{x}_\tau(\theta_t)} ={}& \norm{g_{\tau\mid 0}^*(x_0, \theta_{0:\tau-1}) - g_{\tau\mid 0}^*(x_0, (\theta_t)_{\times \tau})}\nonumber\\
    \leq{}& \sum_{\tau' = 0}^{\tau-1} \norm{g_{\tau\mid 0}^*(x_0, \theta_{0:\tau'}, (\theta_t)_{\times (\tau - \tau' -1)}) - g_{\tau\mid 0}^*(x_0, \theta_{0:\tau'-1}, (\theta_t)_{\times (\tau - \tau')})} \label{thm:OCO-with-parameter-update-trajectory-distance:e3:s1}\\
    \leq{}& \sum_{\tau' = 0}^{\tau-1}\norm{g_{\tau\mid \tau'}^*(x_{\tau'}, \theta_{\tau'}, (\theta_t)_{\times (\tau - \tau' -1)}) - g_{\tau\mid \tau'}^*(x_{\tau'}, (\theta_t)_{\times (\tau - \tau')})} \label{thm:OCO-with-parameter-update-trajectory-distance:e3:s2}\\
    \leq{}& \sum_{\tau' = 0}^{\tau-1} C_{L, g, \theta} \decayfactor^{\tau - \tau'} \norm{\theta_t - \theta_{\tau'}} \label{thm:OCO-with-parameter-update-trajectory-distance:e3:s3}\\
    \leq{}& \frac{C_{L, h, \theta} C_{L, g, \theta} \eta}{1 - \decayfactor} \sum_{\tau' = 0}^{\tau-1} \left(\frac{C_{L, h, \theta}}{1 - \decayfactor}\cdot (t - \tau')\eta + \sum_{\tau'' = \tau'}^{t-1} \norm{\zeta_{\tau''}}\right) \label{thm:OCO-with-parameter-update-trajectory-distance:e3:s4}\\
    \leq{}& \frac{C_{L, h, \theta} C_{L, \phi, \theta} \decayfactor}{(1 - \decayfactor)^2} \left((t - \tau) + \frac{1}{1 - \decayfactor}\right) \cdot \eta\nonumber\\ 
    &+ \frac{C_{L, \phi, \theta} \decayfactor}{1 - \rho}\cdot \left(\sum_{\tau' = \tau}^{t-1} \norm{\zeta_{\tau'}} + \sum_{\tau' = 0}^{\tau-1} \rho^{\tau - \tau'}\norm{\zeta_{\tau'}}\right), \nonumber
\end{align}
\end{subequations}
where we use the triangle inequality in \eqref{thm:OCO-with-parameter-update-trajectory-distance:e3:s1};
we use the definition of multi-step dynamics in \eqref{thm:OCO-with-parameter-update-trajectory-distance:e3:s2};
we use \Cref{lemma:smooth-multi-step-dynamics} in \eqref{thm:OCO-with-parameter-update-trajectory-distance:e3:s3};
we use \eqref{thm:OCO-with-parameter-update-trajectory-distance:e2} in \eqref{thm:OCO-with-parameter-update-trajectory-distance:e3:s4}.

Similarly, since $x_t \in B_n(0, R_S + C \norm{x_0})$ and we also see that $\hat{x}_t(\theta_t) \in B_n(0, R_S + C \norm{x_0})$, we obtain that
\begin{subequations}\label{thm:OCO-with-parameter-update-trajectory-distance:e4}
\begin{align}
    \abs{h_t(x_t, u_t, \theta_t) - F_t(\theta_t)} ={}& \abs{h_t(x_t, u_t, \theta_t) - h_t(\hat{x}_t(\theta_t), \hat{u}_t(\theta_t), \theta_t)}\nonumber\\
    \leq{}& L_h\left(\norm{x_t - \hat{x}_t(\theta_t)} + \norm{u_t - \hat{u}_t(\theta_t)}\right) \label{thm:OCO-with-parameter-update-trajectory-distance:e4:s1}\\
    \leq{}& L_h (1 + L_{\psi, x} + L_{f, x}) \norm{x_t - \hat{x}_t(\theta_t)} \label{thm:OCO-with-parameter-update-trajectory-distance:e4:s3}\\
    \leq{}& \frac{C_{L, h, \theta} C_{L, \phi, \theta}L_h (1 + L_{\psi, x} + L_{f, x}) \decayfactor}{(1 - \decayfactor)^3} \cdot \eta\nonumber\\
    &+ \frac{C_{L, \phi, \theta} L_h (1 + L_{\psi, x} + L_{f, x}) \decayfactor}{1 - \rho} \cdot \sum_{\tau = 0}^{t-1} \rho^{t-\tau} \norm{\zeta_\tau}, \label{thm:OCO-with-parameter-update-trajectory-distance:e4:s4}
\end{align}
\end{subequations}
where we use \Cref{assump:Lipschitz-and-smoothness} in \eqref{thm:OCO-with-parameter-update-trajectory-distance:e4:s1} and \eqref{thm:OCO-with-parameter-update-trajectory-distance:e4:s3};
we use \eqref{thm:OCO-with-parameter-update-trajectory-distance:e3} in \eqref{thm:OCO-with-parameter-update-trajectory-distance:e4:s4}.

\subsection{Proof of Theorem \ref{thm:OCO-with-parameter-update-gradient-bias:full}}\label{appendix:proof:thm:OCO-with-parameter-update-gradient-bias:full}
To simplify the notation, we adopt the shorthand notations $\hat{x}_\tau(\theta) \coloneqq g_{\tau\mid 0}^*(x_0, \theta_{\times \tau})$ and $\hat{u}_\tau(\theta) \coloneqq \pi_\tau(\hat{x}_\tau(\theta), \theta)$ throughout the proof.

We use the triangle inequality to do the decomposition
\begin{align}\label{thm:OCO-with-parameter-update-gradient-bias:e1}
    &\norm{\left.\frac{\partial h_{t\mid 0}^*}{\partial \theta_\tau} \right|_{x_0, \theta_{0:t}} - \left.\frac{\partial h_{t\mid 0}^*}{\partial \theta_\tau} \right|_{x_0, (\theta_t)_{\times (t + 1)}}}\nonumber\\
    ={}&\norm{\left.\frac{\partial h_{t\mid \tau}^*}{\partial \theta_\tau} \right|_{x_\tau, \theta_{\tau:t}} - \left.\frac{\partial h_{t\mid \tau}^*}{\partial \theta_\tau} \right|_{\hat{x}_\tau(\theta_t), (\theta_t)_{\times (t - \tau + 1)}}} \nonumber\\
    \leq{}& \norm{\left.\frac{\partial h_{t\mid \tau}^*}{\partial \theta_\tau} \right|_{x_\tau, \theta_\tau, (\theta_t)_{\times (t - \tau)}} - \left.\frac{\partial h_{t\mid \tau}^*}{\partial \theta_\tau} \right|_{\hat{x}_\tau(\theta_t), (\theta_t)_{\times (t - \tau + 1)}}}\nonumber\\
    &+ \sum_{\tau' = \tau + 1}^{t-1} \norm{\left.\frac{\partial h_{t\mid \tau}^*}{\partial \theta_\tau} \right|_{x_\tau, \theta_{\tau:\tau'}, (\theta_t)_{\times (t - \tau')}} - \left.\frac{\partial h_{t\mid \tau}^*}{\partial \theta_\tau} \right|_{x_\tau, \theta_{\tau:\tau'-1}, (\theta_t)_{\times (t - \tau' + 1)}}}.
\end{align}

Note that we can apply \Cref{coro:smooth-multi-step-costs} to bound each term in \eqref{thm:OCO-with-parameter-update-gradient-bias:e1}. For the first term in \eqref{thm:OCO-with-parameter-update-gradient-bias:e1}, since $x_\tau, \hat{x}_\tau(\theta_t), x_{\tau+1} \in B_n(0, \bar{R}_C)$, we see that
\begin{subequations}\label{thm:OCO-with-parameter-update-gradient-bias:e2}
\begin{align}
    &\norm{\left.\frac{\partial h_{t\mid \tau}^*}{\partial \theta_\tau} \right|_{x_\tau, \theta_\tau, (\theta_t)_{\times (t - \tau)}} - \left.\frac{\partial h_{t\mid \tau}^*}{\partial \theta_\tau} \right|_{\hat{x}_\tau(\theta_t), (\theta_t)_{\times (t - \tau + 1)}}}\nonumber\\
    \leq{}& \decayfactor^{t - \tau} \left(C_{\ell, h, (\theta, x)} \norm{x_\tau - \hat{x}_\tau(\theta_t)} + C_{\ell, h, (\theta, \theta)} \norm{\theta_t - \theta_\tau}\right) \label{thm:OCO-with-parameter-update-gradient-bias:e2:s1}\\
    \leq{}& \frac{(1 - \decayfactor)C_{L, h, \theta} C_{\ell, h, (\theta, x)} + \decayfactor C_{L, h, \theta} C_{L, \phi, \theta} C_{\ell, h, (\theta, \theta)}}{(1 - \decayfactor)^2} \cdot (t - \tau) \decayfactor^{t - \tau} \cdot \eta\nonumber\\
    &+ \frac{\decayfactor C_{L, h, \theta} C_{L, \phi, \theta} C_{\ell, h, (\theta, x)}}{(1 - \decayfactor)^3}\cdot \decayfactor^{t - \tau} \cdot \eta\nonumber\\
    &+ \frac{C_{L, \phi, \theta} C_{\ell, h, (\theta, x)} \decayfactor}{1 - \rho}\cdot \left(\sum_{\tau' = \tau}^{t-1} \norm{\zeta_{\tau'}} + \sum_{\tau' = 0}^{\tau-1} \rho^{\tau - \tau'}\norm{\zeta_{\tau'}}\right) \cdot \rho^{t-\tau}, \label{thm:OCO-with-parameter-update-gradient-bias:e2:s2}
\end{align}
\end{subequations}
where we use \Cref{coro:smooth-multi-step-costs} in \eqref{thm:OCO-with-parameter-update-gradient-bias:e2:s1} and \Cref{thm:OCO-with-parameter-update-trajectory-distance} in \eqref{thm:OCO-with-parameter-update-gradient-bias:e2:s2}.

For any $\tau' \in [\tau+1: t-1]$, since $x_{\tau'}, x_{\tau'+1} \in B_n(0, \bar{R}_C)$, we see that
\begin{subequations}\label{thm:OCO-with-parameter-update-gradient-bias:e3}
\begin{align}
    &\norm{\left.\frac{\partial h_{t\mid \tau}^*}{\partial \theta_\tau} \right|_{x_\tau, \theta_{\tau:\tau'}, (\theta_t)_{\times (t - \tau')}} - \left.\frac{\partial h_{t\mid \tau}^*}{\partial \theta_\tau} \right|_{x_\tau, \theta_{\tau:\tau'-1}, (\theta_t)_{\times (t - \tau' + 1)}}}\nonumber\\
    ={}& \norm{\left(\left.\frac{\partial h_{t\mid \tau'}^*}{\partial x_{\tau'}} \right|_{x_{\tau'}, \theta_{\tau'}, (\theta_t)_{\times (t - \tau')}} - \left.\frac{\partial h_{t\mid \tau'}^*}{\partial x_{\tau'}} \right|_{x_{\tau'}, (\theta_t)_{\times (t - \tau' + 1)}}\right)\left.\frac{\partial g_{\tau'\mid \tau}^*}{\partial \theta_\tau}\right|_{x_\tau, \theta_{\tau:\tau'-1}}}\nonumber\\
    \leq{}& \norm{\left.\frac{\partial h_{t\mid \tau'}^*}{\partial x_{\tau'}} \right|_{x_{\tau'}, \theta_{\tau'}, (\theta_t)_{\times (t - \tau')}} - \left.\frac{\partial h_{t\mid \tau'}^*}{\partial x_{\tau'}} \right|_{x_{\tau'}, (\theta_t)_{\times (t - \tau' + 1)}}} \cdot \norm{\left.\frac{\partial g_{\tau'\mid \tau}^*}{\partial \theta_\tau}\right|_{x_\tau, \theta_{\tau:\tau'-1}}}\nonumber\\
    \leq{}& C_{\ell, h, (x, \theta)} \decayfactor^{t - \tau'} \norm{\theta_t - \theta_{\tau'}} \cdot C_{L, \phi, \theta} \decayfactor^{\tau' - \tau} \label{thm:OCO-with-parameter-update-gradient-bias:e3:s1}\\
    \leq{}& {C_{\ell, h, (x, \theta)}C_{L, \phi, \theta}} \cdot \decayfactor^{t - \tau} \cdot \left(\frac{C_{L, h, \theta}}{1 - \decayfactor}\cdot (t - \tau')\eta + \sum_{\tau'' = \tau'}^{t-1} \norm{\zeta_{\tau''}}\right), \label{thm:OCO-with-parameter-update-gradient-bias:e3:s2}
\end{align}
\end{subequations}
where we use \Cref{lemma:smooth-multi-step-dynamics} and \Cref{coro:smooth-multi-step-costs} in \eqref{thm:OCO-with-parameter-update-gradient-bias:e3:s1}; we use \Cref{thm:OCO-with-parameter-update-trajectory-distance} in \eqref{thm:OCO-with-parameter-update-gradient-bias:e3:s2}. Substituting \eqref{thm:OCO-with-parameter-update-gradient-bias:e2} and \eqref{thm:OCO-with-parameter-update-gradient-bias:e3} into \eqref{thm:OCO-with-parameter-update-gradient-bias:e1} finishes the proof of the first inequality.

For the second inequality, recall that $G_t$ and $\nabla F_t(\theta_t)$ are given by
\[G_t \coloneqq \sum_{\tau = 0}^{t} \left.\frac{\partial h_{t\mid 0}^*}{\partial \theta_{t - \tau}}\right|_{x_0, \theta_{0:t}}, \nabla F_t(\theta_t) = \sum_{\tau=0}^t \left.\frac{\partial h_{t\mid 0}^*}{\partial \theta_{t - \tau}}\right|_{x_0, (\theta_t)_{\times (t+1)}}.\]
Therefore, we see that
\begin{subequations}\label{thm:OCO-with-parameter-update-gradient-bias:e4}
\begin{align}
    \norm{G_t - \nabla F_t(\theta_t)}
    ={}& \norm{\sum_{\tau = 0}^{t} \left.\frac{\partial h_{t\mid 0}^*}{\partial \theta_{t - \tau}}\right|_{x_0, \theta_{0:t}} - \sum_{\tau=0}^t \left.\frac{\partial h_{t\mid 0}^*}{\partial \theta_{t - \tau}}\right|_{x_0, (\theta_t)_{\times (t+1)}}}\nonumber\\
    \leq{}& \sum_{\tau = 0}^{t} \norm{\left.\frac{\partial h_{t\mid 0}^*}{\partial \theta_{t - \tau}}\right|_{x_0, \theta_{0:t}} - \left.\frac{\partial h_{t\mid 0}^*}{\partial \theta_{t - \tau}}\right|_{x_0, (\theta_t)_{\times (t+1)}}}\label{thm:OCO-with-parameter-update-gradient-bias:e4:s1}\\
    \leq{}& \sum_{\tau = 0}^{t} \left(\hat{C}_0 + \hat{C}_1 \tau + \hat{C}_2 \tau^2\right) \decayfactor^\tau \eta\\
    &+ \sum_{\tau = 0}^{t-1} \left(\hat{C}_3 \sum_{\tau' = \tau}^{t-1} \norm{\zeta_{\tau'}} + \hat{C}_4 \sum_{\tau' = \tau}^{t-1} (\tau'-\tau) \norm{\zeta_{\tau'}}\right) \cdot \rho^{t-\tau}\nonumber\\
    &+ \hat{C}_5 \sum_{\tau = 0}^{t-1} \sum_{\tau' = 0}^{\tau-1} \rho^{t-\tau'} \norm{\zeta_{\tau'}}\label{thm:OCO-with-parameter-update-gradient-bias:e4:s2}\\
    \leq{}& \left(\frac{\hat{C}_0}{1 - \decayfactor} + \frac{\hat{C}_1 + \hat{C}_2}{(1 - \decayfactor)^2} + \frac{\hat{C}_2}{(1 - \decayfactor)^3}\right) \eta \nonumber\\
    &+ \sum_{\tau = 0}^{t-1}\left(\frac{\hat{C}_3}{1 - \rho}
    + \frac{\hat{C}_4}{(1 - \rho)^2} + \hat{C}_5 (t - \tau)\right) \rho^{t-\tau} \norm{\zeta_\tau},\nonumber
\end{align}
\end{subequations}
where we use the triangle inequality in \eqref{thm:OCO-with-parameter-update-gradient-bias:e4:s1}; we use the first inequality in \Cref{thm:OCO-with-parameter-update-gradient-bias:full} that we have shown and \Cref{coro:smooth-multi-step-costs} in \eqref{thm:OCO-with-parameter-update-gradient-bias:e4:s2}.

%%%%%%%%%%%%%%%%%%%%%%%%%%%%%%%%%%%%%%%%%%%%%%%%%%%%%%%%%%%%%%%%%%%%%%%%%%%%%%%
%%%%%%%%%%%%%%%%%%%%%%%%%%%%%%%%%%%%%%%%%%%%%%%%%%%%%%%%%%%%%%%%%%%%%%%%%%%%%%%

\end{document}